\documentclass[11pt,twoside]{article}

\usepackage[utf8]{inputenc}
\usepackage{epsfig,amsfonts,color}
\usepackage{amsmath, bbm, mathabx}
\usepackage[normalem]{ulem}
\usepackage{amssymb, palatino, geometry,url}
\usepackage[colorlinks=true,linkcolor=blue,citecolor=blue,urlcolor=blue]{hyperref}
\usepackage[font=small,labelfont=bf]{caption}
\usepackage{subcaption}
\usepackage{multirow}
\usepackage{tikz}
\usetikzlibrary{positioning, arrows.meta, calc, fit, backgrounds, shapes.geometric}
\usepackage{pgfplots}
\pgfplotsset{compat=1.18}
\usepgfplotslibrary{fillbetween}
\tikzset{
  flowbox/.style={draw, rounded corners=3pt, align=center,
    minimum height=1.6cm, text width=3.2cm,
    inner sep=7pt, font=\small},
  flowdec/.style={draw, diamond, aspect=1.6, align=center,
    inner sep=2pt, font=\small, text width=2.0cm},
  flowarr/.style={-Latex, thick}
}

\usepackage{fancyhdr}
\usepackage{lineno}
\usepackage{float}
\usepackage[section]{placeins}
\usepackage[table]{xcolor}

\usepackage[dvipsnames]{xcolor}
\usepackage{accents}
\usepackage{listings}
\DeclareCaptionFont{white}{\color{white}}
\DeclareCaptionFormat{listing}{{\parbox{\dimexpr\textwidth-2\fboxsep\relax}{#1#2#3}}}
\lstdefinelanguage{julia}
{
keywordsprefix=\@,
morekeywords={
exit,whos,edit,load,is,isa,isequal,typeof,tuple,ntuple,uid,hash,finalizer,convert,promote,
subtype,typemin,typemax,realmin,realmax,sizeof,eps,promote_type,method_exists,applicable,
invoke,dlopen,dlsym,system,error,throw,assert,new,Inf,Nan,pi,im,begin,while,for,in,return,
break,continue,macro,quote,let,if,elseif,else,try,catch,end,bitstype,ccall,do,using,module,
import,export,importall,baremodule,immutable,local,global,const,Bool,Int,Int8,Int16,Int32,
Int64,Uint,Uint8,Uint16,Uint32,Uint64,Float32,Float64,Complex64,Complex128,Any,Nothing,None,
function,type,typealias,abstract,get_node,create_estimation_model,set_solution, Matrix, Dict,
solve, get_solution, solve_ss_problem, create_estimation_problem, addnode, nothing, 
set_optimizer_attribute, optimize!, println, value, termination_status, objective_value, Model,
optimizer_with_attributes, zeros, length, getJuMPmodel, savefig, plot, plot!, scatter!, grid, 
xlabel!, ylabel!, xlabel, ylabel, zip, sum, InfiniteModel, seed!, Normal, mean, var, Infinite,
OrthogonalCollocation, collect, OptiGraph, set_optimizer, eachindex, set_parameter_value, 
Parameter, dual, set_silent, InfiniteLogical, Disjunct, Hull, InfiniteGDPModel, ExaTranscriptionBackend, CUDABackend, simulate, allowscalar, warmstart_backend_start_values, sin,
MBM, Exactly, integral, 
},
morekeywords = [2]{},
alsoletter=!,
sensitive=true,
morecomment=[l]{\#},
morestring=[b]"
}

\newcommand{\threedisks}{%
  \addplot[fill=gray!55, draw=black!75,
           line width=0.7pt, smooth, samples=80, domain=0:360,
           forget plot]
      ({1+cos(x)}, {1+sin(x)});
  \addplot[fill=gray!55, draw=black!75,
           line width=0.7pt, smooth, samples=80, domain=0:360,
           forget plot]
      ({4+cos(x)}, {4.5+sin(x)});
  \addplot[fill=gray!55, draw=black!75,
           line width=0.7pt, smooth, samples=80, domain=0:360,
           forget plot]
      ({7+cos(x)}, {1+sin(x)});
  \node[font=\footnotesize] at (axis cs:1,1) {$W_1$};
  \node[font=\footnotesize] at (axis cs:4,4.5) {$W_2$};
  \node[font=\footnotesize] at (axis cs:7,1) {$W_3$};
}

\usepackage{booktabs}
\pgfplotsset{
  relaxpanel/.style={
    width=\linewidth, height=\linewidth,
    axis equal image,
    xmin=-1, xmax=9, ymin=-1.5, ymax=6.5,
    xlabel=$x$, ylabel=$y$,
    grid=major, grid style={gray!20},
    legend style={at={(0.97,0.97)}, anchor=north east,
                  font=\scriptsize, draw=black!50,
                  fill opacity=0.9},
    title style={font=\small},
    label style={font=\footnotesize},
    tick label style={font=\footnotesize},
  },
  solpanel/.style={
    width=\linewidth, height=0.62\linewidth,
    xlabel={MIP Solve Time [s]},
    ylabel={Solution Progress (5\% MIP-Gap $= 0.95$)},
    ymin=0, ymax=1.03,
    ytick={0,0.2,0.4,0.6,0.8,0.95,1.0},
    grid=major, grid style={gray!18},
    legend style={at={(0.97,0.03)}, anchor=south east,
                  font=\scriptsize, draw=black!50,
                  fill opacity=0.9},
    label style={font=\footnotesize},
    tick label style={font=\footnotesize},
    every axis plot/.append style={line width=1pt},
  },
  timebar/.style={
    width=\linewidth, height=0.55\linewidth,
    ybar stacked, bar width=22pt,
    symbolic x coords={BigM,MBM,Hull,PSplit-3},
    xtick=data,
    ymin=0,
    ylabel={Total Time [s]},
    grid=major, grid style={gray!18}, ymajorgrids=true,
    legend style={at={(0.03,0.97)}, anchor=north west,
                  font=\scriptsize, draw=black!50},
    label style={font=\footnotesize},
    tick label style={font=\footnotesize},
  },
  trajpanel/.style={
    width=\linewidth, height=0.72\linewidth,
    grid=major, grid style={gray!18},
    label style={font=\footnotesize},
    tick label style={font=\footnotesize},
    legend style={font=\scriptsize, draw=black!50,
                  fill opacity=0.9},
    every axis plot/.append style={line width=0.9pt},
  }
}

\title{Solution Methods for Infinite-Dimensional \\
Generalized Disjunctive Programming}
\author{Daniel Nguyen$^1$ and
    Joshua L. Pulsipher$^1$\thanks{Corresponding Author:
    pulsipher@uwaterloo.ca}\\
    {\small $^1$Department of Chemical Engineering}\\
	{\small \;University of Waterloo, 200 University Ave W, Waterloo, ON N2L 3G1, Canada}}
\date{}

\begin{document}

\raggedbottom

\emergencystretch=3em

\maketitle

\begin{abstract}
Generalized disjunctive programming (GDP) expresses mixed
discrete-continuous decisions through Boolean indicators and
disjunctions, and can be systematically solved via a library of methods proposed in the literature. The recent InfiniteGDP abstraction
lifts this modeling layer to continuous domains such as time, space, and
uncertainty, but only the big-M and hull reformulations, the
two endpoints of the relaxation spectrum, have been generalized to the
infinite setting. This work closes this gap by generalizing four other GDP
solution methods to infinite-dimensional optimization: the multiple
big-M reformulation, P-split reformulation, cutting plane reformulation, and the logic-based
outer approximation algorithm. It further proposes MBM-GP, a novel Gaussian-process variant of multiple big-M that learns the big-M function over the infinite domain from a small subset of the subproblem solves. Moreover, these approaches are implemented in the Julia package \texttt{InfiniteDisjunctiveProgramming.jl}. The methods are benchmarked on case studies arising in dynamic and stochastic optimization. The results demonstrate how the generalized solution methods can outperform big-M and hull, with MBM-GP retaining the tightness of multiple big-M at a fraction of its reformulation cost.
\end{abstract}


\section{Introduction}\label{sec:intro}

Infinite-dimensional optimization (InfiniteOpt) problems involve
decision variables, objectives, and constraints indexed over
continuous domains such as time, space, and random scenarios
\cite{pulsipher2022infiniteopt}. InfiniteOpt problems are prominent
across process systems engineering (PSE), arising in dynamic
optimization, model predictive control, parameter estimation for
(partial)-differential equations, and optimization
under uncertainty. They are equally prominent outside of
PSE in domains such as vehicle and trajectory optimization
\cite{salazar2017hybrid, heilmeier2019qss} and energy and process
scheduling \cite{biodiesel}. Pulsipher et al.
\cite{pulsipher2022infiniteopt} propose a unifying modeling
abstraction that organizes these problems around a set of
generalized objects (infinite parameters, measure operators, and
differential operators) and pair it with a suite of direct
transcription methods (e.g., orthogonal
collocation over finite elements, finite differences, Monte Carlo
sampling) that automate the conversion of an InfiniteOpt model into
a discretized form suitable for conventional mathematical programming solvers
\cite{pulsipher2022infiniteopt, biegler2007overview}. The
abstraction relieves the modeler of the previous burden of
hand-discretizing continuous dynamics into a finite-dimensional formulation, but
offers no analogous relief for \emph{discrete} decisions whose
value may vary across the infinite domain. Such decisions must
still be encoded directly as binary variables and algebraic
relationships, leaving the modeler to construct the resulting
mixed-integer program (MIP) themselves.

These mixed discrete-continuous decisions are the natural setting for
generalized disjunctive programming (GDP), first proposed by Raman
and Grossmann \cite{raman1994modelling}. GDP models
discrete structure through Boolean indicator variables tied to
disjuncts, each of which carries its own constraint set, with the
active disjunct selected by logical propositions over the indicators
\cite{grossmann2013systematic}. The framework retains the high-level logical
structure of a model while delegating its conversion into a
MIP to a library of automatic reformulations
\cite{trespalacios2014review, kronqvist2025fifty, lee2001global,
papageorgiou2018pseudo}. The classical big-M and hull
reformulations \cite{grossmann2013systematic, balas1985disjunctive} have
been refined and extended in different ways, including per-pair
tightening of big-M (MBM) \cite{trespalacios2015mbm}, partition-based
interpolation between big-M and hull (P-split)
\cite{kronqvist2022psplit}, and extended cutting planes (CP)
\cite{trespalacios2016cutting}. Alongside these reformulations,
logic-based outer approximation (LOA) \cite{turkay1996logic,
viswanathan1990penalty, bergamini2005logic} takes an algorithmic
route for convex GDPs with nonlinear disjunct constraints: rather
than converting the GDP into a single MIP, it
solves the GDP iteratively through fixed-indicator subproblems and
an accumulating master. Software implementations of GDP and its solution methods include \texttt{pyomo.GDP} \cite{chen2022pyomo}, \texttt{DisjunctiveProgramming.jl }\cite{perez2023disjunctiveprogramming}, and GAMS. 

Despite its maturity in the finite case, the application of GDP to
InfiniteOpt problems has so far seen limited development. Related
work addresses the same class of mixed discrete-continuous
decisions over infinite domains, but expresses the model directly as a
mixed-integer nonlinear program (MINLP) rather than as a GDP.
Mixed-integer dynamic optimization (MIDO) \cite{bansal2003mido}
formulates the problem as an MINLP over differential-algebraic
equations with binary variables, and mixed-integer optimal control
(MIOC) \cite{sager2009direct, hante2013relaxation, sager2012integer, zeile2021dwell}.
Both predate the modern GDP machinery and bypass the disjunctive
modeling layer, so problems whose discrete structure is naturally
expressed through Boolean indicators and logical propositions must
be encoded by hand into binary variables and algebraic
relationships. 

Gondosiswanto and Pulsipher
\cite{gondosiswanto2025infiniteopt} recently bridge this gap on the
modeling side by proposing the InfiniteGDP abstraction. It lifts
disjunctions, logical propositions, and cardinality constraints into the InfiniteOpt setting such that GDP can be used to more intuitively model InfiniteOpt problems (e.g., MPC with switching constraints). Moreover, \cite{gondosiswanto2025infiniteopt} also generalizes 
the big-M, hull, and
indicator reformulations to solve InfiniteGDP problems. Big-M and hull represent two extremes: big-M avoids adding an excessive number of variables/constraints at the cost of weak continuous relaxation, while hull yields the tightest relaxation at the cost of increased problem size. The intermediate reformulations that have driven
much of the empirical progress in finite GDP, including MBM,
P-split, and CP, fall between these two extremes and offer
reasonable tradeoffs, while LOA sidesteps the spectrum altogether by
solving the GDP iteratively. To date, none of these other approaches have been generalized to
InfiniteGDP.

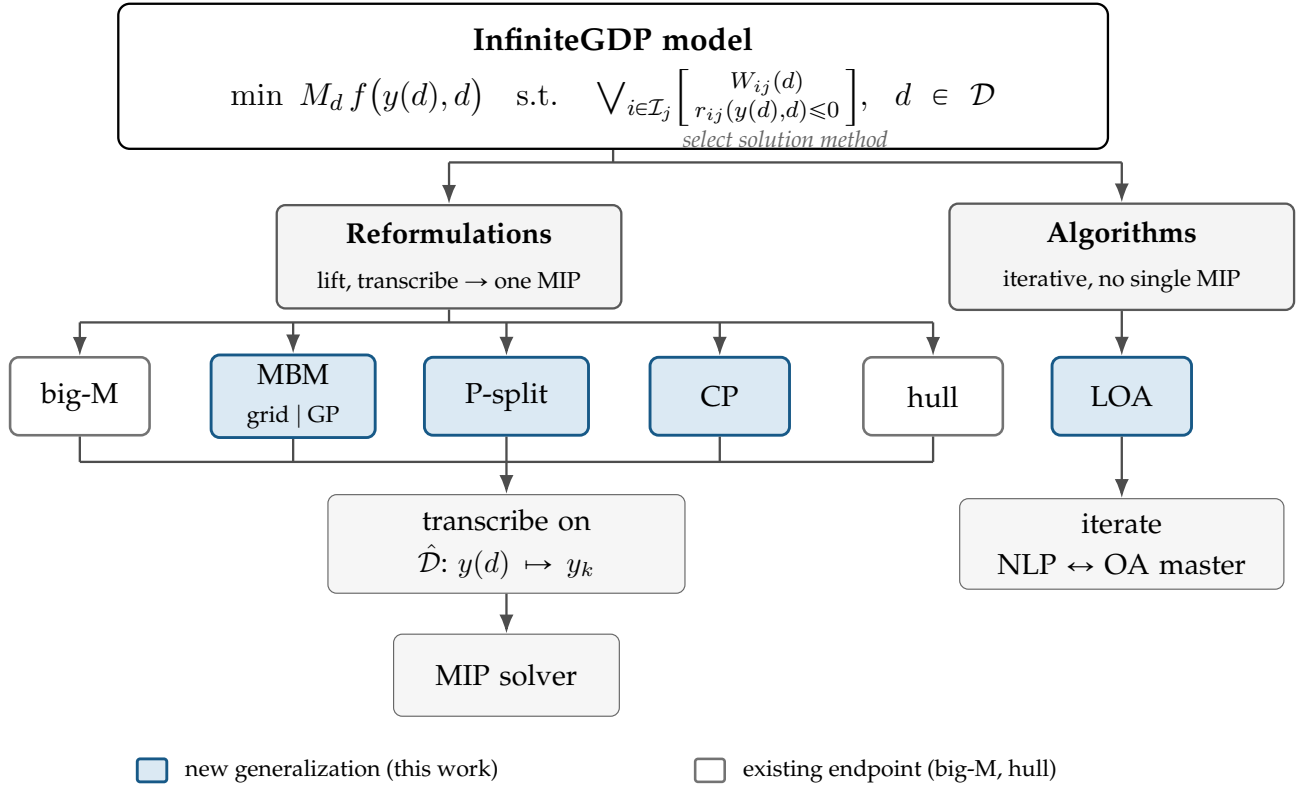
\begin{figure}[H]
\centering
\resizebox{\textwidth}{!}{%
\definecolor{gaaccent}{RGB}{31,119,180}%
\begin{tikzpicture}[
    >=Latex,
    rootbox/.style={draw, rounded corners=4pt, fill=white,
      line width=0.8pt, align=center, inner sep=8pt},
    catbox/.style={draw=black!65, rounded corners=3pt, fill=gray!9,
      line width=0.8pt, align=center, font=\small, inner sep=6pt,
      minimum width=4.2cm},
    methnew/.style={draw=gaaccent!75!black, line width=1pt,
      fill=gaaccent!16, rounded corners=3pt, minimum height=0.95cm,
      minimum width=1.7cm, font=\small, align=center, inner sep=3pt},
    methend/.style={draw=black!55, line width=0.9pt, fill=white,
      rounded corners=3pt, minimum height=0.95cm, minimum width=1.7cm,
      font=\small, align=center, inner sep=3pt},
    procbox/.style={draw=black!60, rounded corners=3pt, fill=gray!7,
      align=center, font=\small, inner sep=5pt, minimum height=1.0cm},
    flow/.style={-Latex, thick, black!70},
    tline/.style={thick, black!70}
  ]
  \node[rootbox, text width=11.5cm] (model) at (8.0, 12.6)
    {\textbf{InfiniteGDP model}\\[3pt]
     $\min\ M_d\, f\bigl(y(d),d\bigr)
       \quad \mathrm{s.t.}\quad
       \bigvee_{i \in \mathcal{I}_j}\!
       \left[\begin{smallmatrix} W_{ij}(d) \\
         r_{ij}(y(d),d) \le 0 \end{smallmatrix}\right]\!,\ \ d \in \mathcal{D}$};

  \node[catbox] (refc) at (6.0, 10.4)
    {\textbf{Reformulations}\\[1pt]
     {\scriptsize lift, transcribe $\to$ one MIP}};
  \node[catbox] (algc) at (14.2, 10.4)
    {\textbf{Algorithms}\\[1pt]
     {\scriptsize iterative, no single MIP}};

  \node[methend] (bigm)   at (1.5, 8.7) {big-M};
  \node[methnew, minimum width=2.0cm] (mbm) at (4.1, 8.7)
    {MBM\\[-1pt]{\scriptsize grid $\mid$ GP}};
  \node[methnew, minimum width=2.0cm] (psplit) at (6.7, 8.7) {P-split};
  \node[methnew] (cp)     at (9.3, 8.7) {CP};
  \node[methend] (hull)   at (11.9, 8.7) {hull};

  \node[methnew] (loa) at (14.2, 8.7) {LOA};

  \node[procbox, text width=4.0cm] (trans) at (6.7, 6.9)
    {transcribe on $\hat{\mathcal{D}}$:\ $y(d)\!\mapsto\! y_k$};
  \node[procbox, text width=2.4cm] (solver) at (6.7, 5.3)
    {MIP solver};
  \node[procbox, text width=3.6cm] (loaloop) at (14.2, 6.9)
    {iterate\\ NLP $\leftrightarrow$ OA master};

  \draw[tline] (model.south) -- (8.0, 11.55);
  \draw[tline] (6.0, 11.55) -- (14.2, 11.55);
  \draw[flow] (6.0, 11.55) -- (refc.north);
  \draw[flow] (14.2, 11.55) -- (algc.north);
  \node[anchor=south, font=\scriptsize\itshape, text=black!60]
    at (10.1, 11.6) {select solution method};

  \draw[tline] (refc.south) -- (6.0, 9.6);
  \draw[tline] (1.5, 9.6) -- (11.9, 9.6);
  \draw[flow] (1.5, 9.6) -- (bigm.north);
  \draw[flow] (4.1, 9.6) -- (mbm.north);
  \draw[flow] (6.7, 9.6) -- (psplit.north);
  \draw[flow] (9.3, 9.6) -- (cp.north);
  \draw[flow] (11.9, 9.6) -- (hull.north);

  \draw[tline] (bigm.south)   -- (1.5, 7.9);
  \draw[tline] (mbm.south)    -- (4.1, 7.9);
  \draw[tline] (psplit.south) -- (6.7, 7.9);
  \draw[tline] (cp.south)     -- (9.3, 7.9);
  \draw[tline] (hull.south)   -- (11.9, 7.9);
  \draw[tline] (1.5, 7.9) -- (11.9, 7.9);
  \draw[flow] (6.7, 7.9) -- (trans.north);
  \draw[flow] (trans.south) -- (solver.north);

  \draw[flow] (algc.south) -- (loa.north);
  \draw[flow] (loa.south) -- (loaloop.north);

  \filldraw[draw=gaaccent!75!black, line width=1pt, fill=gaaccent!16,
            rounded corners=1pt] (2.2,4.0) rectangle (2.56,4.28);
  \node[anchor=west, font=\scriptsize] at (2.65,4.14)
    {new generalization (this work)};
  \draw[draw=black!55, line width=0.9pt, fill=white, rounded corners=1pt]
    (9.0,4.0) rectangle (9.36,4.28);
  \node[anchor=west, font=\scriptsize] at (9.45,4.14)
    {existing endpoint (big-M, hull)};
\end{tikzpicture}%
}
\caption{Graphical abstract. From a single InfiniteGDP model, the six
solution methods form two branches. The reformulation branch (big-M,
MBM, P-split, CP, hull, ordered along the tightness/cost spectrum)
lifts and transcribes the model into one MIP, with the MBM leaf
covering both its grid and Gaussian-process (MBM-GP)
implementations; the algorithm branch is occupied by LOA alone, which
branches separately and solves the model through an iterative procedure. Blue boxes mark the
generalizations introduced in this work; big-M and hull are the pre-existing endpoints.}
\label{fig:graphical_abstract}
\end{figure}

Hence, we propose infinite-dimensional generalizations of MBM, P-split, CP, and LOA to address the aforementioned gap. We additionally propose MBM-GP, a Gaussian-process variant of MBM that is native to the infinite setting: it learns the big-M function over the infinite domain from a small subset of the per-support subproblems, retaining MBM's tightness at a fraction of its reformulation cost.
All four methods expand the library of methods to solve InfiniteGDP problems.
Figure~\ref{fig:graphical_abstract} provides a summary of the proposed InfiniteGDP workflow, highlighting the main contributions of this work. A single InfiniteGDP model branches into six
interchangeable solution methods that split into two families. The
first is the reformulation branch of five methods that lift and
transcribe the model into a single MIP, where MBM admits a choice
between its grid and Gaussian-process (MBM-GP) implementations. The second is the algorithm branch, which
now includes LOA that directly solves the InfiniteGDP without reformulation to a single MIP.

The rest of this paper is structured as follows.
Section~\ref{sec:background} reviews the GDP and InfiniteGDP modeling
abstractions and their existing reformulations.
Section~\ref{sec:reformulations} presents the four generalized solution
methods and the novel MBM-GP variant in detail, together with their
implementation and practical
guidance on choosing among them. Section~\ref{sec:casestudies}
performs benchmarks on three case studies. Finally, Section~\ref{sec:conclusions} provides
concluding remarks and recommendations for future developments.

\section{Background and Notation}\label{sec:background}

In this section, we review relevant GDP background and establish notation. For a more thorough discussion of these topics, please refer to \cite{grossmann2021advanced, kronqvist2022psplit, gondosiswanto2025infiniteopt}.

\subsection{Generalized Disjunctive Programming}\label{sec:bg_gdp}

GDP encodes discrete structure on a finite
decision vector $z \in \mathcal{Z} \subseteq \mathbb{R}^{n_z}$
through Boolean indicator variables and disjunctions \cite{raman1994modelling, balas1985disjunctive,
grossmann2013systematic}. The
generic GDP problem is
\begin{equation}\label{eq:bg_disj}
\begin{aligned}
    \min_{z,\, W} \quad & f(z) \\
    \text{s.t.} \quad & g(z) \leq 0 \\
    & \bigvee_{i \in \mathcal{I}_j}
      \begin{bmatrix}
        W_{ij} \\[2pt]
        r_{ij}(z) \leq 0
      \end{bmatrix},
      \quad j \in \mathcal{J} \\
    & \Omega(W) = \text{True} \\
    & z \in \mathcal{Z},\;
      W_{ij} \in \{\text{True}, \text{False}\}
\end{aligned}
\end{equation}
with objective $f$, global constraints $g$, disjunctions
indexed by $j \in \mathcal{J}$ with disjuncts $i \in \mathcal{I}_j$
carrying Boolean indicators $W_{ij}$ and disjunct constraints
$r_{ij}(z) \leq 0$, and a logical proposition $\Omega(W)$ over
the indicators (built from $\vee$, $\wedge$, $\neg$,
$\Longrightarrow$, $\Leftrightarrow$ and the cardinality sets
\textsc{Exactly}, \textsc{AtLeast}, \textsc{AtMost}). On the
binary realization $w_{ij} \in \{0, 1\}$ the proposition
reduces to algebraic relations \cite{grossmann2013systematic}.

A reformulation maps Problem \eqref{eq:bg_disj} to a finite
MIP; a solution algorithm such as LOA instead
solves Problem \eqref{eq:bg_disj} through a sequence of subproblems without
ever forming a single equivalent program. The classical endpoints
(big-M, hull), intermediate refinements (P-split, MBM), and
decomposition algorithms (CP, LOA) are outlined below.

\paragraph{Big-M.} The simplest reformulation
\cite{grossmann2013systematic} produces the MIP
\begin{equation}\label{eq:bg_bigm}
\begin{aligned}
    \min_{z,\, w} \quad & f(z) \\
    \text{s.t.} \quad & g(z) \leq 0 \\
    & r_{ij}(z) \leq M_{ij}\, (1 - w_{ij}),
      \quad i \in \mathcal{I}_j,\; j \in \mathcal{J} \\
    & \sum_{i \in \mathcal{I}_j} w_{ij} = 1,
      \quad j \in \mathcal{J} \\
    & z \in \mathcal{Z},\; w_{ij} \in \{0, 1\}
\end{aligned}
\end{equation}
where $M_{ij} \geq \max_{z \in \mathcal{Z}} r_{ij}(z)$ is large enough
that the disjunct constraint is slack when $w_{ij} = 0$. Big-M is
structure-preserving and stays linear when $r_{ij}$ is, but yields
a weak continuous relaxation when the $M_{ij}$ are loose.

\paragraph{Hull.} The hull reformulation
\cite{grossmann2013systematic, balas1985disjunctive} disaggregates the
decision $z$ into per-disjunct copies $\nu_{ij}$ and produces
the mixed-integer program
\begin{equation}\label{eq:bg_hull}
\begin{aligned}
    \min_{z,\, \nu_{ij},\, w} \quad & f(z) \\
    \text{s.t.} \quad & g(z) \leq 0 \\
    & z = \sum_{i \in \mathcal{I}_j} \nu_{ij},
      \quad j \in \mathcal{J} \\
    & w_{ij}\, r_{ij}\!\bigl(\nu_{ij}/w_{ij}\bigr) \leq 0,
      \quad i \in \mathcal{I}_j,\; j \in \mathcal{J} \\
    & w_{ij}\, z^L \leq \nu_{ij} \leq w_{ij}\, z^U,
      \quad i \in \mathcal{I}_j,\; j \in \mathcal{J} \\
    & \sum_{i \in \mathcal{I}_j} w_{ij} = 1,
      \quad j \in \mathcal{J} \\
    & z \in \mathcal{Z},\; w_{ij} \in \{0, 1\}
\end{aligned}
\end{equation}
The perspective transform $w_{ij}\, r_{ij}(\nu_{ij}/w_{ij})$ keeps the
disjunct constraint convex when $r_{ij}$ is. Hull yields the
tightest convex relaxation per disjunction at the cost of
duplicated variables and typically nonlinear (perspective)
constraints.

\paragraph{P-split.} P-split \cite{kronqvist2022psplit} partitions
the variables of each disjunct constraint into $P \in \mathbb{Z}_+$
groups $\{V_1,\ldots,V_P\}$, writing
$r_{ij}(z) = \sum_{p=1}^P \tilde{r}_p(z_{V_p}) + c \leq 0$ with
constant $c \in \mathbb{R}$, and introduces auxiliary variables
$v_{ij,p} \in \mathbb{R}$ that envelope each partition expression by
splitting the constraint into
\begin{align}
    \tilde{r}_p(z_{V_p}) &\leq v_{ij,p}, \quad p = 1,\ldots,P,\;
      i \in \mathcal{I}_j,\; j \in \mathcal{J}, \label{eq:psplit_part}\\
    \sum_{p=1}^{P} v_{ij,p} + c &\leq 0,
      \quad i \in \mathcal{I}_j,\; j \in \mathcal{J}. \label{eq:psplit_sum}
\end{align}
The envelope constraints \eqref{eq:psplit_part} are enforced
globally after each $v_{ij,p}$ is bounded via interval arithmetic
on $z_{V_p}^L, z_{V_p}^U$, while hull disaggregation is applied
only to the linking constraint \eqref{eq:psplit_sum}, with $z$
itself left undisaggregated.
$P = 1$ recovers big-M, $P = n_z$ recovers hull, and
intermediate choices of $P$ trade model size for relaxation tightness.
This partitioned form presumes the disjunct constraint is
additively separable across the chosen groups, i.e., expressible
as a sum $\sum_p \tilde{r}_p(z_{V_p})$; non-separable constraints can
still be accommodated by introducing enveloping auxiliary
functions for the coupled terms or by selecting a partition that
keeps interacting variables within the same group.

\paragraph{Multiple Big-M.} MBM \cite{trespalacios2015mbm}
tightens big-M by replacing each $M_{ij}$ with per-pair constants
$M_{ii'}$, one for each pair of disjuncts
$i, i' \in \mathcal{I}_j$ of disjunction $j$, so that each disjunct
constraint $r_{ij}(z) \leq 0$ is replaced with
\begin{equation}\label{eq:mbm_finite}
    r_{ij}(z) \leq \sum_{i' \neq i} M_{ii'}\, w_{i'j},
    \quad i \in \mathcal{I}_j,\; j \in \mathcal{J}.
\end{equation}
Each $M_{ii'}$
comes from a small subproblem maximizing $r_{ij}$ subject only to
disjunct $i'$'s constraint and the variable bounds:
\begin{equation}\label{eq:bg_mbm}
\begin{aligned}
    M_{ii'} = \max_{z} \quad & r_{ij}(z) \\
    \text{s.t.} \quad & r_{i'j}(z) \leq 0 \\
    & z^L \leq z \leq z^U.
\end{aligned}
\end{equation}
The reformulated constraint \eqref{eq:mbm_finite} activates only
the $M$ values relevant to the chosen disjunct, giving a relaxation
strictly between big-M and hull while keeping the model linear
when $r_{ij}, r_{i'j}$ are.

\paragraph{Cutting planes.} The CP algorithm
\cite{trespalacios2016cutting} alternates between a relaxed big-M
master problem on the original variables and a hull-based
separation subproblem that projects the master solution onto the
hull relaxation. At each iteration, let
$\boldsymbol{z}^{\mathrm{rBM}}$ denote the optimal solution of the
current master. The separation subproblem disaggregates
$\boldsymbol{z}$ into per-disjunct copies $\boldsymbol{\nu}_{ij}$
for $i \in \mathcal{I}_j$, $j \in \mathcal{J}$ with relaxed
indicators $w_{ij} \in [0, 1]$ and solves
\begin{equation}\label{eq:cp_sep}
\begin{aligned}
    \min_{\boldsymbol{z},\,\boldsymbol{\nu}_{ij},\,w_{ij}}
    \quad & \|\boldsymbol{z} - \boldsymbol{z}^{\mathrm{rBM}}\|_2^2 \\
    \text{s.t.} \quad
    & \boldsymbol{z} = \sum_{i \in \mathcal{I}_j} \boldsymbol{\nu}_{ij},
      && j \in \mathcal{J} \\
    & w_{ij}\, r_{ij}\bigl(\boldsymbol{\nu}_{ij}/w_{ij}\bigr) \leq 0,
      && i \in \mathcal{I}_j,\; j \in \mathcal{J} \\
    & w_{ij}\, z^L \leq \boldsymbol{\nu}_{ij} \leq w_{ij}\, z^U,
      && i \in \mathcal{I}_j,\; j \in \mathcal{J} \\
    & \sum_{i \in \mathcal{I}_j} w_{ij} = 1,
      \quad w_{ij} \in [0, 1],
      && i \in \mathcal{I}_j,\; j \in \mathcal{J}.
\end{aligned}
\end{equation}
Let $\boldsymbol{z}^{\mathrm{SEP}}$ denote the optimal $\boldsymbol{z}$.
If $\|\boldsymbol{z}^{\mathrm{SEP}} - \boldsymbol{z}^{\mathrm{rBM}}\|_2^2
> \varepsilon$, the separating cut
\begin{equation}\label{eq:cp_cut}
    2\bigl(\boldsymbol{z}^{\mathrm{SEP}}
        - \boldsymbol{z}^{\mathrm{rBM}}\bigr)
    \bigl(\boldsymbol{z} - \boldsymbol{z}^{\mathrm{SEP}}\bigr)
    \geq 0
\end{equation}
is added to the master and the loop repeats. The algorithm
terminates when the projection distance falls below $\varepsilon$
or the master is proven optimal. CP avoids the variable doubling of
hull while recovering hull-like tightness through cuts on the
original variables. The construction assumes convexity: projecting
the master optimum onto the hull and obtaining a valid separating
hyperplane require the disjunct constraints, and hence the
separation subproblem, to be convex, so CP is sound only for
convex GDP. Even if a solution to the separation problem is
found, there is no guarantee that the resulting cut supports the
feasible region.

\paragraph{Logic-based outer approximation.} The logic-based outer
approximation (LOA) algorithm of \cite{turkay1996logic} solves
convex GDPs by alternating between a fixed nonlinear program (NLP)
subproblem~$(S')$ and an outer-approximation (OA) master
problem~$(\hat{M}^b_{OA})$. At each major iteration~$l$ the NLP
fixes all logical variables $W_{ij}$ to a single Boolean
assignment (a combination) and solves for the continuous
variables~$z$. The master accumulates first-order linearizations
of the objective~$f$ and the nonlinear disjunct
constraints~$h_{ij}$ at each NLP solution point~$y^l$, then
selects the next combination. The NLP objective provides an upper
bound~$Z_U$ and the master objective provides a lower
bound~$Z_L$, which converge to the optimum within finitely many
iterations. The algorithm terminates when
$|Z_U - Z_L| \leq \varepsilon$. The augmented penalty variant
\cite{viswanathan1990penalty} relaxes the OA cuts of nonlinear
disjunct constraints with bounded slacks, keeping the master
feasible after NLP infeasibility. The implementation here uses this
augmented-penalty form; because the bounded slacks absorb small
outer-approximation violations, the method tolerates mild
nonconvexity in the disjunct constraints in practice, though it
carries no global guarantee.

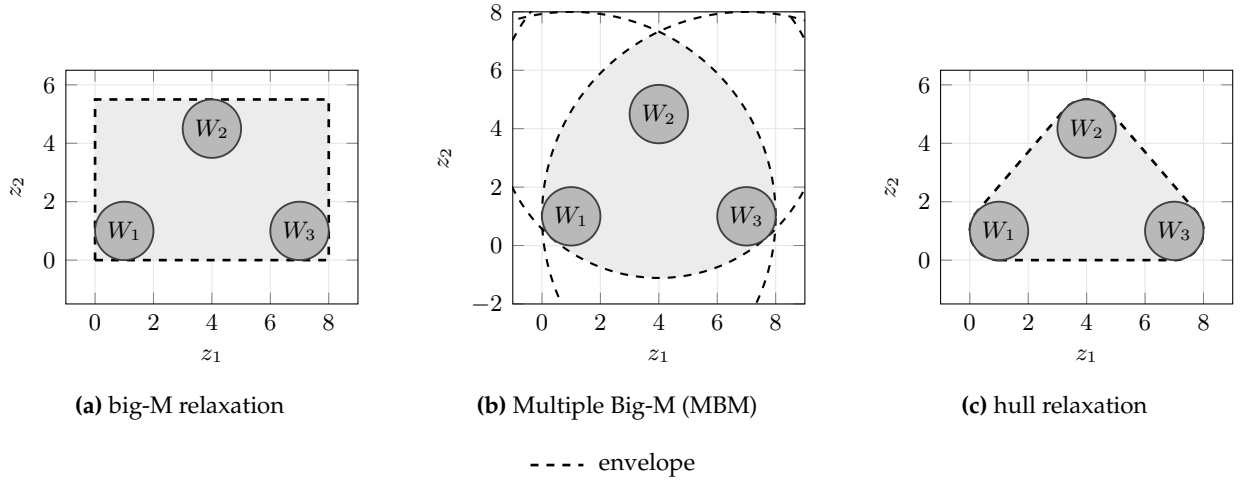
\begin{figure}[!htb]
\centering
\begin{subfigure}[t]{0.32\textwidth}
\centering
\begin{tikzpicture}
\begin{axis}[relaxpanel, xlabel = $z_1$, ylabel = $z_2$]
  \addplot[fill=gray!15, draw=black, dashed, line width=1pt]
      coordinates {(0,0) (8,0) (8,5.5) (0,5.5) (0,0)};
  \threedisks
\end{axis}
\end{tikzpicture}
\caption{big-M relaxation}
\end{subfigure}\hfill
\begin{subfigure}[t]{0.32\textwidth}
\centering
\begin{tikzpicture}
\begin{axis}[relaxpanel, ymin=-2, ymax=8, xlabel = $z_1$, ylabel = $z_2$]
  \draw[fill=gray!15, draw=none]
      (axis cs:0.014, 0.552)
      arc[start angle=224.7, end angle=315.3,
          x radius=5.61, y radius=5.61]
      arc[start angle=356.3, end angle=424.6,
          x radius=7, y radius=7]
      arc[start angle=115.4, end angle=183.7,
          x radius=7, y radius=7]
      -- cycle;
  \addplot[draw=black, dashed, line width=0.8pt,
           smooth, samples=120, domain=0:360, forget plot]
      ({1+7*cos(x)}, {1+7*sin(x)});
  \addplot[draw=black, dashed, line width=0.8pt,
           smooth, samples=120, domain=0:360, forget plot]
      ({4+5.61*cos(x)}, {4.5+5.61*sin(x)});
  \addplot[draw=black, dashed, line width=0.8pt,
           smooth, samples=120, domain=0:360, forget plot]
      ({7+7*cos(x)}, {1+7*sin(x)});
  \threedisks
\end{axis}
\end{tikzpicture}
\caption{Multiple Big-M (MBM)}
\end{subfigure}\hfill
\begin{subfigure}[t]{0.32\textwidth}
\centering
\begin{tikzpicture}
\begin{axis}[relaxpanel, xlabel = $z_1$, ylabel = $z_2$]
  \draw[fill=gray!15, draw=black, dashed, line width=1pt]
      (axis cs:1, 0) -- (axis cs:7, 0)
      arc[start angle=270, end angle=400.6,
          x radius=1, y radius=1]
      -- (axis cs:4.759, 5.151)
      arc[start angle=40.6, end angle=139.4,
          x radius=1, y radius=1]
      -- (axis cs:0.241, 1.651)
      arc[start angle=139.4, end angle=270,
          x radius=1, y radius=1]
      -- cycle;
  \threedisks
\end{axis}
\end{tikzpicture}
\caption{hull relaxation}
\end{subfigure}

\vspace{0.6em}
\centerline{\begin{tikzpicture}[baseline=-0.5ex]
  \draw[dashed, line width=1pt] (0,0) -- (0.7,0);
  \node[right=2pt, font=\footnotesize] at (0.7, 0) {envelope};
\end{tikzpicture}}
\caption{Relaxation envelopes for the three-disk GDP under big-M, MBM, and hull.}
\label{fig:relax_envelopes}
\end{figure}

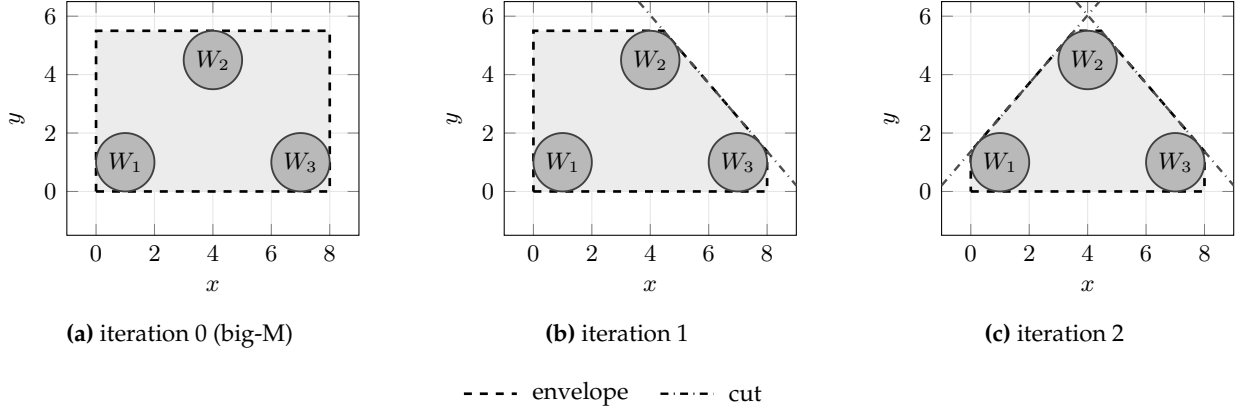
\begin{figure}[H]
\centering
\begin{subfigure}[t]{0.32\textwidth}
\centering
\begin{tikzpicture}
\begin{axis}[relaxpanel]
  \addplot[fill=gray!15, draw=black, dashed, line width=1pt]
      coordinates {(0,0) (8,0) (8,5.5) (0,5.5) (0,0)};
  \threedisks
\end{axis}
\end{tikzpicture}
\caption{iteration 0 (big-M)}
\end{subfigure}\hfill
\begin{subfigure}[t]{0.32\textwidth}
\centering
\begin{tikzpicture}
\begin{axis}[relaxpanel]
  \addplot[fill=gray!15, draw=black, dashed, line width=1pt]
      coordinates {(0,0) (8,0) (8,1.371) (4.458,5.5)
                   (0,5.5) (0,0)};
  \draw[draw=black!70, dash dot, line width=0.9pt]
      (axis cs:3.601, 6.5) -- (axis cs:9, 0.204);
  \threedisks
\end{axis}
\end{tikzpicture}
\caption{iteration 1}
\end{subfigure}\hfill
\begin{subfigure}[t]{0.32\textwidth}
\centering
\begin{tikzpicture}
\begin{axis}[relaxpanel]
  \addplot[fill=gray!15, draw=black, dashed, line width=1pt]
      coordinates {(0,0) (8,0) (8,1.371) (4.458,5.5)
                   (3.542,5.5) (0,1.371) (0,0)};
  \draw[draw=black!70, dash dot, line width=0.9pt]
      (axis cs:3.601, 6.5) -- (axis cs:9, 0.204);
  \draw[draw=black!70, dash dot, line width=0.9pt]
      (axis cs:-1, 0.204) -- (axis cs:4.401, 6.5);
  \threedisks
\end{axis}
\end{tikzpicture}
\caption{iteration 2}
\end{subfigure}

\vspace{0.6em}
\centerline{\begin{tikzpicture}[baseline=-0.5ex]
  \draw[dashed, line width=1pt] (0,0) -- (0.7,0);
  \node[right=2pt, font=\footnotesize] at (0.7, 0) {envelope};
  \draw[dash dot, line width=0.9pt] (2.6,0) -- (3.3,0);
  \node[right=2pt, font=\footnotesize] at (3.3, 0) {cut};
\end{tikzpicture}}
\caption{Cutting plane progression on the three-disk GDP: each
iteration appends a separating hyperplane that tightens the relaxed
big-M master toward the hull.}
\label{fig:cp_envelope}
\end{figure}

\paragraph{Relaxation Tightness.} These methods occupy distinct positions on the relaxation-tightness
spectrum between big-M and hull.
Figure~\ref{fig:relax_envelopes} illustrates this on a three-disk GDP
with disjuncts $D_i: (z_1 - c_{1,i})^2 + (z_2 - c_{2,i})^2 \leq 1$ centered
at $c_1 = (1,1)$, $c_2 = (4, 4.5)$, $c_3 = (7, 1)$. Big-M reduces each
disjunction to its bounding box, hull recovers the exact hull of
the three disks, and MBM sits strictly between by growing each disjunct
constraint to its maximum effective radius
$\sqrt{1 + \max_{i' \neq i} M_{ii'}}$. Because the per-constraint $M_{ii'}$
values cannot trace the curved hull boundary, MBM is strictly weaker
than hull on nonlinear convex disjuncts; CP closes the remaining gap
iteratively. Figure~\ref{fig:cp_envelope} shows this progression on
the same three-disk GDP: starting from the big-M bounding box, each
appended cut \eqref{eq:cp_cut} tightens the relaxed big-M master
toward the hull without disaggregating any variables.

\subsection{The InfiniteGDP Abstraction}\label{sec:bg_infgdp}

An InfiniteOpt problem is posed over a continuous domain $d \in \mathcal{D}$ (time, space, and/or uncertainty, with
$\mathcal{D} = \mathcal{D}_1 \times \cdots \times \mathcal{D}_n$ when several infinite parameters are
present), carrying an objective summarized by a measure operator
$M_d$ (e.g., a space--time integral or an expectation
\cite{gondosiswanto2025infiniteopt}) and constraints that hold for
all $d \in \mathcal{D}$. Infinite decision variables
$y : \mathcal{D} \mapsto \mathcal{Y} \subseteq \mathbb{R}^{n_y}$ carry the
parameter dependence $y(d)$; finite decision variables
$z \in \mathcal{Z} \subseteq \mathbb{R}^{n_z}$ do not. Where the
distinction is immaterial, the methods below write $y(d)$ for the
full decision vector. A transcription method
discretizes the problem onto a support grid
$\hat{\mathcal{D}} = \{\hat{d}_k : k \in \mathcal{K}\} \subset \mathcal{D}$: the measure
becomes a quadrature $\sum_{k} \omega_k\, f(y_k, \hat{d}_k)$ with
weights $\omega_k$, and each constraint is imposed pointwise at
$y_k \equiv y(\hat{d}_k)$, yielding a finite NLP
\cite{pulsipher2022infiniteopt, biegler2007overview}. Let
$\mathcal{T}_{\hat{\mathcal{D}}}$ denote the \emph{transcription map} that carries an
infinite quantity to its vector of support values,
\begin{equation}\label{eq:transcription}
    \mathcal{T}_{\hat{\mathcal{D}}} : y(\cdot) \;\longmapsto\;
    \bigl\{\, y(\hat{d}_k) : k \in \mathcal{K} \,\bigr\},
\end{equation}
and \emph{lifting} for the inverse mapping: replacing a finite
quantity by a $d$-indexed analogue and requiring every relation to hold
for all $d \in \mathcal{D}$. These two operations are incorporated into the methods described below and proposed in Section \ref{sec:reformulations}.

InfiniteGDP \cite{gondosiswanto2025infiniteopt} uses this
infinite-dimensional setting to generalize Problem \eqref{eq:bg_disj} by lifting the indicator variables and disjunctions
into the infinite domain. Each Boolean indicator becomes an infinite variable
$W_{ij} : \mathcal{D} \mapsto \{\text{True}, \text{False}\}$, and its
binary realization $w_{ij}(d) \in \{0, 1\}$ may switch
independently at every point of $\mathcal{D}$. This lifting operation produces the formulation:
\begin{equation}\label{eq:bg_inf_disj}
\begin{aligned}
    \min_{y,\, W} \quad
        & M_d\, f(y(d), d) \\
    \text{s.t.} \quad
        & g(y(d), d) \leq 0,
          \quad d \in \mathcal{D} \\
        & \bigvee_{i \in \mathcal{I}_j}
          \begin{bmatrix}
            W_{ij}(d) \\[2pt]
            r_{ij}(y(d), d) \leq 0
          \end{bmatrix},
          \quad j \in \mathcal{J},\; d \in \mathcal{D} \\
        & \Omega(W(d)) = \text{True},
          \quad d \in \mathcal{D} \\
        & y(d) \in \mathcal{Y},\;
          W_{ij}(d) \in \{\text{True}, \text{False}\}
\end{aligned}
\end{equation}
with disjunctions, propositions, and cardinality constraints
all applied over $\mathcal{D}$.

The work of \cite{gondosiswanto2025infiniteopt} generalizes the
big-M, hull, and indicator reformulations of \eqref{eq:bg_bigm}
and \eqref{eq:bg_hull} to this setting. As a representative case,
the big-M-reformulated InfiniteGDP is
\begin{equation}\label{eq:bg_inf_bigm}
\begin{aligned}
    \min_{y,\, w} \quad
        & M_d\, f(y(d), d) \\
    \text{s.t.} \quad
        & g(y(d), d) \leq 0,
          \quad d \in \mathcal{D} \\
        & r_{ij}(y(d), d) \leq M_{ij}\, (1 - w_{ij}(d)),
          \quad i \in \mathcal{I}_j,\; j \in \mathcal{J},\;
          d \in \mathcal{D} \\
        & \sum_{i \in \mathcal{I}_j} w_{ij}(d) = 1,
          \quad j \in \mathcal{J},\; d \in \mathcal{D} \\
        & y(d) \in \mathcal{Y},\;
          w_{ij}(d) \in \{0, 1\}.
\end{aligned}
\end{equation}
The hull-reformulated version follows analogously by lifting
the disaggregation and perspective constraints of
\eqref{eq:bg_hull} to $d \in \mathcal{D}$. Each reformulation in this
paper seeks to first \emph{lift} the
finite reformulation rule to the infinite domain by replacing each
finite quantity by its $d$-indexed analogue and requiring every
relation to hold for all $d \in \mathcal{D}$; this produces a symbolic,
still-infinite intermediate model. Then transcription via the map
$\mathcal{T}_{\hat{\mathcal{D}}}$ is used to obtain a finite MIP
(Figure~\ref{fig:reform_flow}
sketches this flow). Note that iterative algorithms like LOA require more careful treatment as discussed in Section \ref{sec:reformulations}.

These InfiniteGDP problems can readily be modeled and solved via the Julia
packages \texttt{DisjunctiveProgramming.jl} and its
extension \texttt{InfiniteDisjunctiveProgramming.jl}, which provide a
JuMP-native syntax for the Problem \eqref{eq:bg_inf_disj}
and allow users to select from a library of solution methods; the syntax is further discussed in Section~\ref{sec:implementation}. To date, only the big-M and
hull endpoints of the relaxation spectrum reviewed in
Section \ref{sec:bg_gdp} have been generalized to the InfiniteGDP setting,
alongside the indicator reformulation. The intermediate methods (MBM,
P-split, CP, LOA) have not been generalized and investigated in an InfiniteOpt context which is the focus of this work.

\section{Solution Methods for InfiniteGDPs}\label{sec:reformulations}

\begin{figure}[!htb]
\centering
\resizebox{\textwidth}{!}{%
\begin{tikzpicture}[
    flowarr/.style={-Latex, thick},
    srcbox/.style={draw, rounded corners=3pt, text width=4.6cm,
      minimum height=3.0cm, fill=white, align=center, font=\small,
      inner sep=5pt},
    fbox/.style={draw, rounded corners=3pt, text width=5.8cm,
      minimum height=3.0cm, fill=white, align=center, font=\small,
      inner sep=5pt},
    arrlbl/.style={font=\scriptsize, align=center}
]
  \node[srcbox] (src) at (0, 0)
    {\textbf{InfiniteGDPModel}\\[3pt]
     disjunctions over $d \in \mathcal{D}$\\
     (Eq.~\ref{eq:bg_inf_disj})};

  \node[fbox] (bigm) at (8.5, 4.0)
    {\textbf{big-M reformulated}\\
     {\scriptsize symbolic, infinite}\\[3pt]
     $M_{ij}$ relaxes constraint\\
     when disjunct is off\\
     (Eq.~\ref{eq:bg_inf_bigm})};
  \node[fbox] (bigm_t) at (17.0, 4.0)
    {\textbf{big-M transcribed MIP}\\
     {\scriptsize transcribed, finite}\\[3pt]
     $M_{ij}$ relaxes constraint\\
     when disjunct is off,\\
     at each support $\hat{d}_k \in \hat{\mathcal{D}}$};

  \node[fbox] (hull) at (8.5, 0)
    {\textbf{hull reformulated}\\
     {\scriptsize symbolic, infinite}\\[3pt]
     copy $\nu_{ij}(d)$ vanishes\\
     when disjunct is off\\
     (Eq.~\ref{eq:bg_hull} lifted to $d \in \mathcal{D}$)};
  \node[fbox] (hull_t) at (17.0, 0)
    {\textbf{hull transcribed MIP}\\
     {\scriptsize transcribed, finite}\\[3pt]
     copy $\nu_{ij,k}$ vanishes\\
     when disjunct is off,\\
     at each $\hat{d}_k \in \hat{\mathcal{D}}$};

  \node[fbox] (partition) at (8.5, -4.5)
    {\textbf{P-split: partition}\\
     {\scriptsize symbolic, infinite}\\[3pt]
     auxiliary envelopes\\
     $v_{ij,p}(d)$ per partition\\
     (Section~\ref{sec:psplit})};

  \draw[flowarr] (src.east) --
    node[arrlbl, above, sloped, midway]
      {big-M reform.\\(lift to $d\!\in\!\mathcal{D}$)}
    (bigm.west);
  \draw[flowarr] (src.east) --
    node[arrlbl, above, midway]
      {hull reform.\\(lift to $d\!\in\!\mathcal{D}$)}
    (hull.west);
  \draw[flowarr] (src.east) --
    node[arrlbl, below, sloped, midway,
         text width=2.8cm, align=center]
         {P-split reform.\\(partition $V_1, \ldots, V_P$)}
    (partition.west);

  \draw[flowarr] (partition.north) --
    node[arrlbl, right=2pt, midway, align=left]
         {apply hull\\on linking sum}
    (hull.south);

  \draw[flowarr] (bigm.east) --
    node[arrlbl, above, midway]
      {transcribe on $\hat{\mathcal{D}}$\\$y(d)\!\mapsto\! y_k$}
    (bigm_t.west);
  \draw[flowarr] (hull.east) --
    node[arrlbl, above, midway]
      {transcribe on $\hat{\mathcal{D}}$\\$y(d)\!\mapsto\! y_k$}
    (hull_t.west);

\end{tikzpicture}%
}
\caption{Lift-then-transcribe reformulation flow on an
InfiniteGDPModel, shown for the big-M, hull, and P-split paths.}
\label{fig:reform_flow}
\end{figure}
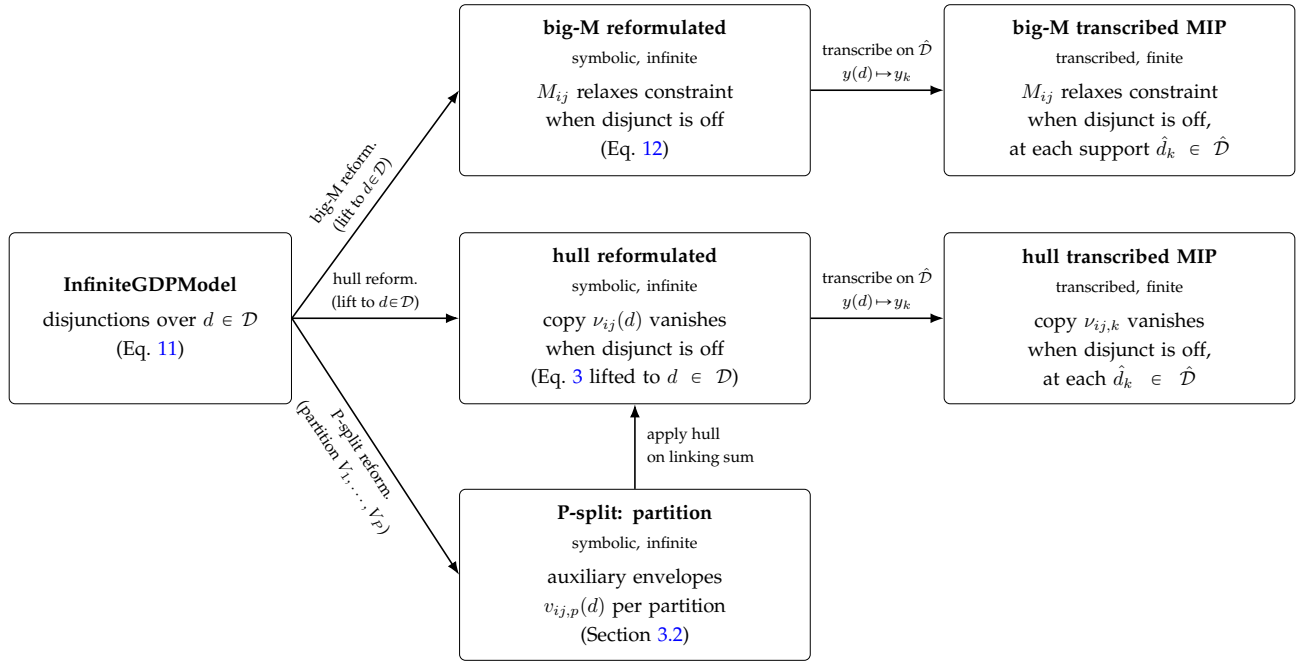

This section generalizes four solution methods to the InfiniteGDP
setting, each of which recovers its established finite form when the
model carries no infinite parameters. Three are reformulations that
fill the spectrum between the big-M and hull endpoints; the fourth,
LOA, is an iterative solution algorithm rather than a reformulation. MBM (Section~\ref{sec:mbm}) tightens big-M with
pairwise $M_{ii'}(d)$ without leaving the big-M form, and
Section~\ref{sec:mbmgp} introduces MBM-GP, a novel Gaussian-process
variant that avoids solving the $M$ subproblems at every support. P-split
(Section~\ref{sec:psplit}) parameterizes a continuous interpolation
between the two extremes via the partition size $P$. CP
(Section~\ref{sec:cp}) seeks hull-tight relaxations through iterative cuts added to big-M that avoid variable doubling.
LOA (Section~\ref{sec:loa}) targets convex GDPs with nonlinear
disjuncts where a direct hull would require a nonlinear MIP solver.
Section~\ref{sec:implementation} then shows how the methods are exposed
in software, and Section~\ref{sec:guidance} provides some general guidance and discussion on method selection. Each reformulation method broadly follows a lift-then-transcribe template described in Figure~\ref{fig:reform_flow}: the finite reformulation is first
lifted to a symbolic infinite model over $d \in \mathcal{D}$, then transcribed on
$\hat{\mathcal{D}}$ through the map $\mathcal{T}_{\hat{\mathcal{D}}}$ of
\eqref{eq:transcription}.

A recurring theme in what follows is that the total cost of solving
an InfiniteGDP splits into two components. The \emph{reformulation
time} covers everything needed to build the final program handed to
the solver, such as solving the per-pair $M$ subproblems of MBM or
generating cuts in CP. The \emph{solver time} is the cost of the
resulting MIP solve itself. The methods of this section occupy
different positions in this trade-off: a cheap reformulation
typically leaves a looser model that the solver must work harder to
close, while an expensive reformulation front-loads work to hand the
solver a tighter model. The reformulation is a one-time cost with respect to
the disjunctive structure, so in applications that repeatedly
re-solve the same model with updated data, such as model predictive
control or multi-objective optimization, it can be reused across
re-solves while the solver time recurs at every solve. The
benchmarks of Section~\ref{sec:casestudies} therefore report the two
components separately.

\subsection{Multiple Big-M}\label{sec:mbm}

The finite MBM reformulation of \cite{trespalacios2015mbm} replaces
each disjunct constraint with the per-pair form
\eqref{eq:mbm_finite}, where each constant $M_{ii'}$ is obtained
from the subproblem \eqref{eq:bg_mbm} of
Section~\ref{sec:bg_gdp}.
This construction lifts to the infinite setting: each
disjunct constraint $r_{ij}(y(d)) \leq 0$ must hold for all
$d \in \mathcal{D}$, so each per-pair big-M parameter becomes a
continuous function over the infinite domain,
$M_{ii'} : \mathcal{D} \mapsto \mathbb{R}$, and the reformulated
constraint
\begin{equation}\label{eq:mbm_inf}
    r_{ij}(y(d)) \leq \sum_{i' \neq i} M_{ii'}(d)\, w_{i'j}(d),
    \quad i \in \mathcal{I}_j,\; j \in \mathcal{J},\;
    d \in \mathcal{D},
\end{equation}
remains fully infinite. The function $M_{ii'}(d)$ is anchored by
pointwise data: transcribing through $\mathcal{T}_{\hat{\mathcal{D}}}$,
the subproblem
\begin{equation}\label{eq:mbm_inf_sub}
\begin{aligned}
    M_{ii',k} = \max_{y_k} \quad & r_{ij}(y_k) \\
    \text{s.t.} \quad & r_{i'j}(y_k) \leq 0 \\
    & y^L \leq y_k \leq y^U,
\end{aligned}
\quad k \in \mathcal{K},
\end{equation}
is solved for each $\hat{d}_k \in \hat{\mathcal{D}}$, where
$y_k \equiv y(\hat{d}_k)$. The per-support construction is only
needed when the subproblem actually varies over the domain. If
neither $r_{ij}$ nor $r_{i'j}$ depends on the infinite parameter $d$,
whether directly or through known continuous functions over the
infinite domain, then \eqref{eq:mbm_inf_sub} is the same program at
every support and $M_{ii'}(d)$ collapses to a scalar obtained from
a single solve. Otherwise, any continuous function over
$\mathcal{D}$ that matches or upper bounds the pointwise values is
a valid choice of $M_{ii'}(d)$. One natural implementation is a
multilinear interpolant of $\{M_{ii',k}\}$ built on the support
grid.
Alternatively, Section~\ref{sec:mbmgp} shows that a Gaussian process can supply a
valid $M_{ii'}(d)$ from subproblem solves at only a small subset of
the supports. This structure preserves the
tightness of the pointwise bound while avoiding a single conservative
scalar $M$ derived from global variable bounds. In practice, the
model structure often correlates the infinite parameters, causing
$M_{ii',k}$ to vary substantially across $\hat{\mathcal{D}}$ and
making $M_{ii'}(d)$ much tighter than the worst-case value.

\begin{figure}[H]
\centering
\resizebox{\textwidth}{!}{%
\definecolor{disji}{RGB}{31,119,180}
\definecolor{disjj}{RGB}{214,96,32}
\begin{tikzpicture}[node distance=1.0cm and 1.6cm,
    inputbox/.style={draw, rounded corners=3pt, text width=4.5cm,
      minimum height=2.0cm, fill=white, align=center, font=\small,
      inner sep=5pt},
    inputi/.style={inputbox, draw=disji!70!black, line width=0.7pt,
      fill=disji!10},
    inputj/.style={inputbox, draw=disjj!70!black, line width=0.7pt,
      fill=disjj!10},
    subbox/.style={draw, rounded corners=3pt, text width=3.4cm,
      minimum height=2.6cm, fill=white, align=center, font=\small,
      inner sep=5pt},
    outbox/.style={draw, rounded corners=3pt, fill=white,
      align=center, font=\small, inner sep=5pt}]

  \node[inputi] (disj_i)
    {disjunct $i$ constraint\\[2pt]
     $r_{ij}(y(d)) \leq 0$};
  \node[inputj, below=0.5cm of disj_i] (disj_j)
    {disjunct $i' \neq i$ region\\[2pt]
     $r_{i'j}(y(d)) \leq 0$};

  \coordinate (input_mid) at
    ($(disj_i.east)!0.5!(disj_j.east)$);

  \node[draw, rounded corners=3pt, fill=gray!8, align=center,
        font=\scriptsize, inner sep=5pt, text width=2.2cm,
        minimum height=3.0cm, anchor=west]
    at ($(input_mid) + (1.7, 0)$) (trans)
    {\textbf{transcribe} on $\hat{\mathcal{D}}$\\[6pt]
     $y(d)\!\mapsto\! y_k$};

  \node[subbox, anchor=west]
    at ($(input_mid) + (5.4, 0)$) (sub)
    {$M_{ii',k} = \max_{y_k}\;$%
     {\setlength{\fboxsep}{1.5pt}%
       \colorbox{disji!30}{$r_{ij}(y_k)$}}\\[2pt]
     s.t.\ {\setlength{\fboxsep}{1.5pt}%
       \colorbox{disjj!30}{$r_{i'j}(y_k) \leq 0$}},\\
     $y^L \leq y_k \leq y^U$};

  \node[font=\small, anchor=north, align=center] (forEachLabel)
    at ($(sub.south) + (0, -0.10)$)
    {one subproblem per\\support $\hat{d}_k \in \hat{\mathcal{D}}$};

  \begin{scope}[on background layer]
    \node[subbox, anchor=center] (back)
      at ($(sub.center) + (0.80, 1.00)$) {};
    \node[subbox, anchor=center] (mid)
      at ($(sub.center) + (0.40, 0.50)$) {};
  \end{scope}

  \node[font=\scriptsize, anchor=north west, inner sep=6pt]
    at (sub.north west) {$k\!=\!1$};
  \node[font=\scriptsize, anchor=north west, inner sep=4pt]
    at (mid.north west) {$\cdots$};
  \node[font=\scriptsize, anchor=north west, inner sep=4pt]
    at (back.north west) {$k\!=\!K$};

  \begin{scope}[on background layer]
    \node[draw, rounded corners=5pt,
          fit=(sub)(back)(forEachLabel),
          inner sep=8pt] (container) {};
  \end{scope}

  \node[outbox, right=3.4cm of container] (Mt)
    {\begin{tikzpicture}
       \draw[->] (0,0) -- (4.4, 0)
         node[right, font=\tiny, inner sep=1pt] {$d$};
       \draw[->] (0,0) -- (0, 3.7);
       \node[font=\small, anchor=north west]
         at (0.15, 3.65) {$M_{ii'}(d)$};
       \draw[thick] plot[smooth, tension=0.7] coordinates {
         (0.30, 0.80) (0.62, 1.10) (0.95, 1.45) (1.30, 1.95)
         (1.65, 2.50) (2.00, 2.30) (2.35, 2.00) (2.70, 2.25)
         (3.05, 2.65) (3.40, 2.95) (3.70, 3.20) (3.95, 2.40)
         (4.20, 1.40)
       };
       \foreach \x/\y in {0.30/0.80, 0.62/1.10, 0.95/1.45, 1.30/1.95,
                          1.65/2.50, 2.00/2.30, 2.35/2.00, 2.70/2.25,
                          3.05/2.65, 3.40/2.95, 3.70/3.20, 3.95/2.40,
                          4.20/1.40}
         \fill (\x, \y) circle (1.4pt);
       \node[font=\tiny, anchor=north] at (0.30, 0) {$k\!=\!1$};
       \node[font=\tiny, anchor=north] at (4.20, 0) {$k\!=\!K$};
     \end{tikzpicture}};

  \draw[flowarr, draw=disji!70!black] (disj_i.east) --
    node[above, font=\scriptsize, midway] {objective}
    (trans.west |- disj_i.east);
  \draw[flowarr, draw=disjj!70!black] (disj_j.east) --
    node[below, font=\scriptsize, midway] {feasibility}
    (trans.west |- disj_j.east);
  \draw[flowarr] (trans.east) -- (container.west |- trans.east);
  \draw[flowarr] (container) --
    node[above, font=\scriptsize, midway, align=center]
      {lift to $d\!\in\!\mathcal{D}$\\$\{M_{ii',k}\}\!\mapsto\! M_{ii'}(d)$}
    (Mt);

\end{tikzpicture}%
}
\caption{Multiple big-M in the infinite setting: the per-pair
subproblem is solved at every support point $\hat{d}_k$ and the
resulting $M_{ii',k}$ realize the continuous big-M function
$M_{ii'}(d)$ over the infinite domain.}
\label{fig:flow_mbm}
\end{figure}
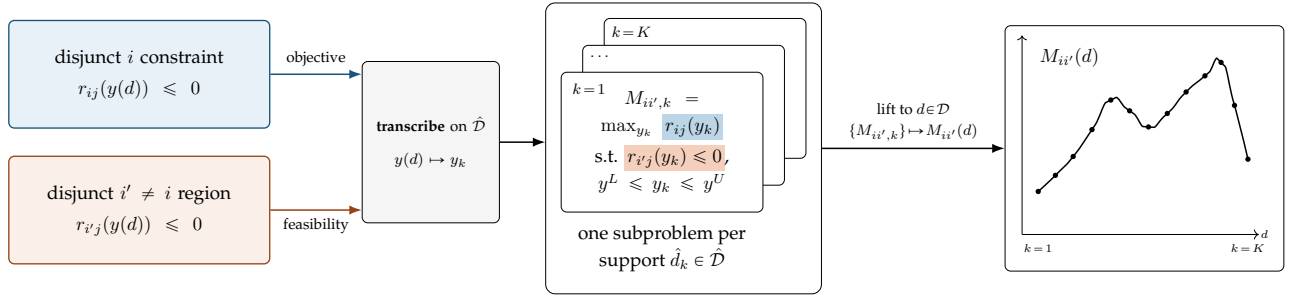

Each subproblem inherits the support grid of the original
InfiniteGDP model, so any known
continuous functions over the infinite domain appearing in
$r_{ij}, r_{i'j}$ evaluate at the matching point $\hat{d}_k$.
Figure~\ref{fig:flow_mbm} summarizes this per-pair, per-support
construction of $M_{ii'}(d)$.

\subsubsection{Gaussian Process Multiple Big-M (MBM-GP)}
\label{sec:mbmgp}

Viewing $M_{ii'}(d)$ as a continuous function over the infinite
domain enables an alternative to solving \eqref{eq:mbm_inf_sub} at
every support point. The grid construction above requires one
subproblem solve per constraint pair per support, so its cost grows
with the resolution of $\hat{\mathcal{D}}$, and on nonconvex models
each subproblem is itself a nonconvex NLP. MBM-GP instead
solves \eqref{eq:mbm_inf_sub} at only a small subset
$\tilde{\mathcal{D}} \subset \hat{\mathcal{D}}$ of the supports,
fits a Gaussian process (GP) \cite{rasmussen2006gaussian} to the
computed values, and takes the
upper bound of the GP prediction as the big-M function at every
unsolved support:
\begin{equation}\label{eq:mbmgp_ucb}
    M_{ii'}(d) = \mu(d) + \kappa\,\sigma(d),
\end{equation}
where $\mu(d)$ and $\sigma(d)$ are the GP posterior mean and
standard deviation and $\kappa \geq 0$ is a confidence factor. Taking the
upper bound ensures validity: a prediction that undershoots
the true pointwise value would make the relaxed constraint too
tight and cut off feasible points, whereas overshooting only
loosens the relaxation. Using \eqref{eq:mbmgp_ucb} therefore keeps
MBM-GP a valid choice of $M_{ii'}(d)$ in the sense described above
while solving only a fraction of the subproblems.

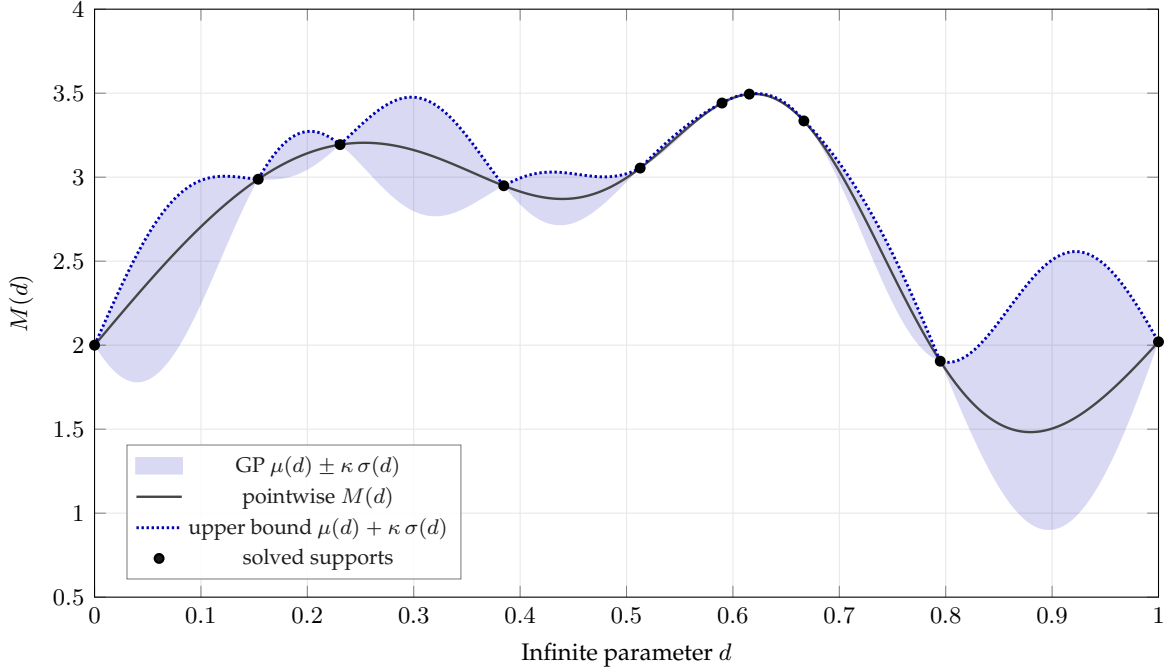
\begin{figure}[H]
\centering
\begin{tikzpicture}
\begin{axis}[trajpanel, width=0.92\textwidth,
  height=0.55\textwidth, xlabel={Infinite parameter $d$},
  ylabel={$M(d)$}, xmin=0, xmax=1, ymin=0.5, ymax=4,
  legend style={at={(0.03,0.03)}, anchor=south west}]
  \addplot[name path=ucbp, draw=none, forget plot]
    table[col sep=comma, x=t, y=ucb] {data/mbmgp_demo_fit.csv};
  \addplot[name path=lcbp, draw=none, forget plot]
    table[col sep=comma, x=t, y=lcb] {data/mbmgp_demo_fit.csv};
  \addplot[blue!70!black, opacity=0.15]
    fill between[of=ucbp and lcbp];
  \addplot[gray!55!black] table[col sep=comma, x=t, y=M]
    {data/mbmgp_demo_true.csv};
  \addplot[blue!70!black, densely dotted, line width=1.1pt]
    table[col sep=comma, x=t, y=ucb] {data/mbmgp_demo_fit.csv};
  \addplot[only marks, mark=*, mark size=1.6pt, black]
    table[col sep=comma, x=t, y=M] {data/mbmgp_demo_pts.csv};
  \legend{GP $\mu(d) \pm \kappa\,\sigma(d)$, pointwise $M(d)$,
    upper bound $\mu(d) + \kappa\,\sigma(d)$, solved supports}
\end{axis}
\end{tikzpicture}
\caption{MBM-GP on an example pointwise $M$ profile over $40$
supports, of which $10$ are solved (markers) and keep their exact
values. The shaded region is the GP fitted to the solved supports
with $\kappa = 2.5$; its upper edge (dotted) is the bound
\eqref{eq:mbmgp_ucb} that the unsolved supports take, and it lies
above the pointwise $M$ values throughout.}
\label{fig:mbmgp_selection}
\end{figure}

The scalar collapse of Section~\ref{sec:mbm} applies unchanged:
a pair whose subproblem does not depend on the infinite parameter
is solved once, and the GP is reserved for pairs whose
$M_{ii'}(d)$ genuinely varies. For those, an implementation must
make two choices: which supports to solve, and how strongly the
prediction is inflated to an upper bound.
\begin{itemize}
\item \emph{Point selection.} A uniform criterion fixes
  $\tilde{\mathcal{D}}$ a priori, e.g., as evenly spaced supports
  across the domain. An adaptive criterion instead builds
  $\tilde{\mathcal{D}}$ sequentially. It solves the support with
  the largest variance \eqref{eq:mbmgp_ucb} and refits the GP,
  so the solves concentrate where the predicted $M$ is most uncertain. Sampling can continue until the posterior
  variance is closed by a prescribed amount. This sequential approach is similar to Bayesian optimization with an an exploration focused acquisition function. Under either
  criterion, solved supports keep their exact values and grid MBM
  is recovered as $\tilde{\mathcal{D}} \to \hat{\mathcal{D}}$.
\item \emph{Upper bound.} The confidence factor $\kappa$ in
  \eqref{eq:mbmgp_ucb} sets how far the prediction is inflated
  above the posterior mean. A larger $\kappa$ makes an undershoot,
  and with it an invalid $M$, less likely, at the price of extra
  slack $\kappa\,\sigma(d)$ at the unsolved supports. A smaller
  $\kappa$ keeps the relaxation tighter but relies more heavily on
  the quality of the GP fit.
\end{itemize}
Figure~\ref{fig:mbmgp_selection} illustrates the construction on
an example $M$ profile. Of $40$ supports, $10$ are solved and keep
their exact values. The GP fitted to those solves supplies the
upper bound \eqref{eq:mbmgp_ucb} at the remaining supports, and
the bound lies above the pointwise values throughout.

For the case studies of Section~\ref{sec:casestudies}, we use the
adaptive criterion, seeding the GP with the first, middle, and
last supports and solving roughly a quarter of the supports, and
we set $\kappa = 2.5$. Under the GP posterior this corresponds to
a roughly $99\%$ one-sided confidence level
($\Phi(2.5) \approx 0.994$), so the prediction overshoots the true
value with high probability, which preserves validity at the cost
of mild extra slack. The benchmarks
find that MBM-GP recovers objectives comparable to grid MBM at a
small fraction of the reformulation cost.

\subsection{P-Split}\label{sec:psplit}

The finite P-split formulation of \cite{kronqvist2022psplit} splits
each disjunct constraint into the globally enforced envelopes
\eqref{eq:psplit_part} and the disjunct linking constraint
\eqref{eq:psplit_sum} of Section~\ref{sec:bg_gdp}, with hull
disaggregation applied only to the latter.

The split lifts to the infinite setting: the auxiliary variable
for a partition whose defining expression depends on infinite
parameters must itself be infinite, so each $v_{ij,p}$ is lifted to an
function inheriting its parameter dependence from
$\tilde{r}_p(y_{V_p}(d))$:
\begin{align}
    \tilde{r}_p(y_{V_p}(d)) &\leq v_{ij,p}(d),
      \quad p = 1,\ldots,P,\; i \in \mathcal{I}_j,\;
      j \in \mathcal{J},\; d \in \mathcal{D},
      \label{eq:psplit_inf_part} \\
    \sum_{p=1}^{P} v_{ij,p}(d) + c(d) &\leq 0,
      \quad i \in \mathcal{I}_j,\; j \in \mathcal{J},\;
      d \in \mathcal{D}, \label{eq:psplit_inf_sum}
\end{align}
where $c(d)$ captures any terms that are known continuous
functions over the infinite domain. The hull reformulation of
\cite{gondosiswanto2025infiniteopt} handles the infinite hull
disaggregation of the linking constraint \eqref{eq:psplit_inf_sum},
and the bounds on $v_{ij,p}(d)$ follow from the global bounds on $y_{V_p}$
over $\mathcal{D}$. Transcription on $\hat{\mathcal{D}}$ is then inherited from that
hull machinery, which discretizes the disaggregated linking
constraint through $\mathcal{T}_{\hat{\mathcal{D}}}$; P-split adds no
transcription step of
its own.

To make explicit which quantities are disaggregated, we write out
this infinite hull reformulation of the linking constraint
\eqref{eq:psplit_inf_sum}. It introduces per-disjunct copies
$\nu_{ij,p}^{i'}(d)$ of each auxiliary function, one for every
disjunct $i' \in \mathcal{I}_j$ of the disjunction, and yields
\begin{subequations}\label{eq:psplit_inf_hull}
\begin{align}
    & v_{ij,p}(d) = \sum_{i' \in \mathcal{I}_j} \nu_{ij,p}^{i'}(d),
      && p = 1,\ldots,P,\; i \in \mathcal{I}_j,\;
      d \in \mathcal{D},
      \label{eq:psplit_hull_agg} \\
    & \sum_{p=1}^{P} \nu_{ij,p}^{i}(d) + c(d)\, w_{ij}(d) \leq 0,
      && i \in \mathcal{I}_j,\; d \in \mathcal{D},
      \label{eq:psplit_hull_link} \\
    & w_{i'j}(d)\, v_{ij,p}^L \leq \nu_{ij,p}^{i'}(d)
      \leq w_{i'j}(d)\, v_{ij,p}^U,
      && p = 1,\ldots,P,\; i, i' \in \mathcal{I}_j,\;
      d \in \mathcal{D},
      \label{eq:psplit_hull_bnd}
\end{align}
\end{subequations}
for each disjunction $j \in \mathcal{J}$, where
$v_{ij,p}^L, v_{ij,p}^U$ are the interval-arithmetic bounds on the
auxiliary functions. Because the linking constraint is linear in
$v_{ij,p}(d)$, its perspective transform in
\eqref{eq:psplit_hull_link} is exact and introduces no
nonlinearity. In contrast with the full hull reformulation
\eqref{eq:bg_hull}, only the
auxiliary functions $v_{ij,p}(d)$ acquire disaggregated copies,
while the original infinite variables $y(d)$ are never duplicated.
They appear in the globally enforced envelopes
\eqref{eq:psplit_inf_part} rather than in the disjunctive
structure, so the nonlinear partition expressions also generate no
perspective terms of the form $w\,\tilde{r}_p(\cdot/w)$. Larger $P$
disaggregates more auxiliary functions and tightens the relaxation
toward hull, while $y(d)$ remains undisaggregated at every $P$.
This is the precise sense in which P-split interpolates between the
big-M and hull endpoints.

\subsection{Cutting Planes}\label{sec:cp}

The finite CP algorithm of Section~\ref{sec:bg_gdp} alternates
between the relaxed big-M master and the hull-based separation
subproblem \eqref{eq:cp_sep}, appending the separating cut
\eqref{eq:cp_cut} until the projection distance falls below the
tolerance $\varepsilon$ (Figure~\ref{fig:cp_envelope}).

In the infinite setting, the master is the big-M reformulation
of the InfiniteGDP, with disjunct constraints
$r_{ij}(y(d)) \leq M_{ij} (1 - w_{ij}(d))$ for $d \in \mathcal{D}$, and the
separation subproblem is the hull reformulation of
\cite{gondosiswanto2025infiniteopt}. The two are coupled through
the separation objective, which measures the distance between
$\boldsymbol{y}$ and the master solution
$\boldsymbol{y}^{\mathrm{rBM}}$ over the infinite domain. The
natural symbolic form of this objective is the measured distance
$M_d\, \|\boldsymbol{y}(d) - \boldsymbol{y}^{\mathrm{rBM}}(d)\|_2^2$.
Transcription carries the measure to its quadrature, so the
objective becomes the weighted sum in \eqref{eq:cp_inf_sep}, with
$\omega_k$ the quadrature weights from the transcription.

\begin{figure}[!htb]
\centering
\resizebox{\textwidth}{!}{%
\begin{tikzpicture}[
    cpbox/.style={draw, rounded corners=3pt, text width=2.4cm,
      minimum height=1.4cm, fill=white, align=center, font=\small,
      inner sep=4pt},
    origbox/.style={draw, rounded corners=3pt, text width=2.4cm,
      minimum height=3.0cm, fill=white, align=center, font=\small,
      inner sep=4pt},
    solidarr/.style={-Latex, thick},
    dashedarr/.style={-Latex, thick, dashed},
    steplabel/.style={font=\scriptsize, align=center, midway}
]
  \node[origbox] (gdp) at (0, 0) {InfiniteGDP\\model};

  \node[cpbox] (sep) at ($(gdp.east) + (5.0, 1.6)$) {Separation\\subproblem\\\scriptsize (hull copy,\\\scriptsize binaries relaxed)};
  \node[cpbox] (rBM) at ($(gdp.east) + (5.0, -1.6)$) {big-M\\master\\\scriptsize (BigM in place,\\\scriptsize binaries relaxed)};

  \draw[solidarr] (gdp.east |- sep.west) --
    node[steplabel, above] {(1) copy + hull
      reformulate\\+ relax integrality\\(over $d\!\in\!\mathcal{D}$)}
    (sep.west);
  \draw[solidarr] (gdp.east |- rBM.west) --
    node[steplabel, below] {(2) BigM reformulate
      in place\\+ relax integrality\\(over $d\!\in\!\mathcal{D}$)}
    (rBM.west);

  \draw[dashedarr] ([xshift=-8pt]rBM.north) --
    node[steplabel, left]
      {(3) master solution\\(sets separation objective)}
    ([xshift=-8pt]sep.south);
  \draw[dashedarr] ([xshift=8pt]sep.south) --
    node[steplabel, right]
      {(4) new cut\\(added to master)}
    ([xshift=8pt]rBM.north);

  \node[cpbox, right=3.4cm of rBM] (final)
    {transcribed\\MIP\\\scriptsize (to solver)};
  \draw[solidarr] (rBM.east) --
    node[steplabel, below]
      {(5) transcribe on $\hat{\mathcal{D}}$\\$y(d)\!\mapsto\! y_k$}
    (final.west);

  \begin{scope}[shift={($(sep.east) + (1.6, -0.8)$)}]
    \draw[solidarr] (0, 0.35) -- (0.7, 0.35);
    \node[anchor=west, font=\scriptsize] at (0.75, 0.35)
      {Model transformation};
    \draw[dashedarr] (0, -0.05) -- (0.7, -0.05);
    \node[anchor=west, font=\scriptsize] at (0.75, -0.05)
      {Information transfer};
  \end{scope}
\end{tikzpicture}%
}
\caption{Cutting plane algorithm in the InfiniteGDP setting: model
transformations (solid) and per-iteration information flow (dashed)
between a hull-reformulated separation copy and an in-place big-M
master.}
\label{fig:cp_flow}
\end{figure}
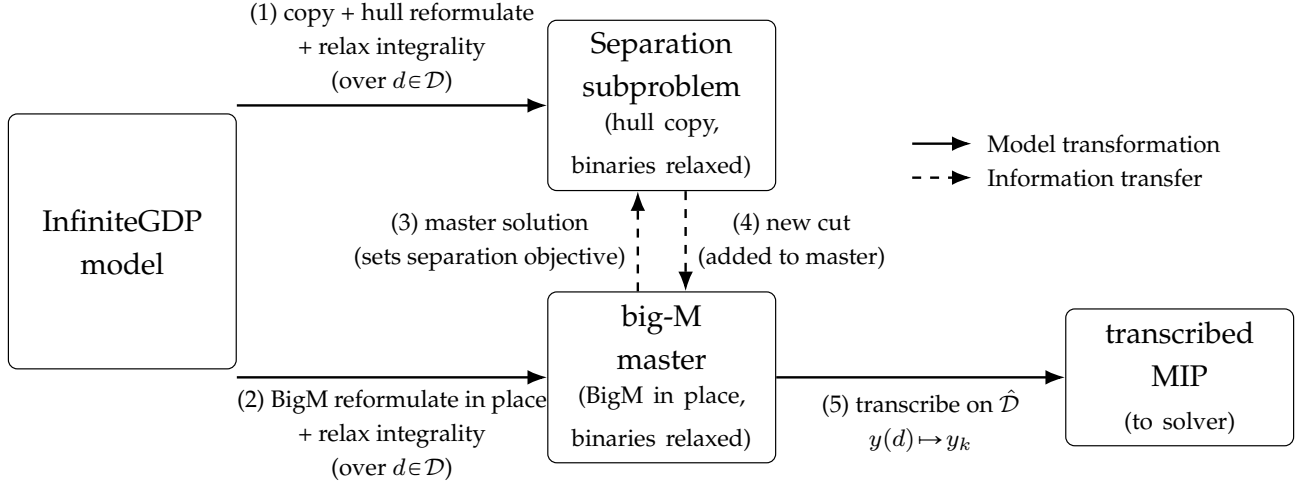

\begin{figure}[H]
\centering
\resizebox{\textwidth}{!}{%
\begin{tikzpicture}[
    cpbox/.style={draw, rounded corners=3pt, text width=3.6cm,
      minimum height=1.8cm, fill=white, align=center, font=\small,
      inner sep=4pt},
    cpsmall/.style={draw, rounded corners=2pt, text width=3.4cm,
      minimum height=0.8cm, fill=white, align=center, font=\scriptsize,
      inner sep=3pt},
    cpbig/.style={draw, rounded corners=3pt, text width=5.6cm,
      minimum height=2.6cm, fill=white, align=center, font=\small,
      inner sep=4pt},
    cpdiamond/.style={draw, diamond, aspect=1.6, align=center,
      inner sep=2pt, font=\scriptsize, text width=1.4cm}]

  \node[cpbox] (solveRbm)
    {solve master,\\obtain $\boldsymbol{y}^{\mathrm{rBM}}(d)$};

  \node[cpsmall, right=2.8cm of solveRbm, yshift=1.3cm] (sep1)
    {$\boldsymbol{y}^{\mathrm{rBM}}_1$};
  \node[cpsmall, below=0.15cm of sep1] (sep2)
    {$\boldsymbol{y}^{\mathrm{rBM}}_2$};
  \node[cpsmall, draw=none, fill=none, below=0.05cm of sep2]
        (sepDots) {$\vdots$};
  \node[cpsmall, below=0.05cm of sepDots] (sepK)
    {$\boldsymbol{y}^{\mathrm{rBM}}_K$};

  \begin{scope}[on background layer]
    \node[draw, rounded corners=5pt,
          fit=(sep1)(sepK), inner sep=8pt,
          label={[font=\scriptsize, align=center, yshift=3pt]%
            above:{transcribe on $\hat{\mathcal{D}}$\\
            $\boldsymbol{y}^{\mathrm{rBM}}(d)\!\mapsto\!
            \boldsymbol{y}^{\mathrm{rBM}}_k$}}]
          (sepContainer) {};
  \end{scope}

  \node[cpbig, right=1.2cm of sepContainer] (sepBig)
    {Separation subproblem\\[3pt]
     $\min \sum_{k=1}^{K} \omega_k
        \|\boldsymbol{y}_k
        - \boldsymbol{y}^{\mathrm{rBM}}_k\|^2$\\
     s.t.\ hull reformulation};

  \node[cpdiamond, right=1.0cm of sepBig] (gap)
    {gap $< \varepsilon$?};

  \node[cpbox, right=0.7cm of gap, text width=3.0cm,
        minimum height=1.8cm] (done)
    {solve master\\(with all cuts)};

  \node[cpbox, below=1.2cm of sepContainer, text width=6.0cm,
        minimum height=1.4cm] (cut)
    {add cut to master\\[2pt]
     {\scriptsize one global cut, aggregated over $\hat{\mathcal{D}}$}\\
     (Eq.~\ref{eq:cp_inf_cut})};

  \draw[flowarr] (solveRbm.east)
    -- (sepContainer.west |- solveRbm.east);
  \draw[flowarr] (sepContainer.east) -- (sepBig.west);
  \draw[flowarr] (sepBig) -- (gap);
  \draw[flowarr] (gap) --
    node[above, font=\scriptsize] {yes} (done);
  \draw[flowarr] (gap.south) |-
    node[font=\scriptsize, pos=0.05, right=2pt] {no}
    (cut.east);
  \draw[flowarr] (cut.west) -| (solveRbm.south);

\end{tikzpicture}%
}
\caption{Per-iteration mechanics of the cutting plane loop: the
transcribed master solution $\{\boldsymbol{y}^{\mathrm{rBM}}_k\}$
parameterizes the separation subproblem, whose minimizer yields one
global cut appended to the master.}
\label{fig:cp_iteration}
\end{figure}
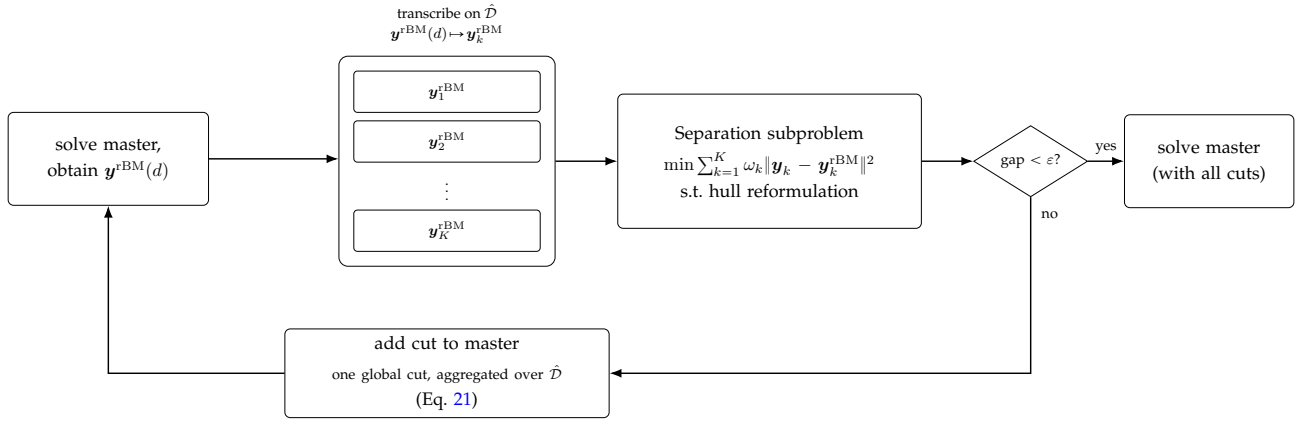

Lifting places both the master and the separation subproblem over
$d \in \mathcal{D}$; the transcription map $\mathcal{T}_{\hat{\mathcal{D}}}$ of
\eqref{eq:transcription}
then carries them to the support grid, sending each decision
variable $\boldsymbol{y}(d)$ to
$\{\boldsymbol{y}(\hat{d}_k) : k \in \mathcal{K}\}$. At
iteration $\ell$, the master solution
$\{\boldsymbol{y}^{\mathrm{rBM}}(\hat{d}_k) : k \in \mathcal{K}\}$
parameterizes the hull-based separation
\begin{equation}\label{eq:cp_inf_sep}
\begin{aligned}
    \min_{\substack{\boldsymbol{y}(\hat{d}_k),\,
        \boldsymbol{\nu}_{ij}(\hat{d}_k),\\ w_{ij}(\hat{d}_k)}}
    \quad & \sum_{k \in \mathcal{K}} \omega_k\,
      \|\boldsymbol{y}(\hat{d}_k)
       - \boldsymbol{y}^{\mathrm{rBM}}(\hat{d}_k)\|_2^2 \\
    \text{s.t.} \quad
    & \boldsymbol{y}(\hat{d}_k)
      = \sum_{i \in \mathcal{I}_j} \boldsymbol{\nu}_{ij}(\hat{d}_k),
      && j \in \mathcal{J},\; k \in \mathcal{K} \\
    & w_{ij}(\hat{d}_k)\, r_{ij}\bigl(\boldsymbol{\nu}_{ij}(\hat{d}_k)
        / w_{ij}(\hat{d}_k)\bigr) \leq 0,
      && i \in \mathcal{I}_j,\; j \in \mathcal{J},\, k \in \mathcal{K} \\
    & w_{ij}(\hat{d}_k)\, y^L \leq \boldsymbol{\nu}_{ij}(\hat{d}_k)
      \leq w_{ij}(\hat{d}_k)\, y^U,
      && i \in \mathcal{I}_j,\; j \in \mathcal{J},\, k \in \mathcal{K} \\
    & \sum_{i \in \mathcal{I}_j} w_{ij}(\hat{d}_k) = 1,
      \quad w_{ij}(\hat{d}_k) \in [0, 1],
      && i \in \mathcal{I}_j,\; j \in \mathcal{J},\, k \in \mathcal{K}.
\end{aligned}
\end{equation}
The optimal $\boldsymbol{y}(\hat{d}_k)$ is the per-support projection
$\boldsymbol{y}^{\mathrm{SEP}}(\hat{d}_k)$ of
$\boldsymbol{y}^{\mathrm{rBM}}(\hat{d}_k)$ onto the hull, and
all supports $k \in \mathcal{K}$ are solved jointly. Since
\eqref{eq:cp_inf_sep} projects under the quadrature-weighted
distance, the optimality condition of the projection supplies a
separating hyperplane in the same weighted geometry. The
separating cut is the quadrature of the measured cut
\begin{equation}\label{eq:cp_inf_cut_measure}
    M_d\, 2\bigl(\boldsymbol{y}^{\mathrm{SEP}}(d)
        - \boldsymbol{y}^{\mathrm{rBM}}(d)\bigr)
    \bigl(\boldsymbol{y}(d)
        - \boldsymbol{y}^{\mathrm{SEP}}(d)\bigr)
    \geq 0,
\end{equation}
which transcribes to
\begin{equation}\label{eq:cp_inf_cut}
    \sum_{k \in \mathcal{K}} \omega_k\,
    2\bigl(\boldsymbol{y}^{\mathrm{SEP}}(\hat{d}_k)
        - \boldsymbol{y}^{\mathrm{rBM}}(\hat{d}_k)\bigr)
    \bigl(\boldsymbol{y}(\hat{d}_k)
        - \boldsymbol{y}^{\mathrm{SEP}}(\hat{d}_k)\bigr)
    \geq 0
\end{equation}
and is appended to the transcribed master, which acts on the transcribed variables
that $\mathcal{T}_{\hat{\mathcal{D}}}$ produces. The weights in
\eqref{eq:cp_inf_cut} match those in the separation objective,
which keeps the transcribed cut consistent with the weighted
geometry of the projection. The cut is a supporting hyperplane
of the hull at $\boldsymbol{y}^{\mathrm{SEP}}$ whose normal
$\boldsymbol{y}^{\mathrm{rBM}} - \boldsymbol{y}^{\mathrm{SEP}}$ is
taken in the weighted geometry. The weighting also makes the
termination test resolution-independent. As supports are
added, the weighted projection distance converges to the measured
distance, so the tolerance $\varepsilon$ keeps its meaning under
grid refinement, whereas an unweighted sum grows with
$|\mathcal{K}|$. Figure~\ref{fig:cp_flow} shows how the master and separation
subproblem are derived from the InfiniteGDP model, and
Figure~\ref{fig:cp_iteration} details the per-iteration cut loop.

When the original InfiniteGDP model contains no coupling between
support points (e.g., a two-stage stochastic program with no first stage decisions), a stronger form of cut becomes available. Coupling
here means any constraint that links the decision at one point of
the domain to the decision at another, such as derivative
constraints or measure constraints that aggregate over
$\mathcal{D}$. When no such constraint is present, the transcribed
feasible region factors into a Cartesian product across supports:
the constraints at support $k$ involve only
$\boldsymbol{y}(\hat{d}_k)$, and the hull of a Cartesian product is
the product of the per-support hulls. The separation subproblem
\eqref{eq:cp_inf_sep} then decouples into $|\mathcal{K}|$
independent projections, and each block of its solution,
$\boldsymbol{y}^{\mathrm{SEP}}(\hat{d}_k)$, is the projection of
$\boldsymbol{y}^{\mathrm{rBM}}(\hat{d}_k)$ onto the hull at that
support alone. Each projection therefore supplies its own valid
separating hyperplane at its support, exactly as in the finite
algorithm, and lifting this family back to the infinite domain
yields the pointwise infinite cut
\begin{equation}\label{eq:cp_inf_cut_pointwise}
    2\bigl(\boldsymbol{y}^{\mathrm{SEP}}(d)
        - \boldsymbol{y}^{\mathrm{rBM}}(d)\bigr)
    \bigl(\boldsymbol{y}(d)
        - \boldsymbol{y}^{\mathrm{SEP}}(d)\bigr)
    \geq 0,
    \quad d \in \mathcal{D},
\end{equation}
which holds at every $d \in \mathcal{D}$ rather than only in aggregate. The
pointwise form is strictly stronger. Aggregating
\eqref{eq:cp_inf_cut_pointwise} over $\mathcal{K}$ with the
weights $\omega_k$ recovers \eqref{eq:cp_inf_cut}, but the
converse fails because individual inner products may be negative
as long as the weighted total is non-negative. The weights are
immaterial within each pointwise cut, since scaling by
$\omega_k > 0$ leaves its halfspace unchanged.
Each separation solve thus contributes $|\mathcal{K}|$ cuts instead
of one. This cut family has no finite-GDP counterpart: it exists
only because the infinite domain replicates the disjunctive
structure at every support point. The summed cut
\eqref{eq:cp_inf_cut} remains the valid fallback whenever
cross-support coupling is present.


\subsection{Logic-Based Outer Approximation}\label{sec:loa}

The finite LOA algorithm of Section~\ref{sec:bg_gdp} alternates
between the fixed-combination of an NLP subproblem~$(S')$ and the OA
master~$(\hat{M}^b_{OA})$, accumulating first-order linearizations
until the bound gap $|Z_U - Z_L|$ falls below $\varepsilon$.

In the infinite setting, $W_{ij}(d)$ varies independently at each
point in the domain, so lifting the combination yields a
\emph{trajectory}: a map from the support grid to disjunct selections,
$\{i_j^*(k)\}_{k \in \mathcal{K}}$ for each disjunction~$j$.
Transcribing through the map $\mathcal{T}_{\hat{\mathcal{D}}}$ of
\eqref{eq:transcription}, the NLP subproblem and master are both
posed on the resulting transcribed representation, in the same
manner as the cutting plane method of Section~\ref{sec:cp}.

\paragraph{NLP subproblem.}
Given a fixed Boolean trajectory $W_{ij}^l(\hat{d}_k) \in
\{\text{True}, \text{False}\}$ for all $i \in \mathcal{I}_j$,
$j \in \mathcal{J}$, $k \in \mathcal{K}$, the infinite NLP
subproblem is
\begin{equation}\label{eq:loa_nlp}
\begin{aligned}
    (S'_\infty) : \quad \min_{y_k} \quad
        & \sum_{k \in \mathcal{K}} \omega_k\, f(y_k) \\
    \text{s.t.} \quad
        & h_{ij}(y_k) \leq 0,
          && (i, j) : W_{ij}^l(\hat{d}_k)
          = \text{True},\; k \in \mathcal{K} \\
        & g(y_k) \leq 0,
          && k \in \mathcal{K} \\
        & y_k \in \mathcal{Y},
          && k \in \mathcal{K},
\end{aligned}
\end{equation}
where $\omega_k$ are the quadrature weights from the transcription.
When there is no coupling across support points (no global
constraints linking different~$k$), \eqref{eq:loa_nlp} decomposes
into $|\mathcal{K}|$ independent subproblems. With coupling (e.g.,
derivative or measure constraints linking support points), it
remains a single NLP. If
$(S'_\infty)$ is infeasible for a given trajectory, a feasibility
restoration subproblem introduces a scalar slack~$u$ and solves
$\min\, u$ subject to $h_{ij}(y_k) \leq u$, $g(y_k) \leq u$. The
resulting~$y^l$ is still used for OA cut generation
\cite{viswanathan1990penalty}.

\paragraph{OA master.}
After $L$ NLP solves, let $K'_{L,j}$ denote the linearization set
for disjunction~$j$: the iterations $l$ in which $W_{ij}^l =
\text{True}$ for some~$i$, following the notation of
\cite{turkay1996logic}. In what follows, $\nabla$ denotes the
gradient with respect to the decision variables, so that
$\nabla f(y_k^l)$ is the gradient of $f$ evaluated at the NLP
solution point $y_k^l$ and
$f(y_k^l) + \nabla f(y_k^l)^\top (y_k - y_k^l)$ is the
corresponding first-order linearization. The infinite OA master is
\begin{subequations}\label{eq:loa_master}
\begin{align}
    \min_{y_k,\,w_{ij,k},\,\alpha_{oa},\,\sigma_{ij,k}^l}
    \quad & \alpha_{oa}
      + \mu \sum_{l, i, j, k} \sigma_{ij,k}^l
    \label{eq:loa_master_obj} \\
    \text{s.t.} \quad
    & \alpha_{oa} \geq \sum_{k \in \mathcal{K}} \omega_k
      \bigl[ f(y_k^l) + \nabla f(y_k^l)^\top
      (y_k - y_k^l) \bigr],
      && l = 1, \ldots, L
      \label{eq:loa_oa_obj} \\
    & g(y_k^l) + \nabla g(y_k^l)^\top
      (y_k - y_k^l) \leq 0,
      && l = 1, \ldots, L,\;
      k \in \mathcal{K}
      \label{eq:loa_oa_global} \\
    & h_{ij}(y_k^l) + \nabla h_{ij}(y_k^l)^\top
      (y_k - y_k^l) \leq M(1 - w_{ij,k})
      + \sigma_{ij,k}^l,
      && l \in K'_{L,j},\;
      i \in \mathcal{I}_j,\; j \in \mathcal{J},\;
      k \in \mathcal{K}
      \label{eq:loa_oa_disj} \\
    & \sum_{i \in \mathcal{I}_j} w_{ij,k} = 1,
      && j \in \mathcal{J},\;
      k \in \mathcal{K},
      \label{eq:loa_exactly1} \\
    & y_k \in \mathcal{Y}, \quad
      w_{ij,k} \in \{0, 1\}, \quad
      \sigma_{ij,k}^l \in [0, \bar{\sigma}], && i \in \mathcal{I}_j, j \in \mathcal{J}, k \in \mathcal{K}.
      \label{eq:loa_domains}
\end{align}
\end{subequations}
Following the augmented-penalty OA of
\cite{viswanathan1990penalty}, a non-negative slack
$\sigma_{ij,k}^l \in [0, \bar{\sigma}]$ relaxes each disjunct OA
cut~\eqref{eq:loa_oa_disj}, and the master
objective~\eqref{eq:loa_master_obj} penalizes the slack with weight
$\mu$. This keeps the master always feasible while driving slacks to
zero as cuts accumulate. Feasibility matters when the NLP at
iteration $l$ was infeasible and $y_k^l$ came from feasibility
restoration.
The structural difference from the finite
master~$\hat{M}^b_{OA}$ is that binary variables $w_{ij,k}$ exist
at each support point~$k$ independently, enabling per-support mode
selection. The three cut types scope differently over the support
grid. The objective OA cut~\eqref{eq:loa_oa_obj} is a single
constraint that aggregates over all~$k$, reflecting the measure
structure of the objective. The disjunct OA
cuts~\eqref{eq:loa_oa_disj} are pointwise: one constraint per
support~$k$, each gated by its own binary~$w_{ij,k}$, reflecting
that disjunct constraints hold for all $d \in \mathcal{D}$. The no-good
cut~\eqref{eq:loa_nogood} likewise aggregates over all supports.
Disjunct OA cuts are added only for iterations~$l \in K'_{L,j}$ in
which disjunction~$j$ had an active disjunct at support~$k$. Affine disjunct constraints
do not require OA linearization and enter the master exactly via
big-M in~\eqref{eq:loa_oa_disj}.

After exploring a trajectory $\{i_j^*(k)\}_{k,j}$, the no-good cut
\begin{equation}\label{eq:loa_nogood}
    \sum_{k \in \mathcal{K}} \sum_{j \in \mathcal{J}}
    \bigl(1 - w_{i_j^*(k),\, j,\, k}\bigr) \geq 1
\end{equation}
keeps the master from re-selecting the identical binary trajectory.
The number of distinct trajectories per disjunction is
$|\mathcal{I}_j|^{|\mathcal{K}|}$, so
\eqref{eq:loa_nogood} is weaker per cut than its finite
counterpart, but the OA
cuts~\eqref{eq:loa_oa_obj}--\eqref{eq:loa_oa_disj} remain the
primary driver of convergence.

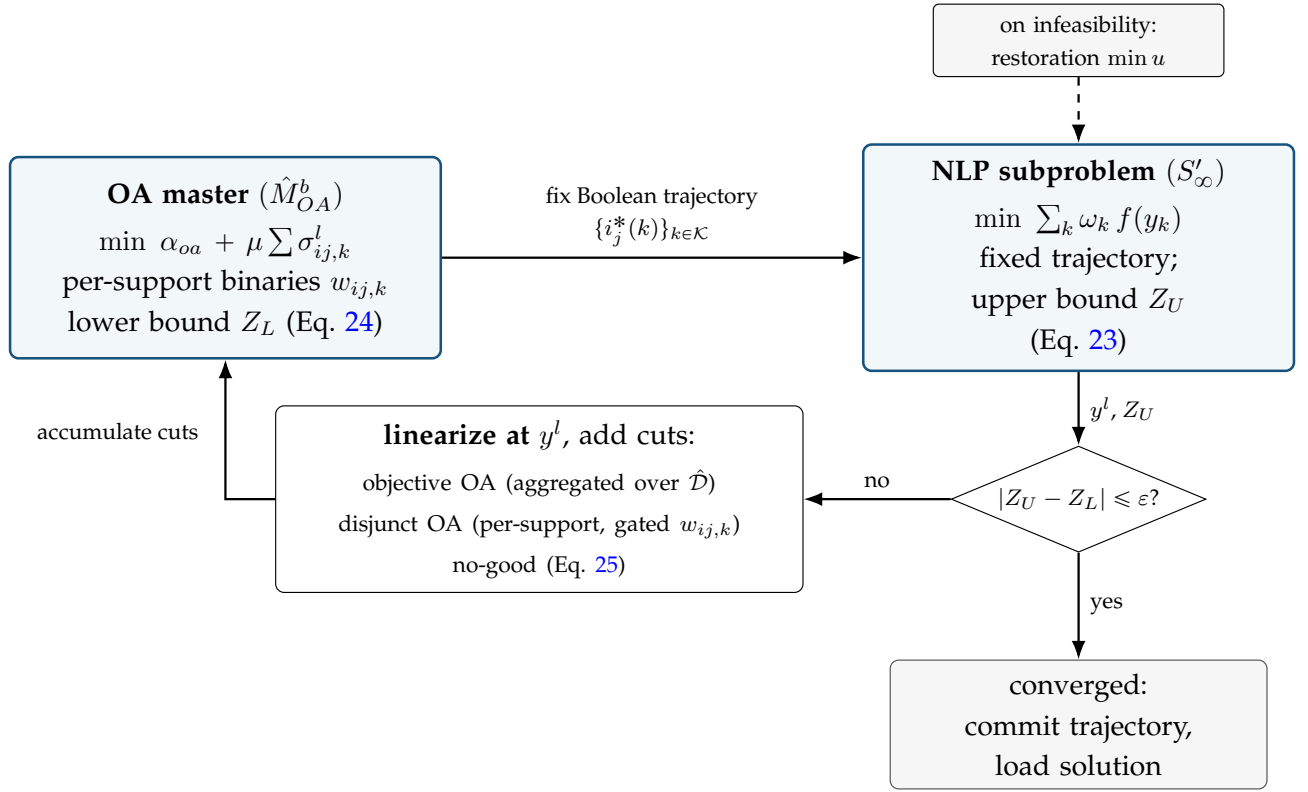
\begin{figure}[!htb]
\centering
\resizebox{\textwidth}{!}{%
\definecolor{gaaccent}{RGB}{31,119,180}%
\hyphenpenalty=10000\exhyphenpenalty=10000\relax
\begin{tikzpicture}[
    loabig/.style={draw=gaaccent!70!black, line width=0.9pt,
      rounded corners=3pt, text width=5.0cm, minimum height=2.5cm,
      fill=gaaccent!6, align=center, font=\small, inner sep=5pt},
    loawide/.style={draw, rounded corners=3pt, text width=6.2cm,
      minimum height=2.1cm, fill=white, align=center, font=\small,
      inner sep=5pt},
    loasmall/.style={draw, rounded corners=2pt, text width=3.4cm,
      minimum height=0.9cm, fill=gray!7, align=center,
      font=\scriptsize, inner sep=3pt},
    loadiamond/.style={draw, diamond, aspect=2.4,
      inner sep=2pt, font=\scriptsize},
    donebox/.style={draw=black!70, rounded corners=3pt,
      text width=4.4cm, minimum height=1.6cm, fill=gray!7,
      align=center, font=\small, inner sep=4pt}]

  \node[loabig] (master) at (0, 3.7)
    {\textbf{OA master} $(\hat{M}^b_{OA})$\\[2pt]
     $\min\ \alpha_{oa} + \mu\sum \sigma^l_{ij,k}$\\
     per-support binaries $w_{ij,k}$\\
     lower bound $Z_L$ (Eq.~\ref{eq:loa_master})};
  \node[loabig] (nlp) at (10.6, 3.7)
    {\textbf{NLP subproblem} $(S'_\infty)$\\[2pt]
     $\min\ \sum_{k} \omega_k\, f(y_k)$\\
     fixed trajectory; upper bound $Z_U$\\
     (Eq.~\ref{eq:loa_nlp})};

  \node[loasmall] (restore) at (10.6, 6.4)
    {on infeasibility: restoration $\min u$};
  \draw[flowarr, dashed] (restore.south) -- (nlp.north);

  \node[loadiamond] (gap) at (10.6, 0.7)
    {$|Z_U - Z_L| \le \varepsilon$?};
  \node[donebox] (done) at (10.6, -2.1)
    {converged:\\ commit trajectory,\\ load solution};

  \node[loawide] (cuts) at (3.9, 0.7)
    {\textbf{linearize at} $y^l$, add cuts:\\[2pt]
     {\scriptsize objective OA (aggregated over $\hat{\mathcal{D}}$)}\\
     {\scriptsize disjunct OA (per-support, gated $w_{ij,k}$)}\\
     {\scriptsize no-good (Eq.~\ref{eq:loa_nogood})}};

  \draw[flowarr] (master.east) --
    node[above, font=\scriptsize, align=center]
      {fix Boolean trajectory\\ $\{i^*_j(k)\}_{k \in \mathcal{K}}$}
    (nlp.west);
  \draw[flowarr] (nlp.south) --
    node[right, font=\scriptsize] {$y^l$, $Z_U$} (gap.north);
  \draw[flowarr] (gap.south) --
    node[right, font=\scriptsize] {yes} (done.north);
  \draw[flowarr] (gap.west) --
    node[above, font=\scriptsize] {no} (cuts.east);
  \draw[flowarr] (cuts.west) -| (master.south);
  \node[font=\scriptsize, anchor=east] at (-0.2, 1.55)
    {accumulate cuts};

\end{tikzpicture}%
}
\caption{Per-iteration mechanics of the logic-based outer
approximation loop. The OA master selects a Boolean trajectory
$\{i^*_j(k)\}$ over the support grid and supplies the lower bound
$Z_L$; the NLP subproblem solves at that fixed trajectory, yielding
the point $y^l$ and the upper bound $Z_U$ (via feasibility
restoration when the NLP is infeasible). Linearizing at $y^l$ appends
the objective, disjunct, and no-good cuts to the master, and the loop
repeats until $|Z_U - Z_L| \le \varepsilon$.}
\label{fig:loa_loop}
\end{figure}

Both subproblems are posed on the support grid $\hat{\mathcal{D}}$, in the
same manner as the cutting plane method of
Section~\ref{sec:cp}. The NLP subproblem~\eqref{eq:loa_nlp}
contains the active disjunct constraints, global constraints,
and objective of the original InfiniteGDP model. The
master~\eqref{eq:loa_master} is the big-M reformulation of the
transcribed model with the nonlinear objective replaced by
$\alpha_{oa}$, accumulating
cuts~\eqref{eq:loa_oa_obj}--\eqref{eq:loa_nogood} as iterations
proceed. Figure~\ref{fig:loa_loop} summarizes the resulting
NLP--master loop. For convex $f$, $g$, $h_{ij}$, the algorithm converges finitely
(Theorem~1 of \cite{turkay1996logic}). The proof extends to the
transcribed infinite case since it reduces to a finite MINLP after
discretization.

\subsection{Implementation}\label{sec:implementation}

The three generalized reformulations, including the MBM-GP variant,
are implemented in \texttt{DisjunctiveProgramming.jl} and
\texttt{InfiniteDisjunctiveProgramming.jl} (the
InfiniteOpt extension) as subtypes of
\texttt{AbstractReformulationMethod}, the same dispatch hierarchy that
already carries the \texttt{BigM}, \texttt{hull}, and \texttt{Indicator}
reformulations of \cite{gondosiswanto2025infiniteopt}. A model is built and stored as an \texttt{InfiniteGDPModel} using the API discussed in \cite{gondosiswanto2025infiniteopt}. The
reformulation is then simply chosen and/or updated at each solve through the
\texttt{gdp\_method} keyword on \texttt{optimize!}, making it easy to switch between
\texttt{BigM()}, \texttt{MBM()}, \texttt{PSplit(P)},
\texttt{CuttingPlanes()}, and \texttt{Hull()}.

Since LOA behaves like a solver rather than a reformulation, it is implemented as an optimizer under the \texttt{DisjunctiveAlgorithms.jl} package. Hence, \texttt{DisjunctiveAlgorithms.jl} acts as a solver that is able to directly solve GDP formulations using algorithms such as LOA. Selecting LOA
therefore is done by using \texttt{DisjunctiveAlgorithms.Optimizer} instead of another optimizer.

Code Snippet~\ref{lst:background} demonstrates the syntax on a small
three-mode heating problem. Over a 24-hour horizon, the heating
power $P(t) \in [0, 5]$ operates in one of three modes: high
($P \geq 4$), medium ($2 \leq P \leq 3$), or off ($P = 0$). The
mode is chosen independently at every support point through the
\texttt{Exactly(1)} cardinality constraint, and the objective
minimizes the integrated power. To learn more, please visit \url{https://github.com/infiniteopt/DisjunctiveProgramming.jl}.

\begin{figure}[!htp]
\begin{minipage}[t]{\linewidth}
\begin{scriptsize}
\lstset{language=Julia, breaklines = true}
\begin{lstlisting}[label=lst:background, caption={Simple three-mode
heating model in
\texttt{InfiniteDisjunctiveProgramming.jl}; only the
\texttt{gdp\_method} argument changes to select a solution method.}]
using DisjunctiveProgramming, InfiniteOpt, Gurobi

# Define model and variables
model = InfiniteGDPModel(Gurobi.Optimizer)
@infinite_parameter(model, t in [0, 24], num_supports = 24)
@variable(model, 0 <= P <= 5, Infinite(t))
@variable(model, W[1:3], InfiniteLogical(t))

# Define the infinite disjunction
@constraint(model, P >= 4.0,        Disjunct(W[1]))   # high
@constraint(model, 2.0 <= P <= 3.0, Disjunct(W[2]))   # medium
@constraint(model, P == 0,          Disjunct(W[3]))   # off
@disjunction(model, W)
@constraint(model, W in Exactly(1))

# Define the objective, reformulate, and solve
@objective(model, Min, integral(P, t))
optimize!(model, gdp_method = MBM())   # or BigM(), PSplit(2), Hull(), CuttingPlanes(), Indicator(), ...
\end{lstlisting}
\end{scriptsize}
\end{minipage}
\end{figure}

\subsection{Method Selection}
\label{sec:guidance}

Because the methods are interchangeable at solve time, the
practical question is which to select for a given model. The choice is governed by three factors: the
tightness of the continuous relaxation, the cost split between the
reformulation (pre-solve) phase and the solver, and the convexity of the
disjunct constraints. Table~\ref{tab:method_guidance} summarizes the
resulting recommendations, which the benchmarks of
Section~\ref{sec:casestudies} empirically demonstrate. Note that the performance of each approach is problem-dependent, and so this section only seeks to provide general recommendations that serve as a rule of thumb for practitioners. 

\begin{table}[H]
\centering
\caption{Relaxation/cost profile of the six solution methods and when
each might be preferred. MBM-GP is listed as an implementation variant of
MBM rather than a separate method.}
\label{tab:method_guidance}
\small
\begin{tabular}{@{}l p{0.46\textwidth} p{0.30\textwidth}@{}}
\toprule
Method & Relaxation / cost profile & Prefer when \\
\midrule
big-M & loosest relaxation; cheapest to build; stays linear
  & a fast pre-solve matters more than the gap \\
MBM & tighter, trajectory-adaptive; still linear; per-pair
  $M$ cost up front & solve cost dominates; a tight single-tree model
  is wanted \\
\quad MBM-GP & MBM tightness at solved supports, GP upper bound
  elsewhere; a fraction of MBM's pre-solve cost & MBM is
  wanted but its per-pair $M$ pre-solve cost is prohibitive \\
P-split & tunable big-M$\leftrightarrow$hull via $P$; larger model
  & a tunable relaxation; convex (McCormick) partitions \\
hull & tightest convex relaxation; variable doubling, perspective NLP
  & a strong relaxation justifies the larger model \\
CP & hull-strength without doubling; iterative
  & \emph{convex} GDP only (separation must be convex) \\
LOA & OA for nonlinear convex disjuncts; avoids a nonlinear MIP
  & convex GDP with nonlinear disjunct constraints \\
\bottomrule
\end{tabular}
\end{table}

Big-M is the cheapest reformulation to build and stays linear when the
disjuncts are, but its relaxation is the loosest, so it leaves a large
gap when the solver is time-limited. Multiple big-M removes most of that
looseness while remaining a single-tree linear model: the per-pair
$M_{ii'}(d)$ adapt to the trajectory, giving the tightest linearizing
relaxation, at the cost of an up-front per-pair $M$ computation that
shifts work into the pre-solve phase. This makes MBM a good balanced
choice for nonconvex InfiniteGDPs: it trades a higher up-front
reformulation cost for a tighter continuous relaxation that can accelerate convergence, so the balance pays off most when the solve would otherwise
dominate, as on a time-limited spatial branch-and-bound. The pre-solve investment is also a one-time
cost. In applications that repeatedly re-solve the same disjunctive
structure with updated global parameters, such as model predictive control or
multi-objective optimization, the computed $M_{ii'}(d)$ functions are
reused across re-solves, so the reformulation cost amortizes while
the tighter continuous relaxation pays off at every solve. When the up-front cost
is prohibitive, the MBM-GP variant of Section~\ref{sec:mbmgp}
solves only a fraction of the per-support subproblems and is the
cheaper option at fine discretizations. Hull yields the tightest convex relaxation per
disjunction but doubles the variables and introduces perspective
nonlinearities, so it builds quickly yet leaves a larger transcribed
model for the solver. P-split interpolates between big-M and hull through
the partition size $P$, offering a tunable position on the spectrum; it
requires convex partition expressions, so for bilinear couplings the
nonconvex products are enveloped with McCormick relaxations to keep the
partitioned formulation a valid relaxation.

Cutting planes and LOA are decomposition algorithms with a convexity
requirement. CP recovers hull-strength tightness without variable
doubling, but its separation subproblem inherits the convexity of the
original disjunct constraints: a valid separating cut is only guaranteed
when the separation problem is convex. CP is therefore sound only on convex
GDP; on a nonconvex model it can generate cuts that remove feasible,
even optimal, points (demonstrated directly in
Section~\ref{sec:biodiesel}). LOA targets convex GDPs whose disjunct
constraints are nonlinear, where a direct hull would require a nonlinear
MIP solver; its augmented-penalty form tolerates mild nonconvexity in
practice but carries no global guarantee. In summary, MBM is a
good balanced choice for nonconvex problems, trading a higher
reformulation cost for a tighter, still-linear master, and MBM-GP
recovers most of that tightness when the pre-solve budget is small.
Hull and
P-split become useful when a stronger relaxation justifies a larger model. Big-M is practical when weak relaxations do not strongly influence solution time. One should generally restrict the use of CP
and LOA to the convex regime for which they are designed.

\section{Case studies}\label{sec:casestudies}

This section evaluates the methods on three case studies. For each, the
formulation is presented first and the methods are then benchmarked
against the guidance of Section~\ref{sec:guidance}. Two are natively infinite-dimensional
dynamic-optimization problems with nonconvex bilinear couplings: a
minimum lap-time control problem with four-way powertrain mode selection
parameterized over arc length (Section~\ref{sec:f1}), and a 24-hour
biodiesel production schedule with a temperature-dependent heat-pump
coefficient of performance parameterized over time
(Section~\ref{sec:biodiesel}). The third is an event-constrained
IEEE 14-bus capacity design problem (Section~\ref{sec:powergrid}) drawn
from \cite{ovalle2025event} which is a linear stochastic program that yields a convex InfiniteGDP. Unless noted otherwise, the F1
problem is discretized at $N = 100$ support points, the biodiesel
problem at $N = 300$, and the powergrid problem at $N = 200$, and all
reformulations target a $5\%$ relative MIP gap. A cross-study summary
follows in Section~\ref{sec:summary}. For solvers, \texttt{Gurobi}~v12.0.2
is used for all MIP and MIQCP solves and the LOA subproblems, and
\texttt{Ipopt} (via \texttt{Ipopt.jl}~v1.15.0) is used for the local NLP
separation solves on the nonconvex biodiesel study. Other package versions
include \texttt{JuMP}~v1.31.1, \texttt{InfiniteOpt.jl}~v0.6.3,
\texttt{DisjunctiveProgramming.jl} v0.7.0, and \texttt{AbstractGPs.jl}~v0.5.24,
running on Julia~1.11.5. The results are collected using a Windows~11
machine with an Intel(R) i9-12900K CPU @ 3.20~GHz and 64~GB of memory.
The source code for
all case studies is available at \url{https://github.com/PULSI-Opt/InfGDP-case-studies}.

\subsection{Minimum lap-time with powertrain mode selection}
\label{sec:f1}

This case study considers the minimum-time traversal of a closed racing circuit by
a hybrid Formula~1 vehicle operating under the 2026 FIA Power Unit
Regulations. Arc length $s \in [0, L]$ along the track centerline
parameterizes the vehicle state, with $L$ the total circuit length. Track geometry enters through the curvature
$\kappa(s)$, a known continuous function over the arc-length domain,
and the corresponding friction-limited cornering
speed $\bar{v}(\kappa(s))$, which accounts for aerodynamic downforce.
At every arc length $s$ the powertrain operates in exactly one of
four modes $i \in \mathcal{I} = \{\text{clip}, \text{ice},
\text{boost}, \text{regen}\}$: partial throttle (clipping), full
internal combustion engine, ICE plus MGU-K deployment, and
regenerative braking. The four-mode classification follows the
time-optimal hybrid race-car policy of \cite{salazar2017hybrid},
where these modes arise as Pontryagin switching cases of a
continuous policy; they are recast here as an explicit arc-length
disjunction so that the active mode is a Boolean decision at
every support point. Each mode constrains the motor
force $F(s)$,
the electrical battery power $P_b(s)$, the motor power $P_m(s) =
F(s)v(s)$, and the effective propulsive force $F_e(s)$ delivered to
the tires through the drivetrain.

The objective is the lap time
\begin{equation}\label{eq:f1_obj}
    \min \;\; t(L),
\end{equation}
which accumulates along the circuit through the kinematic
relation
\begin{equation}\label{eq:f1_time}
    v(s)\,\tfrac{dt}{ds} = 1, \quad s \in [0, L].
\end{equation}

Vehicle dynamics follow Newton's second law along the track
centerline, expressed in arc-length form as in
\cite{perantoni2014f1},
\begin{equation}\label{eq:f1_long}
    m\,v(s)\,\tfrac{dv}{ds}
        = F_e(s) - \tfrac{1}{2}\rho C_d A\, v(s)^2 - f_r m g,
    \quad s \in [0, L],
\end{equation}
with vehicle mass $m$, drag area coefficient $C_d A$, and
rolling resistance $f_r$. The state of charge drains with
battery power, in the same arc-length form used by
\cite{limebeer2014ers} for F1 energy recovery systems,
\begin{equation}\label{eq:f1_soc}
    E_b\,v(s)\,\tfrac{d\,\mathrm{SOC}}{ds} = -P_b(s),
    \quad s \in [0, L],
\end{equation}
where $E_b$ is the battery capacity, and motor power is
mechanically tied to motor force by
\begin{equation}\label{eq:f1_pm}
    P_m(s) = F(s)\,v(s),
    \quad s \in [0, L].
\end{equation}

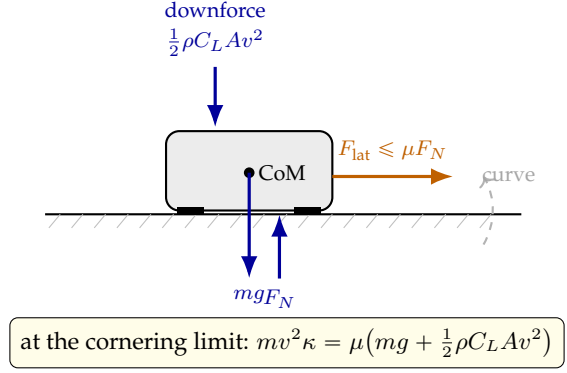
\begin{figure}[H]
\centering
\begin{tikzpicture}[
    vforce/.style={-Latex, very thick, blue!60!black},
    lforce/.style={-Latex, very thick, orange!75!black},
    aux/.style={dashed, gray!60}
]
  \draw[thick] (-2.7, 0) -- (3.6, 0);
  \foreach \x in {-2.5,-2.2,...,3.5} {
    \draw[gray!50, thin] (\x, 0) -- ({\x-0.15}, -0.18);
  }

  \draw[thick, rounded corners=5pt, fill=gray!15]
    (-1.1, 0.05) rectangle (1.1, 1.1);

  \fill[black] (-0.95, 0) rectangle (-0.6, 0.1);
  \fill[black] (0.6, 0) rectangle (0.95, 0.1);

  \fill (0, 0.55) circle (2pt);
  \node[font=\scriptsize, anchor=west, inner sep=2pt]
    at (0.05, 0.55) {CoM};

  \draw[vforce] (0, 0.55) -- (0, -0.85)
    node[at end, below, font=\scriptsize] {$m g$};

  \draw[vforce] (-0.45, 1.95) -- (-0.45, 1.15)
    node[at start, above, font=\scriptsize, align=center]
    {downforce\\ $\tfrac{1}{2}\rho C_L A v^2$};

  \draw[vforce] (0.4, -0.85) -- (0.4, 0)
    node[at start, below, font=\scriptsize] {$F_N$};

  \draw[lforce] (1.1, 0.5) -- (2.7, 0.5)
    node[midway, above=2pt, font=\scriptsize]
    {$F_{\text{lat}} \leq \mu F_N$};

  \draw[aux, ->, thick] (3.1, -0.4) arc[start angle=-30,
    end angle=30, radius=0.9];
  \node[font=\scriptsize, gray!70, anchor=west]
    at (2.95, 0.5) {curve};

  \node[draw, rounded corners=3pt, fill=yellow!10,
    align=center, font=\footnotesize, inner sep=4pt]
    at (0.5, -1.7)
    {at the cornering limit:
       $m v^2 \kappa
       = \mu\bigl(m g + \tfrac{1}{2}\rho C_L A v^2\bigr)$};
\end{tikzpicture}
\caption{Free-body diagram for the closed-form cornering speed
\eqref{eq:f1_corner} (rear view). Vertical forces (blue) set the
tire normal load $F_N$; the friction-circle bound caps the lateral
tire force (orange) at $\mu F_N$, which at the cornering limit
balances the centripetal demand $m v^2 \kappa$.}
\label{fig:f1_cornering}
\end{figure}

The cornering speed comes from a rigid-body lateral force
balance with aerodynamic downforce, in closed form,
\begin{equation}\label{eq:f1_corner}
    \bar{v}(\kappa) = \sqrt{\frac{\mu g}
    {\kappa \;-\; \tfrac{\mu \rho C_L A}{2 m}}},
\end{equation}
obtained by setting the required centripetal force
$m v^2 \kappa$ equal to the friction-limited tire force
$\mu(m g + \tfrac{1}{2}\rho C_L A v^2)$ on a load that combines
static weight and downforce, with downforce area coefficient
$C_L A$. This point-mass closed form replaces the full Pacejka tire
model with load transfer of \cite{perantoni2014f1}, and
$\bar{v}(\kappa(s))$ enters the model as a known continuous function
over the arc-length domain, evaluated on the support grid.
Figure~\ref{fig:f1_cornering} sketches the corresponding
free-body diagram.

The state envelopes combine cornering and SOC limits,
\begin{align}
    v_{\min} \leq v(s) &\leq \min\{v_{\max},
        \bar{v}(\kappa(s))\},
    && s \in [0, L], \label{eq:f1_vbounds} \\
    \mathrm{SOC}_{\min} \leq \mathrm{SOC}(s)
        &\leq \mathrm{SOC}_{\max},
    && s \in [0, L]. \label{eq:f1_socbounds}
\end{align}
A longitudinal acceleration envelope caps how fast the speed can
change, keeping the trajectory within the available grip and
powertrain limits,
\begin{equation}\label{eq:f1_accel}
    a_{\min} \leq v(s)\,\tfrac{dv}{ds} \leq a_{\max},
    \quad s \in [0, L],
\end{equation}
where $v\,\tfrac{dv}{ds} = \tfrac{dv}{dt}$ is the longitudinal
acceleration. The bounds are asymmetric. Aerodynamic downforce grows
the tire normal load with speed, so braking sustains far higher
deceleration than acceleration, which is limited by traction and
engine power. The envelope bounds
longitudinal grip and complements the lateral cornering limit in
\eqref{eq:f1_corner}.

At every arc length the powertrain selects exactly one of four
modes through the disjunction
\begin{equation}\label{eq:f1_disj}
    \bigvee_{i \in \mathcal{I}}
      \begin{bmatrix}
        W_i(s) \\[4pt]
        \underline{F}_i \leq F(s) \leq \bar{F}_i \\
        \underline{P}_m^i \leq P_m(s) \leq \bar{P}_m^i \\
        F_e(s) = \eta_{\mathrm{dt}}^i\, F(s) \\
        \underline{P}_b^i \leq P_b(s) \leq \bar{P}_b^i(s) \\
        g_i(P_m(s), P_b(s)) \leq 0
      \end{bmatrix},
    \quad s \in [0, L].
\end{equation}
Drivetrain efficiency $\eta_{\mathrm{dt}}^i$ equals
$\eta_{\mathrm{dt}}$ for the propulsion modes and $1$ for the
regenerative mode. The mode-dependent coupling $g_i$ links motor
and battery power: for clipping and ICE the battery is inactive
($P_b = 0$); for boost it enforces the band $P_m -
\bar{P}_{\text{ice}} \leq \eta_{\mathrm{mguk}} P_b \leq P_m -
\underline{P}_{\text{ice}}$; for regeneration it enforces
$P_b = \eta_{\mathrm{mguk}} P_m$ with $P_m \leq 0$. Mode bounds are
drawn from the quasi-steady-state simulation of
\cite{heilmeier2019qss}, scaled to the 2026 400~kW ICE and
$\pm 350$~kW MGU-K limits. The battery bounds are constant in $s$
for all modes except boost. Its deployment ceiling follows a
lap-periodic energy-management schedule,
\begin{equation}\label{eq:f1_depcap}
    \bar{P}_b^{\text{boost}}(s)
        = \min\bigl\{\bar{P}_{\text{mguk}},\,
          P_b^{\text{cap}}(s)\bigr\},
    \qquad
    P_b^{\text{cap}}(s) = 200 + 150\cos(2\pi s / L)~\text{kW},
\end{equation}
which keeps the cap within $[50, 350]$~kW so the boost region
never empties. The schedule is not necessarily realistic, but introducing it makes the boost disjunct's constraint data explicitly position-dependent such that this more rigorously benchmarks the proposed methods.

Exactly one mode is active at every support point,
\begin{equation}\label{eq:f1_exactly}
    \mathrm{EXACTLY}\bigl(1, W_{\text{clip}}(s),
        W_{\text{ice}}(s),
        W_{\text{boost}}(s), W_{\text{regen}}(s)\bigr)
        = \text{True},
    \quad s \in [0, L],
\end{equation}
and the lap closes,
\begin{equation}\label{eq:f1_bc}
    v(0) = v(L), \quad
    \mathrm{SOC}(0) = \mathrm{SOC}(L), \quad
    t(0) = 0,
\end{equation}
with $W_i : [0, L] \mapsto \{\text{True}, \text{False}\}$ for
$i \in \mathcal{I}$.

The full formulation reads
\begin{equation*}
\begin{aligned}
\min \;\; & \text{lap time}
    && \eqref{eq:f1_obj} \\
\text{s.t.} \;\;
    & \text{kinematics and vehicle dynamics}
    && \eqref{eq:f1_time},\,
        \eqref{eq:f1_long}\text{--}\eqref{eq:f1_pm} \\
    & \text{state envelopes}
    && \eqref{eq:f1_corner}\text{--}\eqref{eq:f1_accel} \\
    & \text{mode disjunction and cardinality}
    && \eqref{eq:f1_disj},\,\eqref{eq:f1_exactly} \\
    & \text{boundary conditions}
    && \eqref{eq:f1_bc}.
\end{aligned}
\end{equation*}
The distinguishing structural feature is the spatial
parameterization of the logical variables: $W_i(s)$ is infinite
and the active mode may switch at every support point along the
circuit. Unlike the tank changeover in
\cite{gondosiswanto2025infiniteopt}, where mode selection is
constant over each operating period, here the selection is
fully spatially resolved. The continuous dynamics contain three
bilinear terms ($v \cdot dv/ds$, $v \cdot d\mathrm{SOC}/ds$, and
$F \cdot v$), which Gurobi's spatial branch-and-bound handles.

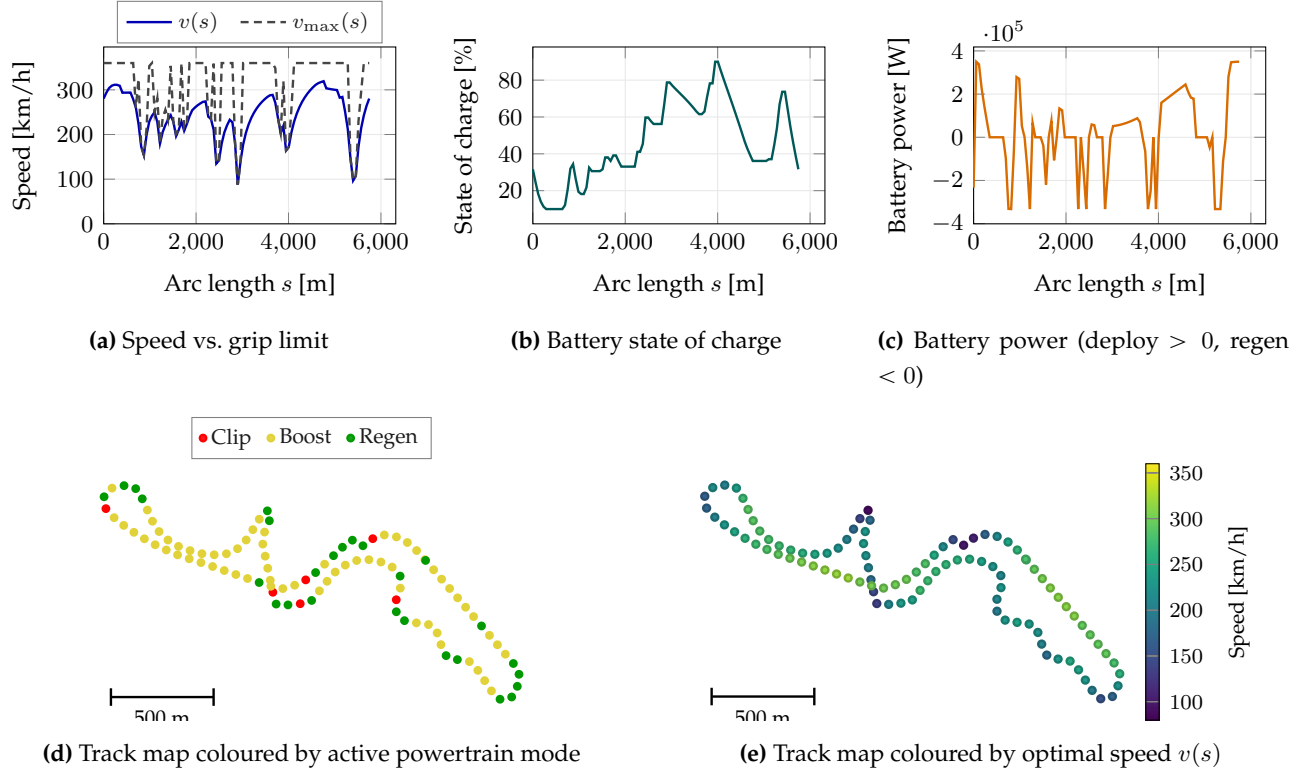
\begin{figure}[H]
\centering
\begin{subfigure}[t]{0.32\textwidth}
\centering
\begin{tikzpicture}
\begin{axis}[trajpanel, xlabel={Arc length $s$ [m]},
  ylabel={Speed [km/h]}, xmin=0, ymin=0,
  legend style={at={(0.5,1.02)}, anchor=south, legend columns=2,
  /tikz/every even column/.append style={column sep=6pt}}]
  \addplot[blue!70!black] table[col sep=comma, x=s_m, y=v_kmh]
    {data/traj_f1.csv};
  \addplot[gray!55!black, densely dashed] table[col sep=comma,
    x=s_m, y=v_max_kmh] {data/traj_f1.csv};
  \legend{$v(s)$, $v_{\max}(s)$}
\end{axis}
\end{tikzpicture}
\caption{Speed vs.\ grip limit}
\label{fig:f1_speed}
\end{subfigure}\hfill
\begin{subfigure}[t]{0.32\textwidth}
\centering
\begin{tikzpicture}
\begin{axis}[trajpanel, xlabel={Arc length $s$ [m]},
  ylabel={State of charge [\%]}, xmin=0]
  \addplot[teal!70!black] table[col sep=comma, x=s_m, y=SOC_pct]
    {data/traj_f1.csv};
\end{axis}
\end{tikzpicture}
\caption{Battery state of charge}
\label{fig:f1_soc}
\end{subfigure}\hfill
\begin{subfigure}[t]{0.32\textwidth}
\centering
\begin{tikzpicture}
\begin{axis}[trajpanel, xlabel={Arc length $s$ [m]},
  ylabel={Battery power [W]}, xmin=0]
  \addplot[orange!85!black] table[col sep=comma, x=s_m, y=P_batt_W]
    {data/traj_f1.csv};
\end{axis}
\end{tikzpicture}
\caption{Battery power (deploy $>0$, regen $<0$)}
\label{fig:f1_pbatt}
\end{subfigure}

\vspace{0.7em}

\begin{subfigure}[t]{0.48\textwidth}
\centering
\begin{tikzpicture}
\begin{axis}[trajpanel, axis equal image, hide axis,
  legend style={at={(0.5,1.02)}, anchor=south,
  legend columns=3,
  /tikz/every even column/.append style={column sep=4pt},
  font=\scriptsize}]
  \addplot[only marks, mark=*, mark size=1.2pt, red]
    table[col sep=comma, x=x_m, y=y_m]
    {data/traj_f1_track_Clip.csv};
  \addplot[only marks, mark=*, mark size=1.2pt, yellow!85!black]
    table[col sep=comma, x=x_m, y=y_m]
    {data/traj_f1_track_Boost.csv};
  \addplot[only marks, mark=*, mark size=1.2pt, green!60!black]
    table[col sep=comma, x=x_m, y=y_m]
    {data/traj_f1_track_Regen.csv};
  \legend{Clip, Boost, Regen}
  \draw[line width=0.9pt, |-|]
    (axis cs:-1500,-650) -- (axis cs:-1000,-650)
    node[midway, below, font=\scriptsize, yshift=-1pt] {500 m};
\end{axis}
\end{tikzpicture}
\caption{Track map coloured by active powertrain mode}
\label{fig:f1_track}
\end{subfigure}\hfill
\begin{subfigure}[t]{0.48\textwidth}
\centering
\begin{tikzpicture}
\begin{axis}[trajpanel, axis equal image, hide axis,
  colorbar, colorbar style={width=5pt,
    ylabel={Speed [km/h]},
    ylabel style={font=\scriptsize, yshift=-2pt},
    ytick={100,150,200,250,300,350},
    tick label style={font=\scriptsize}},
  colormap name=viridis,
  point meta min=80, point meta max=360]
  \addplot[scatter, scatter src=explicit, only marks,
           mark=*, mark size=1.3pt, mark options={draw=none}]
    table[col sep=comma, x=x_m, y=y_m, meta=v_kmh]
    {data/traj_f1.csv};
  \draw[line width=0.9pt, |-|]
    (axis cs:-1500,-650) -- (axis cs:-1000,-650)
    node[midway, below, font=\scriptsize, yshift=-1pt] {500 m};
\end{axis}
\end{tikzpicture}
\caption{Track map coloured by optimal speed $v(s)$}
\label{fig:f1_speedmap}
\end{subfigure}
\caption{Solution trajectory for the F1 lap-time problem on the
Suzuka circuit ($N = 100$, MBM reformulation, lap time
$91.61$~s). The trajectory is an illustrative solve with Gurobi's
default settings; the benchmark runs of Table~\ref{tab:summary}
use the configuration stated there and stop at the $5\%$ gap
target. (a) Optimal speed $v(s)$ and grip-limited ceiling
$v_{\max}(s)$ along the lap. (b) Battery state of charge. (c)
Battery power (deployment $>0$, regeneration $<0$). (d) Active
powertrain mode at each support point mapped onto the circuit
geometry. (e) The same support
points coloured by $v(s)$.}
\label{fig:f1_traj}
\end{figure}

This problem is solved on the Suzuka circuit
($L = 5.807$~km), discretized at $N = 100$ support points.
Figure~\ref{fig:f1_traj} shows a representative optimal solution: the
speed tracks the grip-limited ceiling through the corners, the battery
deploys and regenerates around the lap, and the active powertrain mode
switches at individual support points along the circuit, with the
ICE-only mode never selected. With three
bilinear couplings this is a nonconvex MIQCP, so the guidance of
Section~\ref{sec:guidance} predicts that MBM reaches a near-optimal lap
time fastest, hull converges more slowly through a larger model, big-M
stays far from the best-known lap time within the budget, and CP --
whose separation problem is nonconvex here -- fails to produce a usable
incumbent. LOA is likewise outside the convex regime it targets, so it
carries no valid bound on this problem. It is still run as a heuristic
to probe how far its augmented-penalty form degrades under strong
nonconvexity. The benchmark results below confirm these predictions.

\begin{figure}[H]
\centering
\begin{tikzpicture}
\begin{axis}[timebar, symbolic x coords={MBM,MBM-GP,PSplit-3},
             area legend,
             legend style={at={(0.5,1.03)}, anchor=south,
               legend columns=-1, draw=black!50,
               font=\scriptsize}]
  \addplot+[fill=blue!25, draw=blue!55!black]
    table[col sep=comma, x=method, y=reform_s]
    {data/time_f1_n100_reform.csv};
  \addplot+[fill=orange!35, draw=orange!70!black]
    table[col sep=comma, x=method, y=solver_s]
    {data/time_f1_n100_reform.csv};
  \legend{Reformulation, Solver}
\end{axis}
\end{tikzpicture}
\caption{Reformulation time (stacked, under solver time) versus solver
time on the nonconvex F1 problem ($N = 100$); values from
Table~\ref{tab:summary}. Big-M and hull are omitted because they
reach the $1200$~s time limit with no incumbent. The P-split
solver bar is the $1200$~s time limit, reached with a
$16.9\%$-gap incumbent.}
\label{fig:timebar_f1}
\end{figure}
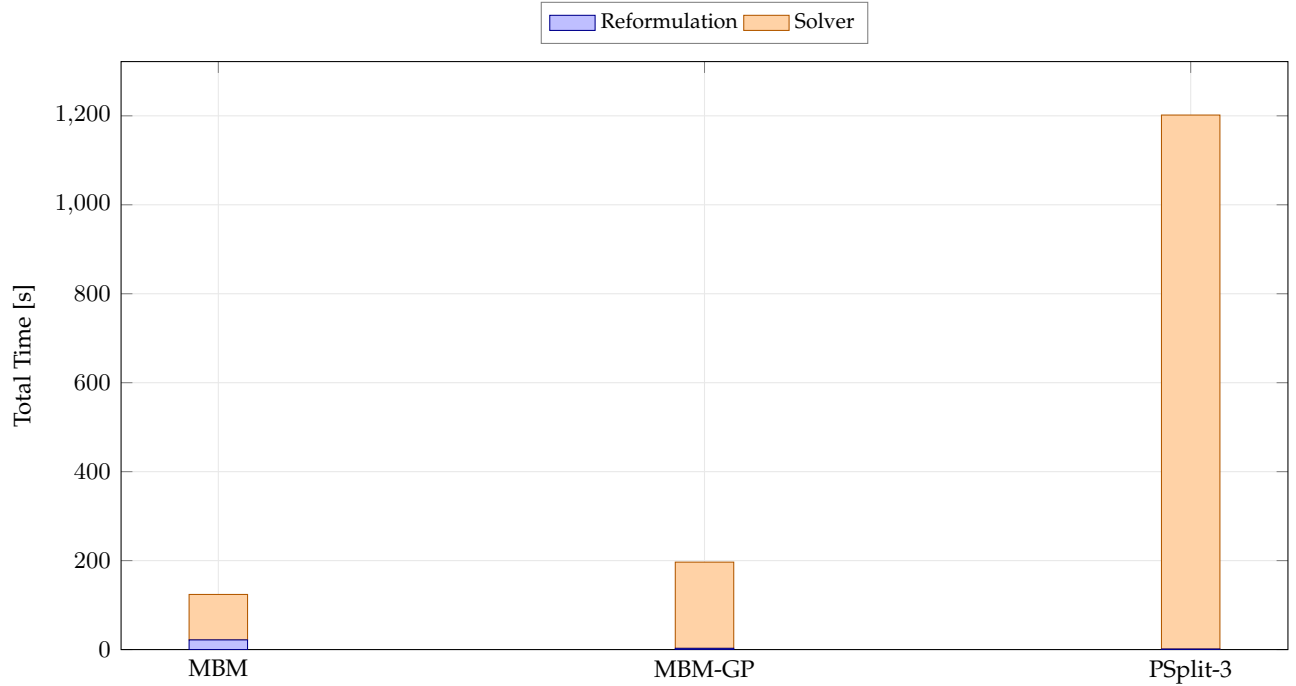

\paragraph{Results.} Progress against solver time is reported through a
solution-progress measure in $[0, 1]$ derived from each method's
running MIP gap: the gap is mapped so that reaching the $5\%$ MIP-gap
stopping tolerance corresponds to $0.95$ (the dashed reference line)
and proving the true optimum corresponds to $1.0$; a method that
plateaus below $0.95$ has not solved to tolerance. Big-M uses
$M = 10^{7}$ and P-split uses $P = 3$ with the McCormick lifting.
MBM-GP uses the adaptive point selection and $\kappa = 2.5$ of
Section~\ref{sec:mbmgp}, here and in the other two studies.
Figure~\ref{fig:solpct_f1} reports the result. The deployment
schedule \eqref{eq:f1_depcap} makes the boost disjunct's constraint
data vary along the lap, which strongly affects the relative
performance of the reformulations. Only the multiple big-M methods solve within the $1200$~s
budget. MBM crosses $0.95$ at $\approx 102$~s. MBM-GP follows at
$\approx 194$~s after a $2.7$~s build: its GP surrogate samples the
varying $M(s)$ curves instead of solving all $9600$ per-support
subproblems, which exact MBM pays $21.7$~s for. P-split finds its
first lap only at $\approx 997$~s and ends at $97.46$~s with a
$16.9\%$ gap. Big-M and hull have no incumbent when the
time limit is reached and are therefore omitted from
Figures~\ref{fig:solpct_f1} and~\ref{fig:timebar_f1}. Hence, for
this case study, feasible solutions are only obtained with the
intermediate reformulations rather than the big-M and hull
endpoints. CP is absent because, as
predicted for a nonconvex separation problem, it produced no valid
incumbent. LOA runs as a heuristic: it exhausts its $1200$~s
budget in set-covering initialization without producing a bound, so
its curve plots the incumbent's distance to the best lap time found,
reaching $164.29$~s ($78\%$ off, Table~\ref{tab:summary}) at
$\approx 251$~s with no further improvement.
Figure~\ref{fig:timebar_f1} splits the total runtime into
its reformulation (pre-solve) and solver components, confirming that
MBM's advantage comes from a tight master that the solver closes
quickly rather than from a low reformulation cost.

\begin{figure}[H]
\centering
\begin{tikzpicture}
\begin{axis}[solpanel, xmin=0, xmax=1250]
  \addplot[const plot, mark=none, blue!70!black, dashed] table
    [col sep=comma, x=time_s, y=sol_frac] {data/solpct_f1_mbm.csv};
  \addplot[const plot, mark=none, green!60!black, densely dotted] table
    [col sep=comma, x=time_s, y=sol_frac] {data/solpct_f1_mbmgp.csv};
  \addplot[const plot, mark=none, teal!70!black] table
    [col sep=comma, x=time_s, y=sol_frac]
    {data/solpct_f1_psplit.csv};
  \addplot[const plot, mark=none, violet!80!black, dashdotted] table
    [col sep=comma, x=time_s, y=sol_frac]
    {data/solpct_f1_loa.csv};
  \legend{MBM, MBM-GP, P-split, LOA$^{*}$}
  \draw[gray!55, dashed, line width=0.7pt]
    (axis cs:0,0.95) -- (axis cs:1250,0.95);
  \draw[gray!40, densely dotted, line width=0.6pt]
    (axis cs:1200,0) -- (axis cs:1200,1.03);
  \node[anchor=north west, font=\scriptsize, gray!60!black]
    at (axis cs:60,0.945) {5\% MIP-Gap Target};
\end{axis}
\end{tikzpicture}
\caption{Solution progress versus solve time on the nonconvex F1
problem ($N = 100$) under the $1200$~s wall budget (dotted vertical
line); the dashed line marks the $5\%$ MIP-gap target ($0.95$).
Big-M and hull are omitted because they find no incumbent within
the time limit; P-split's first incumbent arrives at
$\approx 997$~s. LOA$^{*}$ produces no bound on this problem, so its
curve plots the incumbent's distance to the best known lap time.}
\label{fig:solpct_f1}
\end{figure}
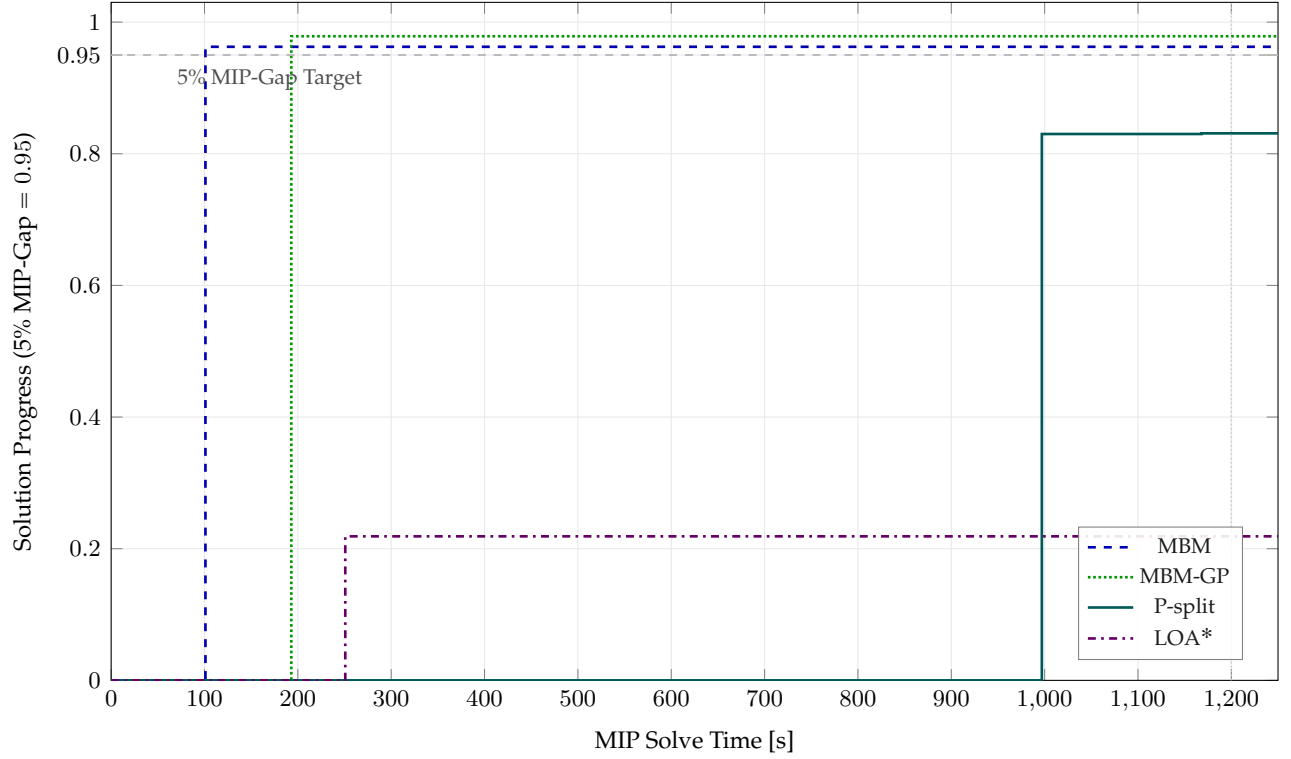

\subsection{Temperature-coupled biodiesel production scheduling}
\label{sec:biodiesel}

\begin{figure}[H]
\centering
\includegraphics[width=0.6\textwidth]{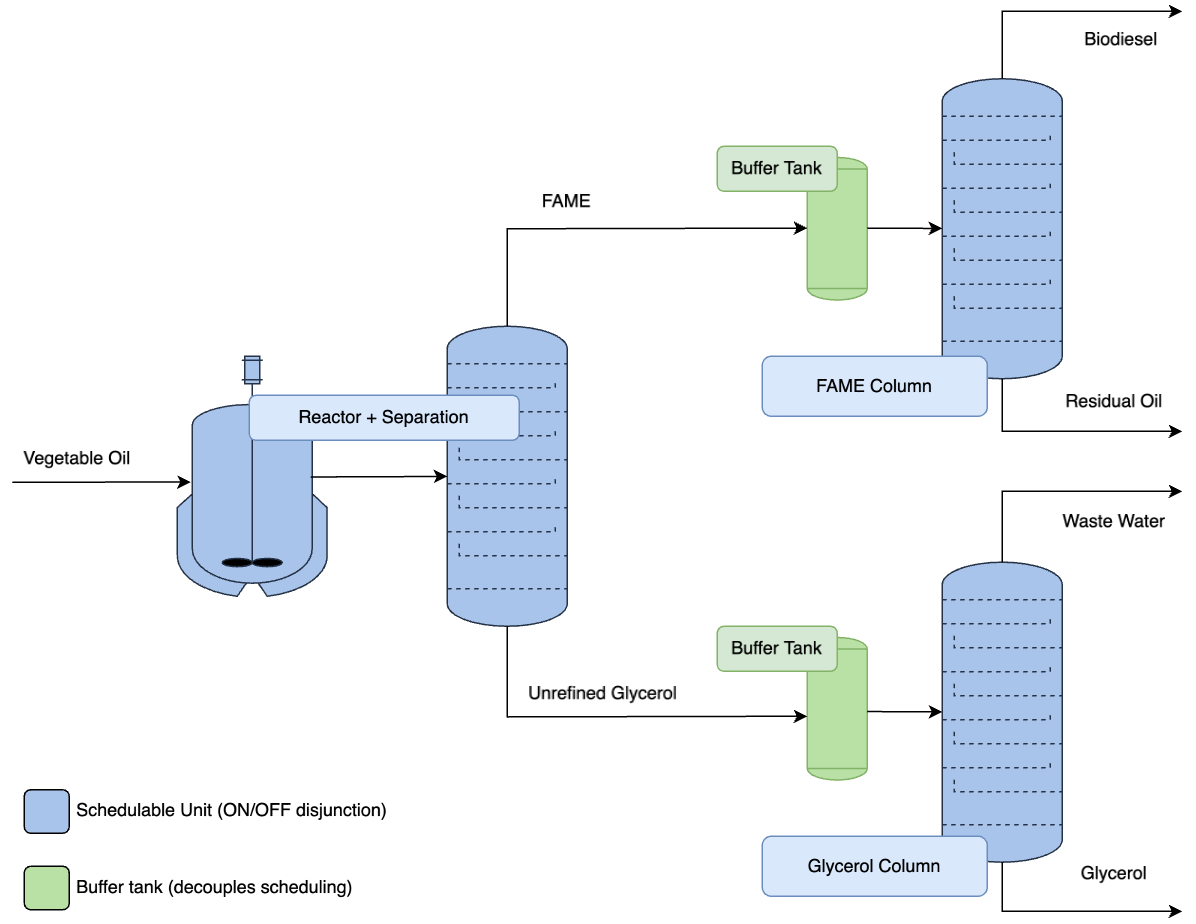}
\caption{Process flowsheet for the electrified biodiesel facility.
Each schedulable unit (blue) carries an ON/OFF disjunction over the
24-hour horizon; buffer tanks (green) decouple the unit schedules.}
\label{fig:biodiesel_flowsheet}
\end{figure}

This case study considers the 24-hour scheduling of a fully electrified biodiesel
facility \cite{biodiesel} with three process units
(Figure~\ref{fig:biodiesel_flowsheet}): a
transesterification reactor, a FAME distillation column, and a
glycerol purification column. Each unit $u \in \mathcal{U} =
\{\mathrm{rx}, \mathrm{f}, \mathrm{g}\}$ can independently switch on
or off over the horizon $t \in [0, 24]$~h in response to a
time-varying electricity price $c_e(t)$ and scenario-weighted product
prices $\bar{p}_f(t), \bar{p}_g(t)$. A heat pump drives the glycerol
column reboiler, with a coefficient of performance (COP) that
degrades with reboiler temperature. This creates a compounding
electrical penalty that couples the temperature setpoint
$T_g(t) \in [T_g^L, T_g^U]$ to both the heat duty and the electrical
power demand of the column. The on/off scheduling structure and the
temperature-dependent glycerol heat-pump COP are introduced here
as additions to the continuous-control formulation of
\cite{biodiesel}, which treats all units as always on and uses a
constant COP for the glycerol cascade. Throughout, the formulation
adopts linear flow-driven surrogates of the unit models in
\cite{biodiesel}, trading kinetic and
stage-by-stage fidelity for tractability at the scheduling
level.

Linear calibrations specify the temperature--purity and
temperature--COP couplings, with the COP form approximating the
near-linear COP--temperature trend reported for industrial
high-temperature heat pumps in \cite{arpagaus2018hthp},
\begin{align}
    \pi_g(T_g) &= \pi^{\min} +
      \frac{\pi^{\max} - \pi^{\min}}{T_g^U - T_g^L}
      (T_g - T_g^L), \label{eq:bio_purity} \\
    \mathrm{COP}_g(T_g) &= \alpha - \beta\,(T_g - T_g^L),
      \label{eq:bio_cop}
\end{align}
with $\pi^{\min} = 0.95$, $\pi^{\max} = 0.995$, $\alpha = 3.75$,
$\beta = 0.02425$~K$^{-1}$, and $T_g \in [350, 430]$~K. The
glycerol reboiler heat duty is modeled as a linear function of
flow with a temperature--flow bilinear correction,
\begin{equation}\label{eq:bio_Qreb}
    Q_g^{\text{reb}}(t) = q_0 + q_1 F_g(t)
      + \kappa_T\,(T_g(t) - T_g^L)\,F_g(t),
\end{equation}
the resulting electrical power follows from the
temperature-dependent COP,
\begin{equation}\label{eq:bio_Preb}
    P_g^{\text{reb}}(t)\cdot\mathrm{COP}_g(T_g(t))
        = Q_g^{\text{reb}}(t),
\end{equation}
and the total glycerol column power sums the reboiler and
condenser duties,
\begin{equation}\label{eq:bio_Pg}
    P_g(t) = P_g^{\text{reb}}(t)
        + Q_g^{\text{cond}}(t)/\mathrm{COP}^{\text{cond}}.
\end{equation}
The two bilinear terms $\kappa_T (T_g - T_g^L) F_g$ in
\eqref{eq:bio_Qreb} and $P_g^{\text{reb}}\cdot T_g$ in
\eqref{eq:bio_Preb} are the only nonconvex couplings in the
formulation. Together they produce a power amplification factor
of $\sim\!2.7$ between the minimum and maximum reboiler
temperatures at full throughput.

Let $W_u(t)$ denote the infinite logical variable indicating
that unit $u$ is active at time $t$. The economic objective
trades electricity cost and switching cost against product
revenue, following the cost structure of \cite{biodiesel} with
unit-level startup and shutdown penalties added,
\begin{equation}\label{eq:bio_obj}
    \min \;\; \int_0^{24} c_e(t)\,P_{\text{tot}}(t)\,dt
      + \sum_{u \in \mathcal{U}} \bigl(
        C_u^{\text{su}} n_u^{\text{su}}
        + C_u^{\text{sd}} n_u^{\text{sd}} \bigr)
      - \int_0^{24} \bigl( \bar{p}_f(t) D_f(t)
        + \bar{p}_g(t) D_g(t) \bigr) dt,
\end{equation}
where $P_{\text{tot}}(t) = P_{\mathrm{rx}}(t)
+ P_{\mathrm{meth}}(t) + P_f(t) + P_g(t)$ is the total plant
power draw, $C_u^{\text{su/sd}}$ are switching costs, and
$n_u^{\text{su/sd}}$ are transition counts computed on the
support grid from successive values of the binary indicators.

The glycerol column power $P_g$ follows the temperature-coupled chain
\eqref{eq:bio_Qreb}--\eqref{eq:bio_Pg}; the remaining power terms are
linear flow-driven surrogates. Each column draws a condenser duty
proportional to its throughput, and the reactor and the FAME and
methanol columns convert their heat duties into electrical power at
constant coefficients of performance,
\begin{subequations}\label{eq:bio_power}
\begin{align}
    & Q_u^{\text{cond}}(t) = c_u\,F_u(t),
      && u \in \{\mathrm{meth}, \mathrm{f}, \mathrm{g}\},
      \label{eq:bio_Qcond} \\
    & P_{\mathrm{rx}}(t) = Q_{\mathrm{rx}}(t)/\eta_{\mathrm{rx}},
      \label{eq:bio_Prx} \\
    & P_u(t) = Q_u^{\text{reb}}(t)/\mathrm{COP}_u^{\text{reb}}
        + Q_u^{\text{cond}}(t)/\mathrm{COP}^{\text{cond}},
      && u \in \{\mathrm{meth}, \mathrm{f}\},
      \label{eq:bio_Pcol}
\end{align}
\end{subequations}
where $\eta_{\mathrm{rx}}$ is the reactor heating efficiency,
$\mathrm{COP}_u^{\text{reb}}$ the reboiler COP of column $u$, and
$\mathrm{COP}^{\text{cond}}$ the shared condenser COP. The methanol
recovery column carries no independent on/off decision: it has no
buffer tank, so its throughput is the crude methanol flow
$F_{\text{meoh}}(t) = \gamma_{\text{meoh}} F_{\text{oil}}(t)$ fixed by
the reactor, and its reboiler duty $Q_{\mathrm{meth}}^{\text{reb}}$
together with the reactor heat duty $Q_{\mathrm{rx}}$ are enforced
within the reactor disjunction \eqref{eq:bio_disj_rx}, so both vanish
when the reactor is off. The glycerol column power \eqref{eq:bio_Pg}
takes its condenser duty $Q_{\mathrm{g}}^{\text{cond}}$ from
\eqref{eq:bio_Qcond}, and the FAME reboiler duty $Q_f^{\text{reb}}$ in
\eqref{eq:bio_Pcol} is the one set by its disjunction
\eqref{eq:bio_disj_f}.

Each unit independently switches on or off through a two-term
disjunction. Linear ON-state surrogates of the unit
models in \cite{biodiesel} are adopted: stoichiometric yields for the
reactor (drawn from the transesterification stoichiometry
reviewed in \cite{ma1999biodiesel}) and linear flow-driven
reboiler duties for the columns. The transesterification
reactor sets crude product flows from oil flow $F_{\text{oil}}$
via the stoichiometric yields $\gamma_c$ for $c \in \mathcal{C}
= \{f, g, \text{meoh}\}$,
\begin{equation}\label{eq:bio_disj_rx}
    \begin{bmatrix}
        W_{\mathrm{rx}}(t) \\[3pt]
        F_{\text{oil}}(t) \geq \underline{F}_{\text{oil}} \\
        F_c(t) = \gamma_c\,F_{\text{oil}}(t),
            \; c \in \mathcal{C} \\
        Q_{\mathrm{rx}}(t) = q_0^{\mathrm{rx}}
            + q_1^{\mathrm{rx}} F_{\text{oil}}(t) \\
        Q_{\mathrm{rx}}(t) \leq Q^{\text{avail}}(t) \\
        Q_{\mathrm{meth}}^{\text{reb}}(t) = q_0^{\mathrm{meth}}
            + q_1^{\mathrm{meth}} F_{\text{meoh}}(t)
    \end{bmatrix} \vee
    \begin{bmatrix}
        \neg W_{\mathrm{rx}}(t) \\[3pt]
        F_{\text{oil}}(t) = 0 \\
        F_c(t) = 0,\; c \in \mathcal{C} \\
        Q_{\mathrm{rx}}(t) = 0 \\
        Q_{\mathrm{meth}}^{\text{reb}}(t) = 0
    \end{bmatrix},
    \quad t \in [0, 24].
\end{equation}
The ON branch caps the electrified reactor duty by a
time-varying utility availability,
\begin{equation}\label{eq:bio_qavail}
    Q^{\text{avail}}(t) = 300 - 80\cos(2\pi t / 24)~\text{kW},
\end{equation}
which is an arbitrary approximation of variable output of the utility supplying the reactor, chosen for the sake of example. The cap stays above $220$~kW, so the ON region
remains feasible at every support, and it makes the reactor
disjunct's constraint data explicitly time-dependent.
The FAME distillation column carries a linear reboiler heat
duty,
\begin{equation}\label{eq:bio_disj_f}
    \begin{bmatrix}
        W_{\mathrm{f}}(t) \\[3pt]
        F_f(t) \geq \underline{F}_f \\
        Q_f^{\text{reb}}(t) = q_0^f + q_1^f\,F_f(t)
    \end{bmatrix} \vee
    \begin{bmatrix}
        \neg W_{\mathrm{f}}(t) \\[3pt]
        F_f(t) = 0 \\
        Q_f^{\text{reb}}(t) = 0
    \end{bmatrix},
    \quad t \in [0, 24].
\end{equation}
The glycerol column carries the temperature-coupled heat duty
\eqref{eq:bio_Qreb}; in the OFF disjunct, the temperature is
held at its lower bound $T_g^L$ so the bilinear COP coupling
\eqref{eq:bio_Preb} is inactive,
\begin{equation}\label{eq:bio_disj_g}
    \begin{bmatrix}
        W_{\mathrm{g}}(t) \\[3pt]
        F_g(t) \geq \underline{F}_g \\
        Q_g^{\text{reb}}(t) = q_0^g + q_1^g F_g(t)
          + \kappa_T (T_g(t) - T_g^L) F_g(t)
    \end{bmatrix} \vee
    \begin{bmatrix}
        \neg W_{\mathrm{g}}(t) \\[3pt]
        F_g(t) = 0,\; Q_g^{\text{reb}}(t) = 0 \\
        T_g(t) = T_g^L,\; P_g^{\text{reb}}(t) = 0
    \end{bmatrix},
    \quad t \in [0, 24].
\end{equation}

Buffer and product tanks decouple the unit schedules \cite{biodiesel},
so the plant can meet downstream demand while turning off
upstream units during high-price intervals,
\begin{equation}\label{eq:bio_tanks}
    \tfrac{dL_u^{\text{buf}}}{dt}
        = F_u^{\text{crude}}(t) - F_u(t),
    \quad
    \tfrac{dL_u^{\text{prod}}}{dt}
        = F_u(t) - D_u(t),
    \quad u \in \{f, g\}.
\end{equation}
A dynamic mass balance tracks the pure glycerol mass in the
product tank from the temperature-dependent purity inflow and
the contract-purity outflow,
\begin{equation}\label{eq:bio_mglycerol}
    \tfrac{dM_g^{\text{pure}}}{dt} =
        \pi_g(t) F_g(t) - \pi^{\text{spec}} D_g(t),
    \quad t \in [0, 24],
\end{equation}
and the blended inventory must meet the contract purity
$\pi^{\text{spec}} = 0.99$ at all times,
\begin{equation}\label{eq:bio_spec}
    M_g^{\text{pure}}(t) \geq \pi^{\text{spec}}\,
        L_g^{\text{prod}}(t),
    \quad t \in [0, 24].
\end{equation}

Tank levels and pure glycerol mass return to their initial
values at the end of the horizon,
\begin{equation}\label{eq:bio_bc}
    L^{\text{buf}}(0) = L^{\text{buf}}(24) = L^0,
    \;\;
    L^{\text{prod}}(0) = L^{\text{prod}}(24) = L^0,
    \;\;
    M_g^{\text{pure}}(0) = M_g^{\text{pure}}(24)
        = \pi^{\text{spec}} L^0,
\end{equation}
and $W_u : [0, 24] \mapsto \{\text{True}, \text{False}\}$ for
$u \in \mathcal{U}$.

The full formulation reads
\begin{equation*}
\begin{aligned}
\min \;\; & \text{electricity cost, switching, revenue}
    && \eqref{eq:bio_obj} \\
\text{s.t.} \;\;
    & \text{temperature couplings}
    && \eqref{eq:bio_purity}\text{--}\eqref{eq:bio_Pg} \\
    & \text{unit power draws}
    && \eqref{eq:bio_power} \\
    & \text{on/off disjunctions}
    && \eqref{eq:bio_disj_rx}\text{--}\eqref{eq:bio_disj_g} \\
    & \text{tank dynamics}
    && \eqref{eq:bio_tanks} \\
    & \text{glycerol mass balance and purity spec}
    && \eqref{eq:bio_mglycerol},\,\eqref{eq:bio_spec} \\
    & \text{boundary conditions}
    && \eqref{eq:bio_bc}.
\end{aligned}
\end{equation*}
The key structural feature is the temperature--COP
coupling within the glycerol disjunction: the OFF disjunct holds
$T_g(t)$ at $T_g^L$, and $T_g(t)$ acts through the bilinear
term inside \eqref{eq:bio_Qreb} and the bilinear COP relation
\eqref{eq:bio_Preb} only when the column is on. Unlike the F1
problem, where a single exactly-one disjunction selects among
four powertrain modes, the biodiesel plant has three
independent binary disjunctions, one per process unit. This
structure is more typical of industrial scheduling problems, in
which the startup/shutdown economics of each unit interact only
through the downstream tank dynamics.

\begin{figure}[H]
\centering
\begin{subfigure}[t]{0.47\textwidth}
\centering
\begin{tikzpicture}
\begin{axis}[trajpanel, name=bioaxA, xlabel={$t$ [h]},
  xmin=0, xmax=24, axis y line*=left,
  ylabel={Electricity price [\$/MWh]}, ymin=0]
  \addplot[gray!55!black] table[col sep=comma, x=t_h,
    y=elec_price] {data/traj_bio_nc.csv};
\end{axis}
\begin{axis}[trajpanel, at=(bioaxA.south east),
  anchor=south east, xmin=0, xmax=24, axis x line=none,
  axis y line*=right, grid=none, ymin=0, ymax=3,
  ytick={0.5,1.5,2.5}, yticklabels={g,f,rx}]
  \addplot[const plot, teal!70!black] table[col sep=comma,
    x=t_h, y expr=\thisrow{Y_glyc}*0.7+0.15]
    {data/traj_bio_nc.csv};
  \addplot[const plot, orange!85!black] table[col sep=comma,
    x=t_h, y expr=\thisrow{Y_fame}*0.7+1.15]
    {data/traj_bio_nc.csv};
  \addplot[const plot, blue!70!black] table[col sep=comma,
    x=t_h, y expr=\thisrow{Y_reactor}*0.7+2.15]
    {data/traj_bio_nc.csv};
\end{axis}
\end{tikzpicture}
\caption{Unit on/off schedules with electricity price}
\label{fig:bio_sched}
\end{subfigure}\hfill
\begin{subfigure}[t]{0.47\textwidth}
\centering
\begin{tikzpicture}
\begin{axis}[trajpanel, name=bioaxB, xlabel={$t$ [h]},
  xmin=0, xmax=24, axis y line*=left,
  ylabel={Reboiler temperature $T_g$ [K]},
  ymin=345, ymax=435,
  legend style={at={(0.5,1.02)}, anchor=south, legend columns=2,
  /tikz/every even column/.append style={column sep=6pt}}]
  \addplot[red!75!black] table[col sep=comma, x=t_h,
    y=T_glyc_reb] {data/traj_bio_nc.csv};
  \addlegendimage{blue!60!black, densely dashed}
  \legend{$T_g(t)$, $\mathrm{COP}_g(t)$}
\end{axis}
\begin{axis}[trajpanel, at=(bioaxB.south east),
  anchor=south east, xmin=0, xmax=24, axis x line=none,
  axis y line*=right, grid=none,
  ylabel={Glycerol heat-pump COP}]
  \addplot[blue!60!black, densely dashed] table[col sep=comma,
    x=t_h, y=COP_glyc] {data/traj_bio_nc.csv};
\end{axis}
\end{tikzpicture}
\caption{Reboiler temperature and heat-pump COP}
\label{fig:bio_temp}
\end{subfigure}

\vspace{0.7em}

\begin{subfigure}[t]{0.62\textwidth}
\centering
\begin{tikzpicture}
\begin{axis}[trajpanel, xlabel={$t$ [h]}, xmin=0, xmax=24,
  ylabel={Inventory level}, legend style={at={(0.5,1.02)},
  anchor=south, legend columns=3,
  /tikz/every even column/.append style={column sep=6pt}}]
  \addplot[blue!70!black] table[col sep=comma, x=t_h,
    y=L_buffer_fame] {data/traj_bio_nc.csv};
  \addplot[blue!40!black, densely dashed] table[col sep=comma,
    x=t_h, y=L_prod_fame] {data/traj_bio_nc.csv};
  \addplot[teal!70!black] table[col sep=comma, x=t_h,
    y=L_buffer_glyc] {data/traj_bio_nc.csv};
  \addplot[teal!45!black, densely dashed] table[col sep=comma,
    x=t_h, y=L_prod_glyc] {data/traj_bio_nc.csv};
  \addplot[orange!85!black, dotted, line width=1.1pt]
    table[col sep=comma, x=t_h, y=M_glyc_pure]
    {data/traj_bio_nc.csv};
  \legend{FAME buffer, FAME product, glyc.\ buffer,
    glyc.\ product, pure glyc.\ mass}
\end{axis}
\end{tikzpicture}
\caption{Tank inventories}
\label{fig:bio_tanks}
\end{subfigure}
\caption{Solution trajectory for the nonconvex biodiesel
scheduling problem ($N = 300$, MBM reformulation, objective
$-213423$ at a $0.14\%$ gap). The trajectory is an illustrative
solve with Gurobi's default settings; the benchmark runs of
Table~\ref{tab:summary} use the configuration stated there.
(a) Unit on/off schedules (rx reactor,
f FAME, g glycerol) against the time-varying electricity price.
(b) Glycerol reboiler temperature $T_g(t)$ and the resulting
heat-pump COP. (c) Buffer and final product inventories with the
purified glycerol mass.}
\label{fig:bio_traj}
\end{figure}
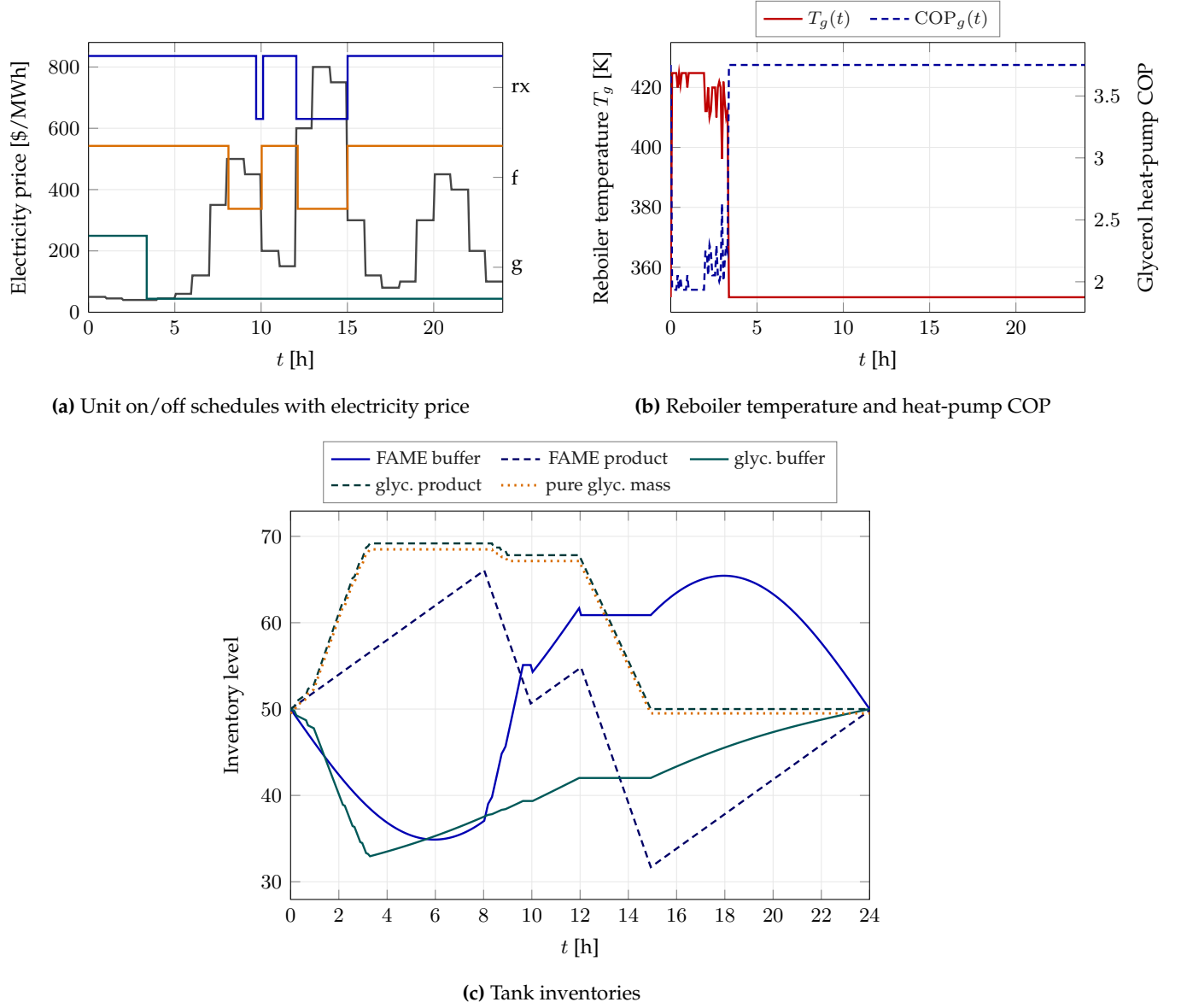

This problem is solved over $[0, 24]$~h, discretized at $N = 300$
support points. Figure~\ref{fig:bio_traj} shows the optimal schedule,
identical across the linearizing methods: the three units shed load
during the electricity
price peaks, the glycerol reboiler temperature tracks the COP trade-off,
and the buffer tanks decouple the unit schedules. Because the only
nonconvexities are the two glycerol bilinear terms, the guidance of
Section~\ref{sec:guidance} predicts that all four linearizing methods
reach the global optimum, differing only in the split between
reformulation and solver time, and that CP -- nonconvex here --
generates invalid cuts. LOA also loses its convexity guarantee on
this problem. It is still run, with the expectation that its
augmented-penalty form tolerates the two bilinear terms and returns a
good local schedule without a valid bound.

\begin{figure}[H]
\centering
\begin{tikzpicture}
\begin{axis}[timebar, symbolic x coords={BigM,MBM,MBM-GP,Hull,PSplit-3},
             area legend,
             legend style={at={(0.5,1.03)}, anchor=south,
               legend columns=-1, draw=black!50,
               font=\scriptsize}]
  \addplot+[fill=blue!25, draw=blue!55!black]
    table[col sep=comma, x=method, y=reform_s]
    {data/time_biodiesel_n300_reform.csv};
  \addplot+[fill=orange!35, draw=orange!70!black]
    table[col sep=comma, x=method, y=solver_s]
    {data/time_biodiesel_n300_reform.csv};
  \legend{Reformulation, Solver}
\end{axis}
\end{tikzpicture}
\caption{Reformulation versus solver time on the nonconvex biodiesel
problem ($N = 300$); axes as in Figure~\ref{fig:timebar_f1}.}
\label{fig:timebar_bio}
\end{figure}
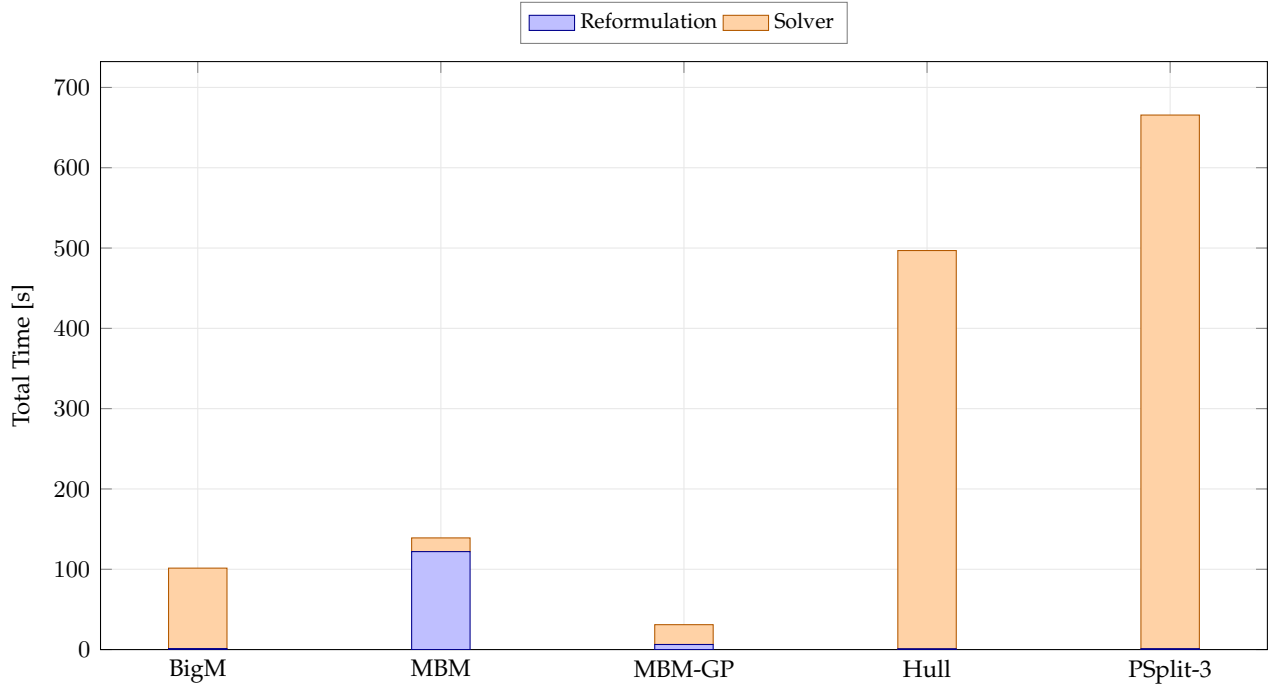

\paragraph{Results.} Big-M uses $M = 10^{3}$ and P-split uses
$P = 3$. The bilinear couplings make the transcribed model a nonconvex
MIQCP, solved by spatial branch-and-bound. All five reformulations
solve the instance to the gap target inside the $1200$~s budget,
and the availability cap \eqref{eq:bio_qavail} makes the reactor
disjunct's $M(t)$ curves genuinely time-varying.
Figure~\ref{fig:solpct_bio} shows the resulting separation: MBM
closes in $17$~s of solver time, big-M needs $100$~s ($0.8\%$ gap
at $-213308$), hull runs $496$~s but lands on the best known
solution $-213682$ with an essentially closed gap ($0.006\%$), and
P-split-3 takes $665$~s to a $1.3\%$ incumbent. P-split lifts the
bilinear terms to McCormick-enveloped auxiliary variables, since
its plain partitioning cannot represent them. The root relaxations
reflect this ordering: MBM and hull start within $0.03\%$ of the
optimum, while big-M and P-split both start with a $2.1\%$ root gap.
The automatic partition leaves a wide-range variable in each group,
so the per-partition envelopes are no tighter than a single big-M,
while the lifted model carries $3\times$ hull's rows, and this
larger model must be processed at every branch-and-bound node. CP is omitted, as its
separation subproblem inherits the nonconvexity. LOA finds no
feasible incumbent within the budget: its NLP subproblems,
nonconvex QCPs at fixed mode schedules, never return a solution, so
it fails here exactly as on F1 and is likewise absent from the
figure. Figure~\ref{fig:timebar_bio} shows the cost split. MBM's
$122$~s reformulation bar reflects its per-pair $M$ computation,
now itself a set of nonconvex subproblems solved at every support,
but the resulting tight master yields the fastest solver time. MBM-GP
keeps most of that tightness at a fraction of the build. It matches
grid MBM's objective ($-213397$ for both) while cutting the
reformulation time from $122$ to $6.3$~s, a $19\times$ reduction.
In this case, the saving comes from the GP fit itself: the
surrogate samples the time-varying $M(t)$ induced by \eqref{eq:bio_qavail} instead of solving all
$\approx 16000$ per-support subproblems. Its $31$~s total is the
fastest of any method, ahead of grid MBM's $139$~s and big-M's
$101$~s.

\begin{figure}[H]
\centering
\begin{tikzpicture}
\begin{axis}[solpanel, xmin=0, xmax=700]
  \addplot[const plot, mark=none, black] table
    [col sep=comma, x=time_s, y=sol_frac]
    {data/solpct_bio_nc_bigm.csv};
  \addplot[const plot, mark=none, blue!70!black, dashed] table
    [col sep=comma, x=time_s, y=sol_frac]
    {data/solpct_bio_nc_mbm.csv};
  \addplot[const plot, mark=none, green!60!black, densely dotted] table
    [col sep=comma, x=time_s, y=sol_frac]
    {data/solpct_bio_nc_mbmgp.csv};
  \addplot[const plot, mark=none, orange!85!black] table
    [col sep=comma, x=time_s, y=sol_frac]
    {data/solpct_bio_nc_hull.csv};
  \addplot[const plot, mark=none, teal!70!black] table
    [col sep=comma, x=time_s, y=sol_frac]
    {data/solpct_bio_nc_psplit.csv};
  \legend{big-M, MBM, MBM-GP, hull, P-split}
  \draw[gray!55, dashed, line width=0.7pt]
    (axis cs:0,0.95) -- (axis cs:700,0.95);
  \node[anchor=north east, font=\scriptsize, gray!60!black]
    at (axis cs:690,0.945) {5\% MIP-Gap Target};
\end{axis}
\end{tikzpicture}
\caption{Solution progress versus solve time on the nonconvex
biodiesel problem ($N = 300$); axes as in
Figure~\ref{fig:solpct_f1}. LOA is absent: it finds no feasible
incumbent within the budget (Table~\ref{tab:summary}).}
\label{fig:solpct_bio}
\end{figure}
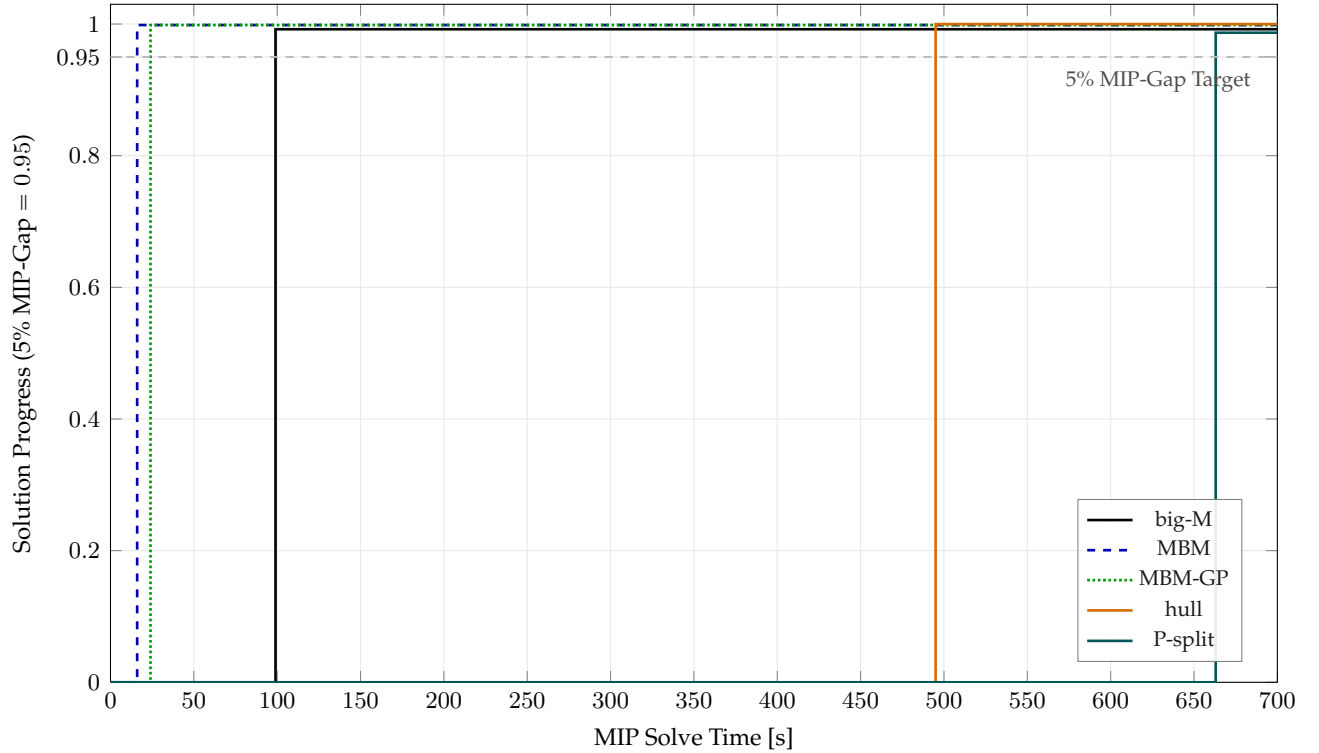

CP fails here exactly as Section~\ref{sec:guidance} anticipates. To
verify the failure directly, a fresh big-M model was solved to a
verified optimum and each cut CP had generated was evaluated at that
optimum. A valid cut must be nonnegative at every feasible point, in
particular at the true optimum; a strictly negative value certifies
that the cut excludes it. All three generated cuts were strictly
violated at the verified big-M optimum, with values $-3.1 \times
10^{7}$, $-8.6 \times 10^{4}$, and $-6.2 \times 10^{4}$ (separation
solved locally with Ipopt, as required because the separation problem
inherits the model's nonconvexity). A global solver would arrive at a
globally optimal separation with enough time, but the resulting cut
could still be invalid. The
consequence compounds with
size: at $N = 300$ the CP run returns a solver error with no
recoverable solution, and where an incumbent is recovered its objective
is $\approx 25\%$ worse than the true optimum. The same failure occurs
on the nonconvex F1 problem, where CP produces no valid incumbent at
all. CP is therefore reported as a benchmarked method only on the
convex powergrid case study below.

\subsection{Event-constrained IEEE 14-bus capacity design}
\label{sec:powergrid}

The third case study is a chance-constrained transmission and
generation capacity design problem on the IEEE 14-bus test
system, adapted verbatim from the event-constrained programming
formulation of \cite{ovalle2025event}. A planner
sizes oversizing margins on five generators and twenty lines so
that, under random nodal demand $\xi \in \mathbb{R}^{11}$, the
event that no component exceeds its rated capacity occurs with
probability at least $\alpha$. Unlike the F1 and biodiesel
problems, the infinite parameter here is a probability
distribution rather than a physical domain; the support grid is
a Monte Carlo sample $\{\hat{\xi}_k\}_{k=1}^N$ drawn from
$\xi \sim \mathcal{N}(\xi^{\mathrm{nom}}, \Sigma)$. All
constraints are linear, so the transcribed reformulations are
MILPs and the cutting plane separation subproblem is convex,
making CP theoretically sound and providing the control case
against which CP's invalid-cut behavior on the two nonconvex
problems is interpreted.

Two design vectors are finite: generator
oversizing $z^{\mathrm{gen}} \in [0, 300]^5$ and line oversizing
$z^{\mathrm{line}} \in [0, 100]^{20}$. Two state vectors are
infinite in $\xi$: generator dispatch
$y^{\mathrm{gen}}(\xi) \in [0, 1000]^5$ and line flow
$y^{\mathrm{line}}(\xi) \in [-1000, 1000]^{20}$.

The objective is the total oversizing cost,
\begin{equation}\label{eq:pg_obj}
    \min_{y,\, z,\, W} \;\;
    \sum_{g=1}^{5} z^{\mathrm{gen}}_g
    + \sum_{\ell=1}^{20} z^{\mathrm{line}}_\ell .
\end{equation}

At each of the fourteen buses, net
injection equals demand. The system is written compactly with the
bus-incidence pattern of \cite{ovalle2025event} as
\begin{equation}\label{eq:pg_balance}
    A^{\mathrm{gen}}\, y^{\mathrm{gen}}(\xi)
    \;+\; A^{\mathrm{line}}\, y^{\mathrm{line}}(\xi)
    \;=\; B\, \xi,
    \qquad \xi \in \Xi,
\end{equation}
where $A^{\mathrm{gen}} \in \{0, 1\}^{14 \times 5}$ and
$A^{\mathrm{line}} \in \{-1, 0, 1\}^{14 \times 20}$ encode which
generator and line is incident to each bus, and
$B \in \{0, 1\}^{14 \times 11}$ maps the eleven uncertain demands
onto their buses.

For each generator
$g$, dispatch is either within or above its design limit
$\bar{p}_g + z^{\mathrm{gen}}_g$ (with
$\bar{p} = (332, 140, 100, 100, 100)$~MW),
\begin{equation}\label{eq:pg_disj_gen}
    \begin{bmatrix}
        W^{\mathrm{gen}}_g(\xi) \\[3pt]
        y^{\mathrm{gen}}_g(\xi)
            \leq \bar{p}_g + z^{\mathrm{gen}}_g
    \end{bmatrix}
    \vee
    \begin{bmatrix}
        \neg W^{\mathrm{gen}}_g(\xi) \\[3pt]
        y^{\mathrm{gen}}_g(\xi)
            \geq \bar{p}_g + z^{\mathrm{gen}}_g
    \end{bmatrix},
    \quad g = 1, \ldots, 5,\; \xi \in \Xi .
\end{equation}

For each line $\ell$,
flow is either within the symmetric band
$[-(\bar{f} + z^{\mathrm{line}}_\ell),\,
\bar{f} + z^{\mathrm{line}}_\ell]$ around zero (with base limit
$\bar{f} = 50$~MW) or it exceeds the band in one of the two
directions,
\begin{equation}\label{eq:pg_disj_line}
    \begin{bmatrix}
        W^{\mathrm{line, low}}_\ell(\xi) \\[3pt]
        y^{\mathrm{line}}_\ell(\xi)
            \leq -(\bar{f} + z^{\mathrm{line}}_\ell)
    \end{bmatrix}
    \vee
    \begin{bmatrix}
        W^{\mathrm{line}}_\ell(\xi) \\[3pt]
        |y^{\mathrm{line}}_\ell(\xi)|
            \leq \bar{f} + z^{\mathrm{line}}_\ell
    \end{bmatrix}
    \vee
    \begin{bmatrix}
        W^{\mathrm{line, high}}_\ell(\xi) \\[3pt]
        y^{\mathrm{line}}_\ell(\xi)
            \geq \bar{f} + z^{\mathrm{line}}_\ell
    \end{bmatrix},
\end{equation}
for $\ell = 1, \ldots, 20$ and $\xi \in \Xi$.

The event being
constrained is that at least $n^{\mathrm{gen}}_{\min}$ generators
and at least $n^{\mathrm{line}}_{\min}$ lines are within their
design limits. Following \cite{ovalle2025event}, two aggregator
indicators $W^{\mathrm{at}}_{\mathrm{gen}}(\xi)$ and
$W^{\mathrm{at}}_{\mathrm{line}}(\xi)$ are introduced, tied to the
counts of within-limit components,
\begin{align}
    \sum_{g} w^{\mathrm{gen}}_g(\xi)
        \geq n^{\mathrm{gen}}_{\min}
        &\;\;\Longleftrightarrow\;\;
        W^{\mathrm{at}}_{\mathrm{gen}}(\xi),
        \label{eq:pg_atleast_gen} \\
    \sum_{\ell} w^{\mathrm{line}}_\ell(\xi)
        \geq n^{\mathrm{line}}_{\min}
        &\;\;\Longleftrightarrow\;\;
        W^{\mathrm{at}}_{\mathrm{line}}(\xi),
        \label{eq:pg_atleast_line}
\end{align}
where $w^{\mathrm{gen}}_g, w^{\mathrm{line}}_\ell$ are the binary
realizations of the disjunct indicators in
\eqref{eq:pg_disj_gen}--\eqref{eq:pg_disj_line}. The event
indicator $W(\xi)$ is their conjunction,
\begin{equation}\label{eq:pg_event}
    \bigl( W^{\mathrm{at}}_{\mathrm{gen}}(\xi)
        \wedge W^{\mathrm{at}}_{\mathrm{line}}(\xi) \bigr)
    \;\;\Longleftrightarrow\;\; W(\xi) .
\end{equation}

The expectation of
$W(\xi)$ under the Monte Carlo sample equals the empirical event
probability, and is required to be at least $\alpha$,
\begin{equation}\label{eq:pg_prob}
    \mathbb{E}_\xi[\, W(\xi) \,]
    \;=\; \frac{1}{N} \sum_{k=1}^{N} w_k
    \;\geq\; \alpha .
\end{equation}

The parameter settings of
\cite{ovalle2025event} are used: $\alpha = 0.9$,
$n^{\mathrm{gen}}_{\min} = 5$, $n^{\mathrm{line}}_{\min} = 20$
(so the atleast logic of
\eqref{eq:pg_atleast_gen}--\eqref{eq:pg_atleast_line} reduces in
this configuration to the conjunction over all components), and
$N = 200$ Monte Carlo samples drawn with seed 42 from
$\xi \sim \mathcal{N}(\xi^{\mathrm{nom}}, \Sigma)$ with
$\xi^{\mathrm{nom}}_{1:6} = (87.3, 50, 25, 28.8, 50, 25)$~MW and
the remaining components zero. The model is taken verbatim from
\texttt{power\_grid.jl} of \cite{ovalle2025event}, with the only
inert change being explicit upper bounds on the count auxiliaries
in \eqref{eq:pg_atleast_gen}--\eqref{eq:pg_atleast_line} that
hull requires and that act as a no-op for big-M.

Because every constraint is linear in $(y, z, W)$, the
transcribed reformulations are MILPs and CP's hull-based
separation subproblem is convex. Per Section~\ref{sec:guidance}, this
is the regime in which CP and LOA are sound, so the study plays a dual
role: it benchmarks the linearizing methods on a discretely-uncertain,
non-dynamic InfiniteGDP, and it supplies the convex control case in
which CP's cuts are valid and LOA's OA bound is legitimate.

\begin{figure}[H]
\centering
\begin{tikzpicture}
\begin{axis}[solpanel, xmin=0, xmax=1250]
  \addplot[const plot, mark=none, black] table
    [col sep=comma, x=time_s, y=sol_frac]
    {data/solpct_pg_bigm.csv};
  \addplot[const plot, mark=none, blue!70!black, dashed] table
    [col sep=comma, x=time_s, y=sol_frac]
    {data/solpct_pg_mbm.csv};
  \addplot[const plot, mark=none, green!60!black, densely dotted] table
    [col sep=comma, x=time_s, y=sol_frac]
    {data/solpct_pg_mbmgp.csv};
  \addplot[const plot, mark=none, orange!85!black] table
    [col sep=comma, x=time_s, y=sol_frac]
    {data/solpct_pg_hull.csv};
  \addplot[const plot, mark=none, teal!70!black] table
    [col sep=comma, x=time_s, y=sol_frac]
    {data/solpct_pg_psplit2.csv};
  \addplot[const plot, mark=none, red!80!black, densely
    dashed] table
    [col sep=comma, x=time_s, y=sol_frac]
    {data/solpct_pg_cp_finalmip.csv};
  \addplot[const plot, mark=none, violet!80!black, dashdotted] table
    [col sep=comma, x=time_s, y=sol_frac]
    {data/solpct_pg_loa.csv};
  \legend{big-M, MBM, MBM-GP, hull, P-split, CP, LOA$^{*}$}
  \draw[gray!55, dashed, line width=0.7pt]
    (axis cs:0,0.95) -- (axis cs:1250,0.95);
  \node[anchor=north east, font=\scriptsize, gray!60!black]
    at (axis cs:1230,0.945) {5\% MIP-Gap Target};
\end{axis}
\end{tikzpicture}
\caption{Solution progress versus solve time on the convex powergrid
case study ($N = 200$); axes as in Figure~\ref{fig:solpct_f1}. The CP
curve is the branch-and-bound trace of its final big-M+cuts MIP; cut
generation occurs in the reformulation phase
(Figure~\ref{fig:timebar_powergrid}). LOA$^{*}$ plots the relative OA
gap of its iteration loop; because every constraint is linear, the OA
master carries the full problem exactly and the optimum is reached
trivially once the set-covering initialization completes, so the
trace measures initialization rather than search.}
\label{fig:solpct_powergrid}
\end{figure}
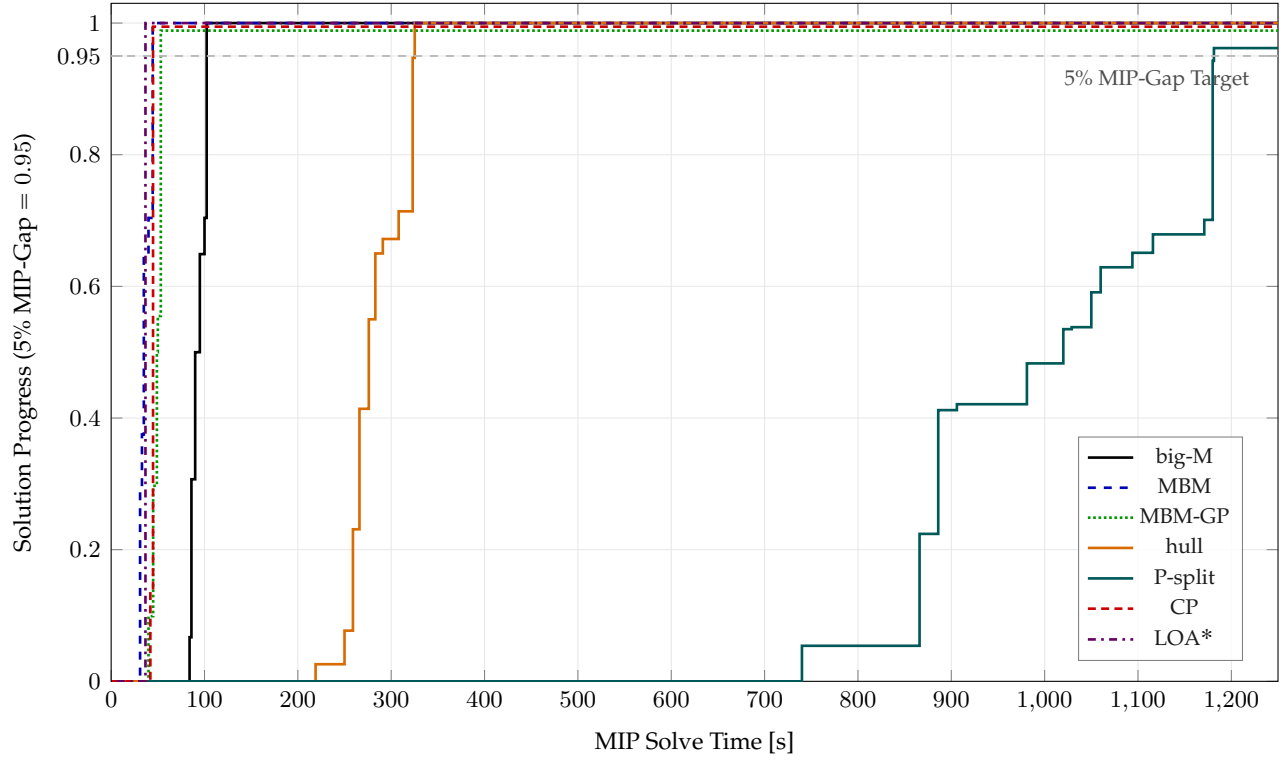

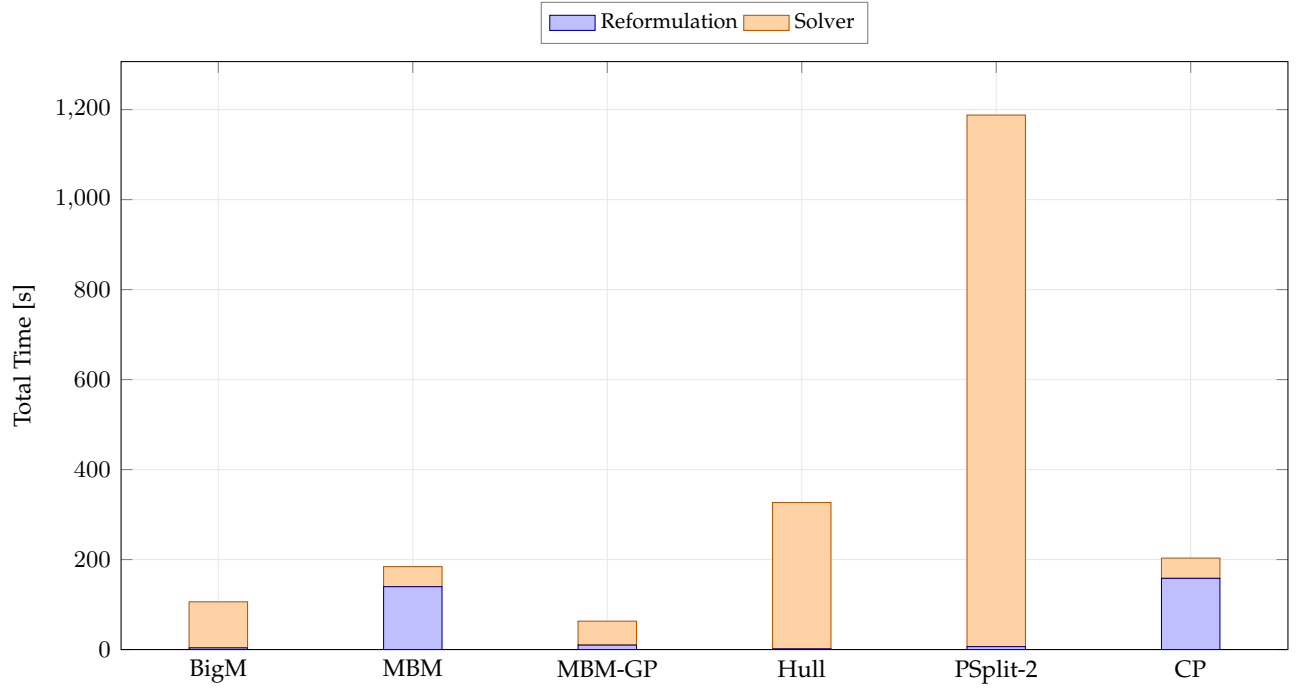
\begin{figure}[H]
\centering
\begin{tikzpicture}
\begin{axis}[timebar,
             symbolic x coords={BigM,MBM,MBM-GP,Hull,PSplit-2,CP},
             area legend,
             legend style={at={(0.5,1.03)}, anchor=south,
               legend columns=-1, draw=black!50,
               font=\scriptsize}]
  \addplot+[fill=blue!25, draw=blue!55!black]
    table[col sep=comma, x=method, y=reform_s]
    {data/time_powergrid_n200_reform.csv};
  \addplot+[fill=orange!35, draw=orange!70!black]
    table[col sep=comma, x=method, y=solver_s]
    {data/time_powergrid_n200_reform.csv};
  \legend{Reformulation, Solver}
\end{axis}
\end{tikzpicture}
\caption{Reformulation versus solver time on the convex powergrid case
study ($N = 200$, $5\%$ MIP gap, Gurobi 12.0.2); axes as in
Figure~\ref{fig:timebar_f1}. Timings are single-trial.}
\label{fig:timebar_powergrid}
\end{figure}

\paragraph{Results.} Big-M uses $M = 10^{5}$, P-split uses
$P = 2$, and CP runs to natural convergence with a separation
tolerance of $10^{-6}$. As
Section~\ref{sec:guidance} predicts for a convex GDP, all methods
--- now including MBM-GP, CP and LOA --- reach the same global
optimum ($36.71$). P-split closes last, needing nearly the full
budget: the root LP bound is zero for every reformulation on this
event-constrained structure, and P-split's lifted transcription
carries $8\times$ big-M's rows, so its bound moves slowest.
Figure~\ref{fig:solpct_powergrid} shows the progress
curves. LOA's curve reflects a
different mechanism than the branch-and-bound traces. Because every
constraint and the objective are linear, the disjunct constraints
enter the master exactly via big-M and the objective and global OA
cuts are exact, so the OA master reduces to a big-M reformulation of
the transcribed model, augmented by the OA and no-good cuts
accumulated during the set-covering seed. LOA thus recovers
essentially the same big-M MILP structure the standalone big-M method
solves, and the optimum is available once the seed identifies the
optimal trajectory and the master bound matches it, with no further
search required.
Complementing this,
Figure~\ref{fig:timebar_powergrid} reports the reformulation/solver
split, with CP's $158$~s reformulation bar accounting for cut
generation. Hull solves in $325$~s and big-M in $102$~s, while
P-split's $1181$~s shows its lifted master paying its size's full
price. MBM-GP again removes most of MBM's
reformulation cost: it builds in $10.0$~s instead of $139.8$~s, solves in a similar amount of time ($53.2$ versus $44.5$~s), and
its $\approx 63$~s total is the fastest of any method on this
study while retaining MBM's per-pair tightness. On this study
every $M(d)$ curve is constant in $d$, so the speedup
comes from the sampler's uniform-$M$ detection rather than the GP
fill. The contrast with the nonconvex F1 and biodiesel results,
where CP produced invalid cuts, confirms that CP's validity is
contingent on convexity of the separation subproblem.

\subsection{Summary}\label{sec:summary}

Table~\ref{tab:summary} collects the objectives, gaps, and the
reformulation/solver time split across the three studies. The
results are consistent with the guidance of
Section~\ref{sec:guidance}. On the
challenging nonconvex F1 problem, whose boost disjunct carries the
position-dependent deployment cap \eqref{eq:f1_depcap}, only the
multiple big-M methods solve within the budget: MBM in $102$~s and
MBM-GP in $194$~s with an 8-fold reduced pre-solve time, while big-M and
hull time out with no incumbent and P-split reaches only a
$16.9\%$-gap lap at the wall. Feasible solutions for this problem
thus come only from the intermediate reformulations rather than
the big-M and hull endpoints. On the biodiesel and powergrid
studies every reformulation closes the gap, and the solver times
spread by more than an order of magnitude in the relaxation
ordering: the MBM family solves fastest ($17$--$53$~s), big-M pays
for its loose relaxation ($100$~s on both), hull pays for its size
($496$ and $325$~s) while returning the tightest certificates,
including the best known biodiesel solution $-213682$ at a
$0.006\%$ gap, and P-split's oversized lifts make it slowest on
both ($665$ and $1181$~s). Across all three studies, MBM-GP reproduces grid MBM's
solution quality to within the gap target while cutting the
reformulation cost by $8$--$19$ times, making it the fastest total
method on biodiesel ($31$~s) and powergrid ($63$~s); on the two
studies where the $M$ curves genuinely vary ($M(s)$ on F1, $M(t)$
on biodiesel) the saving comes from the GP surrogate itself rather
than uniform-$M$ detection. CP matches every method on the convex
powergrid study but is invalid on both nonconvex problems. For the
powergrid problem, LOA converges to the same $36.71$ optimum as
every reformulation, closing its OA gap to numerical precision in
$36.8$~s, faster than any reformulation solves. On the nonconvex
biodiesel and F1 studies its fixed-schedule NLP subproblems are
themselves near-intractable under the benchmark settings: LOA
exhausts its wall budget without an incumbent on biodiesel, and on
F1 its best incumbent ($164.29$~s) is $78\%$ above the best lap
time found. Like CP, LOA is reliable in the convex regime it
targets, but degrades sharply outside of it.

\begin{table}[H]
\centering
\caption{Cross-method summary. Reformulation and solver time are
reported separately; all methods on all three studies use a $5\%$
MIP gap target under a uniform $1200$~s wall budget, with Gurobi's
presolve, cuts, and heuristics disabled so solver time reflects
reformulation quality. Objectives
are lap time [s] for F1, cost [\$] for biodiesel, and total
expansion margin [MW] for powergrid; gap is the solver's terminal
relative MIP gap. LOA is an iterative decomposition, so the
reformulation/solver split does not apply; its Solver column
reports the time at which it reaches its result, and on F1 its Gap
column reports its incumbent's distance from the best-known
objective. A dash marks entries that are not available.}
\label{tab:summary}
\small
\begin{tabular}{llrrrrr}
\toprule
Problem & Method & Status & Objective & Gap [\%]
  & Reform [s] & Solver [s] \\
\midrule
\multirow{7}{*}{F1, $N{=}100$}
 & big-M    & Time Limit & ---      & ---    & $4.5$
   & $1200.2$ \\
 & MBM      & Optimal    & $93.85$  & $3.7$  & $21.7$
   & $102.2$ \\
 & MBM-GP   & Optimal    & $92.24$  & $2.2$  & $2.7$
   & $193.8$ \\
 & hull     & Time Limit & ---      & ---    & $0.7$
   & $1200.2$ \\
 & P-split  & Time Limit & $97.46$  & $16.9$ & $1.2$
   & $1200.6$ \\
 & CP       & ---        & ---      & ---    & ---
   & --- \\
 & LOA      & Time Limit & $164.29$ & $78.1$ & ---
   & $1210.9$ \\
\midrule
\multirow{7}{*}{\shortstack[l]{Biodiesel\\(nonconvex),\\$N{=}300$}}
 & big-M    & Optimal    & $-213308$ & $0.8$ & $1.2$
   & $100.2$ \\
 & MBM      & Optimal    & $-213397$ & $0.2$ & $122.0$
   & $17.1$ \\
 & MBM-GP   & Optimal    & $-213397$ & $0.2$ & $6.3$
   & $24.7$ \\
 & hull     & Optimal    & $-213682$ & $0.0$ & $1.0$
   & $495.9$ \\
 & P-split  & Optimal    & $-213001$ & $1.3$ & $1.1$
   & $664.5$ \\
 & CP       & ---        & ---       & ---   & ---
   & --- \\
 & LOA      & Time Limit & ---       & ---   & ---
   & $1215.8$ \\
\midrule
\multirow{7}{*}{\shortstack[l]{Powergrid,\\$N{=}200$}}
 & big-M    & Optimal    & $36.71$ & $0.0$ & $3.8$
   & $102.3$ \\
 & MBM      & Optimal    & $36.71$ & $0.0$ & $139.8$
   & $44.5$ \\
 & MBM-GP   & Optimal    & $36.71$ & $1.1$ & $10.0$
   & $53.2$ \\
 & hull     & Optimal    & $36.71$ & $0.0$ & $1.6$
   & $325.2$ \\
 & P-split  & Optimal    & $36.71$ & $3.8$ & $6.6$
   & $1181.3$ \\
 & CP       & Optimal    & $36.71$ & $0.5$ & $158.3$
   & $45.0$ \\
 & LOA      & Converged  & $36.71$ & $0.0$ & ---
   & $36.8$ \\
\bottomrule
\end{tabular}
\end{table}

\section{Conclusions and Future Work}\label{sec:conclusions}

We have generalized four GDP solution methods, the multiple big-M,
P-split, and cutting plane reformulations and the logic-based outer
approximation algorithm, to the infinite-dimensional setting. The
reformulations fill the relaxation spectrum between the big-M and
hull endpoints that were previously extended for InfiniteGDP. LOA adds an
iterative alternative that never forms a single mixed-integer
program.
We have also introduced MBM-GP, a Gaussian-process variant of
multiple big-M with no finite counterpart that learns the big-M
function over the infinite domain from a small subset of the
subproblem solves. The methods are implemented in
\texttt{DisjunctiveProgramming.jl}, its extension
\texttt{InfiniteDisjunctiveProgramming.jl}, and the
\texttt{DisjunctiveAlgorithms.jl} optimizer. The benchmark studies of this work demonstrate how the InfiniteGDP MBM method often provides competitive solution times which can be attributed to its tight continuous relaxation; moreover, MBM-GP significantly reduces the upfront pre-solve computational cost when it becomes a bottleneck. Hull and P-split suit problems where
a stronger relaxation justifies a larger model, big-M remains a fast
baseline, and CP and LOA are reliable only on the convex problems
they target. Future work will focus on a global logic-based outer approximation
that rigorously handles nonconvex disjuncts and on GPU-accelerated
solves that exploit the recurrent structure of the transcribed
problems.

\section*{Declaration of generative AI and AI-assisted technologies in the manuscript preparation process}

Authors were assisted by Claude models (Anthropic) in the preparation of this work with figure generation, equation formatting, and grammatical editing. All aspects of this work have been carefully reviewed and edited by the authors such that they take full responsibility for the final manuscript.

\section*{Acknowledgments}

We acknowledge support from the Natural Sciences and Engineering Research Council of Canada under grant RGPIN-202403997 and support provided from the Department of Chemical Engineering at the University of Waterloo.

\bibliography{refs}
\end{document}